\documentclass{arXiv}

\usepackage{amsthm}
\usepackage{todonotes}

\DeclareMathOperator*{\argmin}{arg\,min}

\newcommand{\mat}[3]{#1 \vert_{#2}^{#3}}
\newcommand{\R}{\mathbb{R}}
\newcommand{\T}{\mathbf{T}}
\newcommand{\PSI}{\mathbf{\Psi}}
\newcommand{\PHI}{\mathbf{\Phi}}
\newcommand{\I}{\mathbf{I}}
\newcommand{\J}{\mathbf{J}}

\newtheorem{Theorem}{Theorem}

\newcommand{\innerprod}[2]{\left\langle #1,\, #2 \right\rangle} 
\begin{document}

\title{\huge Tensor-based computation of \newline metastable and coherent sets}

\abstract{Recent years have seen rapid advances in the data-driven analysis of dynamical systems based on Koopman operator theory and related approaches. On the other hand, low-rank tensor product approximations -- in particular the tensor train (TT) format -- have become a valuable tool for the solution of large-scale problems in a number of fields. In this work, we combine Koopman-based models and the TT format, enabling their application to high-dimensional problems in conjunction with a rich set of basis functions or features. We derive efficient algorithms to obtain a reduced matrix representation of the system's evolution operator starting from an appropriate low-rank representation of the data. These algorithms can be applied to both stationary and non-stationary systems. We establish the infinite-data limit of these matrix representations, and demonstrate our methods' capabilities using several benchmark data sets.}


\doi{}
\author{Feliks N\"uske}{Center for Theoretical\\Biological Physics \&\\ Department of Chemistry,\\Rice University,\\Houston, TX, 77005,\\United States  \\ \\ Institute of Mathematics, \\ Universität Paderborn, \\Paderborn 33100,\\Germany }
\author{Patrick Gel\ss}{Department of\\Mathematics and\\Computer Science,\\Freie Universität Berlin,\\Berlin 14195, Germany}
\author{Stefan Klus}{Department of Mathematics,\\University of Surrey,\\
Guildford, GU2 7XH, UK}
\author{Cecilia Clementi}{Center for Theoretical\\Biological Physics \&\\ Department of Chemistry,\\Rice University,\\Houston, TX, 77005,\\United States \\ \\ Department of Physics,\\Freie Universität Berlin,\\Berlin 14195, Germany}
\email{feliks.nueske@uni-paderborn.de}

\maketitle

\section{Introduction}

The data-driven analysis of high-dimensional dynamical systems has been a highly successful research field for several years, with applications in fluid dynamics, control theory, molecular dynamics, and many others. Much of the work along these lines has focused on the infinite-dimensional description of a system using \textit{transfer operators} or \textit{Koopman operators}, see~\cite{Koopman1931,DELLNITZ1999,SCHUETTE1999,Mezic2005,KLUS2016b}. For non-stationary or time-dependent systems, it is often advantageous to consider the \textit{forward-backward operator}. Its spectral analysis leads to the detection of \emph{coherent sets} \cite{FROYLAND2010,FROYLAND2013}. We will summarily refer to all of these operators as \emph{evolution operators} in this paper. A host of different methods for the numerical approximation of evolution operators from simulation or measurement data have been developed. These include \emph{extended dynamic mode decomposition} (EDMD) \cite{WILLIAMS2015, KLUS2016b}, the \emph{variational approach to conformational dynamics} (VAC) \cite{NOE2013,NUESKE2014} and its generalization \emph{variational approach to Markov processes} (VAMP) \cite{WU2020}, as well as \emph{canonical correlation analysis} (CCA) \cite{HARDOON2004,KHM19}. All of these methods are closely related, and revolve, in one way or another, around data-based approximations to the Koopman operator on (usually finite-dimensional) subspaces. For a detailed review and comparison, see \cite{KLUS2018b}.

Much of the appeal of these techniques is due to their formulation as data-driven regression problems, which paves the way for the application of modern machine learning techniques. Examples include kernel-based formulations \cite{WILLIAMS2015b,KHM19,Klus2020a,Klus2020} and combinations with deep learning \cite{Lusch2018,MARDT2018}. A different avenue towards the solution of high-dimensional problems are tensor products, where functions and operators on high-dimensional spaces are approximated in linear spaces of products of simple (often univariate) functions. The expansion coefficients of such a function form a multi-dimensional array, called a \textit{tensor}. As the size of a tensor grows exponentially with the dimension, \textit{low-rank formats} requiring only a manageable number of parameters need to be used. Important examples include the \emph{canonical format} \cite{CARROLL1970}, the \emph{Tucker format} \cite{TUCKER1964}, and the \emph{hierarchical Tucker format} \cite{HACKBUSCH2009}, with the \emph{tensor train (TT) format}~\cite{OSELEDETS2009, OSELEDETS2011} as an important special case of the latter. The common idea behind these formats is to decompose a high-dimensional tensor into a network of lower-dimensional tensors. Several applications of tensor decompositions have shown that it is possible to tackle large-scale problems which cannot be solved using conventional numerical methods, see, e.g., \cite{BECK2000, DOLGOV2015, GELSS2016, ZHANG2017, GELSS2019}, and especially \cite{AFFLECK1987,OSTLUND1995,SZALAY2015} for quantum chemistry applications.

Previous work on low-rank approximation in the context of Koopman operator modeling includes \cite{KLUS2016,NUESKE2016,KLUS2018,LITZINGER2018,GELSS2019}. Importantly, it was shown in~\cite{KLUS2016} that the matrices required for data-driven approximations to the Perron--Frobenius or Koopman operator based on a product basis can directly be written in the canonical format. The resulting generalized eigenvalue problems can be converted to TT format and then solved with the aid of power iteration methods. This, however, requires repeated rank reductions and appropriate estimates of the eigenvalues so that power iteration converges. In \cite{KLUS2018,GELSS2019}, a TT representation of the \emph{data tensor}, which contains the evaluations of a product basis at all data points, was introduced, accompanied by a method to compute an approximate singular value decomposition (global SVD) of the data tensor.

In this study, we build on these results to derive and analyze new tensor decomposition methods for Koopman operator approximation. Our methods proceed by first computing a (compressed) representation of the data tensor. Subsequently, a tensor network needs to be contracted to arrive at a reduced matrix approximation for the Koopman operator. Importantly, the calculation or inversion of Gramian matrices is entirely avoided. The detailed contributions of our study are as follows: 

\begin{itemize}[leftmargin=*]
\item  First, we present a multi-linear version of the \emph{AMUSE} algorithm~\cite{TONG1990} for Koopman approximation, which we call \emph{AMUSEt}. This method uses the global SVD to compute an approximate multi-linear singular value decomposition of the data tensor. By contraction of an appropriate tensor network, a matrix representation of the Koopman operator can then be obtained.
 \item Second, we analyze AMUSEt by showing that the final matrix approximation is indeed the data-driven representation of the Koopman operator on a data-dependent, finite-dimensional subspace. We establish convergence of this representation in the limit of infinite data.
 \item Third, we present an alternative tensor train decomposition of the data tensor corresponding to a product basis. It is based on a method outlined in \cite{OSELEDETS2010}, which provides a \textit{higher-order CUR decomposition}. A detailed description of the latter method's algorithmic realization, including several enhancements, is provided.
 \item Finally, we demonstrate the capabilities of the proposed methods by analyzing benchmark molecular dynamics and fluid dynamics data sets.
\end{itemize}
The remainder of this work is structured as follows: In Section~\ref{sec:basic_concepts}, we introduce the required notation and concepts regarding evolution operators, their numerical approximation, tensor decompositions, and the application of tensor methods in the context of evolution operators. The AMUSEt algorithm is introduced in Section~\ref{sec:tensor_based_EDMD}, while the novel HOCUR-based algorithm is presented in Section~\ref{sec: HOCUR}. The theoretical analysis of AMUSEt follows in Section~ \ref{sec:analysis_amuset}, numerical results for benchmark problems are then shown in Section~\ref{sec:Numerical Examples}. Concluding remarks and open problems follow in Section~\ref{sec:Conclusion and Outlook}.

\section{Basic Concepts}
\label{sec:basic_concepts}

We first recapitulate basic concepts from dynamical systems theory, especially evolution operators and their Galerkin approximation in Section~\ref{subsec:evolution_op_approx}. Afterwards, we change topics and discuss low-rank approximations of tensors in the tensor train format and how these concepts can be applied in the context of dynamical systems in Section~\ref{subsec:TT}. Table~\ref{tab: notation} summarizes the notation that is used throughout the paper.

\begin{table}[h]
\renewcommand{\arraystretch}{1.3}
  \caption{Notation used in this work.}
  \centering
  \begin{tabular}{ll}
    \hline
    \textbf{Symbol} & \textbf{Description} \\
    \hline 
    $\mathbb{H}, L(\mathbb{H})$ & general Hilbert space, space of bounded linear operators \\
    $L^2_{\rho_0},\, L^2_{\rho_1}, \, L^2_{\mu}$ & weighted $L^2$-spaces \\
    $\mathbb{V}, \mathbb{W}$ & finite-dimensional Hilbert spaces or subspaces \\
    $\psi, \phi; \eta, \zeta$ & (orthonormal) bases of $\mathbb{V}$ and $\mathbb{W}$, respectively \\
    $\mathcal{P}_{\mathbb{V}}$ & orthogonal projection onto $\mathbb{V}$ \\
    $\mathcal{T}_\tau , \mathcal{K}_\tau, \mathcal{F}_\tau$ & Perron--Frobenius, Koopman, and forward-backward operator \\
    $\mathcal{T}_\tau(\mathbb{V}, \mathbb{W}) , \mathcal{K}_\tau(\mathbb{V}, \mathbb{W}), \mathcal{F}_\tau(\mathbb{V}, \mathbb{W})$ & Galerkin projections onto finite-dimensional subspaces \\
    $T_\tau(\psi, \phi), K_\tau(\psi, \phi), F_\tau(\psi, \phi)$ & matrix representations of Galerkin projections,\\
    $\Psi(\,\cdot\,)$ & transformed data matrix \\
    $C(\,\cdot\,), A(\,\cdot\,,\,\cdot\,)$ & covariance and cross-covariance matrices \\
    $\mathbb{G}_r, \mathbb{B}_r$ & spectral subspace and associated coefficient vector space \\
    $\mathbf{T} = \left\llbracket \mathbf{T}^{(1)}\right\rrbracket \otimes \dots \otimes \left\llbracket \mathbf{T}^{(p)}\right\rrbracket$ & tensor train of order $p$ \\
    $\mathbf{T}\vert_k =  \mat{\mathbf{T}}{n_1, \dots, n_k}{n_{k+1}, \dots, n_p}$ & mode-$k$ unfolding of $\mathbf{T}$ \\
    $\mathbf{\Psi}(\,\cdot\,)$ & transformed data tensor \\
    $\mathbf{C}(\,\cdot\,), \mathbf{A}(\,\cdot\,,\,\cdot\,)$ & covariance and cross-covariance tensors \\
    $\widehat{\,\cdot\,}$ & data-driven estimates\\
    \hline
  \end{tabular}
  \label{tab: notation}
\end{table}

\subsection{Evolution Operators and Their Approximation}
\label{subsec:evolution_op_approx}

\subsubsection{Evolution Operators}
\label{subsec:evolution_op_general}
The main concern of this study is the analysis of dynamical systems using evolution operators. Let $\mathcal{X}_t \in \R^d$  be a deterministic or stochastic dynamical system. For a positive lag time $\tau$, assume that the densities of $\mathcal{X}_t$ at times $t = 0$ and $t = \tau$ are given by $\rho_0$ and $\rho_1$, respectively. The \emph{Perron--Frobenius operator} $\mathcal{T}_\tau \colon L^2_{\rho_0} \to L^2_{\rho_1}$, and its adjoint, the \emph{Koopman operator} $\mathcal{K}_\tau = \mathcal{T}^*_\tau \colon L^2_{\rho_1} \to L^2_{\rho_0}$, are defined by
\begin{align*}
  \mathcal{T}_\tau f(y) &= \frac{1}{\rho_1(y)}\int f(x) \rho_0(x) p^\tau(x, y)\,\mathrm{d}x, &   \mathcal{K}_\tau f(x) &= \int p^\tau(x, y) f(y) \,\mathrm{d}y,
\end{align*}
where $p^\tau$ is the stochastic transition kernel associated with the process $\mathcal{X}_t$, see \cite{Koopman1931,DELLNITZ1999,SCHUETTE1999,Mezic2005,KLUS2016b} for more details. If the process is stationary and admits an invariant distribution $\mu$, it is convenient to choose $\rho_0 = \rho_1 = \mu$. In particular, if the dynamics $\mathcal{X}_t$ are reversible with respect to $\mu$, the operator $\mathcal{T}_\tau$ is self-adjoint on $L^2_\mu$, and hence identical to $\mathcal{K}_\tau$. For non-reversible or time-dependent systems, it is in general more appropriate to consider the \textit{forward-backward operator}
\begin{equation*}
\mathcal{F}_\tau = \mathcal{T}_\tau^* \mathcal{T}_\tau :\, L^2_{\rho_0} \to L^2_{\rho_0}.
\end{equation*}
One is often interested in certain spectral components of these operators, especially those which are largest in magnitude (also called \emph{leading / dominant} eigenvalues), as they can be used to determine dynamically long-lived structures. For stationary systems, leading eigenpairs of $\mathcal{K}_\tau$ or $\mathcal{T}_\tau$ (which will be close to one in absolute value), can be used for metastability analysis, especially if the system is reversible \cite{DEUFLHARD2005}. If the system is non-stationary, leading singular values and functions of the Koopman operator are typically used instead. If they exist, left singular functions are automatically eigenfunctions of $\mathcal{F}_\tau$, with eigenvalue equal to the square of the singular value. These functions can be used to determine finite-time coherent sets \cite{FROYLAND2010,FROYLAND2013,BK17:coherent,KHM19}.

\subsubsection{Galerkin Projection and Dimensionality Reduction}

\label{subsec:galerkin_proj}
For the purposes of numerical analysis, the linear operators above must be represented on finite-dimensional subspaces. To introduce the notation, let $\mathbb{H}$ be a Hilbert space of functions on $\mathbb{R}^d$ with inner product $\innerprod{\cdot}{\cdot}_\mathbb{H}$. Let $\mathbb{V}$ be a finite-dimensional subspace of dimension $n$, with a basis $\psi = (\psi_1, \ldots, \psi_n)^\top$. The orthogonal projector onto $\mathbb{V}$ is denoted by $\mathcal{P}_\mathbb{V}$. Any function in $\mathbb{V}$ can be represented uniquely by a vector of expansion coefficients with respect to $\psi$ in $\mathbb{R}^n$. More generally, any $r$-dimensional subspace $\mathbb{F}$ of $\mathbb{V}$, spanned by functions $\theta = (\theta_1, \ldots, \theta_r)^\top$, is encoded by a matrix $A \in \mathbb{R}^{n \times r}$, the columns of which contain the expansion coefficients of $\theta$ with respect to the basis $\psi$:
\[ \theta(x)^\top = \psi(x)^\top A. \]
In this context, we will call the column space of $A$ in $\mathbb{R}^n$ the \emph{coefficient vector space} corresponding to $\mathbb{F}$. Finally, the Gramian matrix of a finite basis is denoted by
\begin{equation*}
C(\psi) := \left(\innerprod{\psi_i}{\psi_j}_\mathbb{H} \right)_{ij}.
\end{equation*} 

\noindent Returning to the analysis of evolution operators, let finite-dimensional spaces $\mathbb{V} \subset L^2_{\rho_0}$ and $\mathbb{W} \subset L^2_{\rho_1}$, with bases $\psi = (\psi_{1},\dots, \psi_{n})^\top$ and $\phi = (\phi_{1},\dots, \phi_{n})^\top$, be given. We then consider the Galerkin projections
\begin{equation} \label{eq:Galerkin_Projection_Op}
\begin{split}
\mathcal{K}_\tau(\mathbb{V}, \mathbb{W}) &= \mathcal{P}_\mathbb{V} \mathcal{K}_t \mathcal{P}_\mathbb{W}, \\
\mathcal{T}_\tau(\mathbb{V}, \mathbb{W}) &= \mathcal{P}_\mathbb{W} \mathcal{T}_t \mathcal{P}_\mathbb{V}, \\
\mathcal{F}_\tau(\mathbb{V}, \mathbb{W}) &= \mathcal{K}_\tau(\mathbb{V}, \mathbb{W}) \mathcal{T}_\tau(\mathbb{V}, \mathbb{W}).
\end{split}
\end{equation}
The matrix representations of these operators with respect to the bases $\psi$ and $\phi$ are given by
\begin{equation} \label{eq:Galerkin_Projection_Mat}
\begin{split}
K_\tau(\psi, \phi) &= \big(C(\psi)\big)^{-1} A(\psi, \phi), \\
T_\tau(\psi, \phi) &= \big(C(\phi)\big)^{-1} A(\psi, \phi)^\top, \\
F_\tau(\psi, \phi) &= K_\tau(\psi, \phi) T_\tau(\psi, \phi),
\end{split}
\end{equation}
where $C(\psi)$ and $C(\phi)$ are the Gramians in $L^2_{\rho_0}$ and $L^2_{\rho_1}$, respectively, and the matrix $A(\psi, \phi)$ satisfies
\begin{align*}
A(\psi, \phi)_{ij} &= \langle \psi_{i}, \mathcal{K}_\tau \phi_{j}\rangle_{\rho_0} = \langle \mathcal{T}_\tau\psi_{i}, \phi_{j}\rangle_{\rho_1},
\end{align*}
see  \cite{WILLIAMS2015,KLUS2016b,KLUS2018b}. We will frequently use orthonormal bases to represent the operators introduced above. If $\eta$ and $\zeta$ are orthonormal bases for $\mathbb{V}$ and $\mathbb{W}$, the corresponding matrix representations reduce to
\begin{align*}
K_\tau(\eta, \zeta) &= A(\eta, \zeta) = T_\tau(\eta, \zeta)^\top, & F_\tau(\eta, \zeta) &= A(\eta, \zeta)A(\eta, \zeta)^\top.
\end{align*}
We observe that the matrix representation $K_\tau(\eta, \zeta) = A(\eta, \zeta)$ of the projected Koopman operator for orthonormal basis sets serves as a baseline model from which most relevant quantities can be calculated directly, such as eigenpairs in the stationary case, or singular pairs and forward-backward eigenpairs in the non-stationary case. Therefore, representations of this type will be particularly interesting in what follows. Arbitrary bases can always be transformed into orthonormal ones using spectral decompositions of the Gramians: If $C(\psi) = U_\psi \Sigma^2_\psi U_\psi^\top$, $C(\phi) = U_\phi \Sigma^2_\phi U_\phi^\top$, then we arrive at the following orthonormal bases (\emph{whitening transformation}):
\begin{align*}
\eta(x)^\top &= \psi(x)^\top U_\psi \Sigma_\psi^{-1}, & \zeta(x)^\top &= \phi(x)^\top U_\phi \Sigma_\phi^{-1}.
\end{align*}
By additionally truncating the spectral decompositions after their first $r \leq n$ components, reduced orthonormal bases $\eta^\top_r = \psi^\top (U_{\psi, r} \Sigma_{\psi, r}^{-1})$ and $\zeta^\top_r = \phi^\top (U_{\phi, r} \Sigma_{\phi, r}^{-1})$ can be obtained. The Koopman matrix corresponding to these bases can be calculated according to
\[
M_{\tau, r} := K_\tau(\eta_r, \zeta_r) = (\Sigma_{\psi, r}^{-1} U_{\psi, r}^\top) A(\psi, \phi) (U_{\phi, r} \Sigma_{\phi, r}^{-1}).\]
The matrix $M_{\tau, r}$ is then typically used as an approximate representation of $\mathcal{K}_\tau(\mathbb{V}, \mathbb{W})$. Quantities of interest can be expressed with respect to the full basis by means of the transformations $(U_{\psi, r} \Sigma_{\psi, r}^{-1})$ and $(U_{\phi, r} \Sigma_{\phi, r}^{-1})$. For instance, if $(\sigma_i, v_i, w_i)$ is a singular triplet of $M_{\tau, r}$, then the corresponding singular functions can be expressed in the original bases by the coefficient vectors $\xi_i = (U_{\psi, r} \Sigma_{\psi, r}^{-1})v_i$ and $\chi_i = (U_{\phi, r} \Sigma_{\phi, r}^{-1})w_i$. However, it should be kept in mind that $M_{\tau, r}$ really just represents the Koopman operator on the reduced subspaces spanned by $\eta_r, \, \zeta_r$, which will be analyzed in more detail in Section \ref{sec:analysis_amuset}. 

\subsubsection{Data-driven Approximation}
\label{subsec:data_driven_approx}
Usually, the integrals required for the matrices in \eqref{eq:Galerkin_Projection_Mat} cannot be computed analytically, and are estimated from data instead. Assume we have a pair of $\mathbb{R}^d$-valued random variables $x,\, y$ on a probability space $(\Omega, \mathcal{B}, \mathbb{P})$ such that the joint distribution of $(x, y)$ on $\mathbb{R}^d \times \mathbb{R}^d$ is $\vartheta(x, y) = \rho_0(x) p^\tau(x, y)$. For Hilbert space-valued, $\vartheta$-integrable functions $f \colon \mathbb{R}^d \times \mathbb{R}^d \to \mathbb{H}$, the strong law of large numbers \cite{Chacon1962} then implies that for almost surely any sequence of i.i.d. samples $(x_k, y_k) \in \mathbb{R}^d \times \mathbb{R}^d$, we have
\begin{align}
\label{eq:strong_law_sampling}
\lim_{m \rightarrow \infty} \frac{1}{m}\sum_{k=1}^m f(x_k, y_k) = \mathbb{E}^\vartheta[f].
\end{align}
For such an i.i.d.\ sequence, we can assemble all pairs into data matrices $X,\,Y\in \mathbb{R}^{d\times m}$, where $X=\left[x_1, \dots, x_m\right]$ and $Y=\left[y_1, \dots, y_m\right]$. For finite-dimensional subspaces $\mathbb{V} \subset L^2_{\rho_0}$ , $\mathbb{W} \subset L^2_{\rho_1}$ as above, we then define the \emph{transformed data matrices} in $\mathbb{R}^{n\times m}$ by
\begin{equation}
  \label{eq: transformed data matrices}
  \begin{split}
    \Psi(X) =
    \begin{bmatrix}
      \psi(x_1) & \dots & \psi(x_m)
    \end{bmatrix}
  \end{split}
  \quad \text{and} \quad
  \begin{split}
    \Phi(Y) =
    \begin{bmatrix}
      \phi(y_1) & \dots & \phi(y_m)
    \end{bmatrix}
  \end{split}.
\end{equation}
These give rise to the following empirical estimates of the matrices $C(\psi),\, C(\phi), \, A(\psi, \phi)$:
\begin{align}
\label{eq:empirical_estimates}
\widehat{C}(\psi) &= \frac{1}{m}\Psi(X)\Psi(X)^\top, & \widehat{C}(\phi) &= \frac{1}{m}\Phi(Y)\Phi(Y)^\top, & \widehat{A}(\psi, \phi) &= \frac{1}{m}\Psi(X)\Phi(Y)^\top.
\end{align}
One of the central results of the Koopman approach is \cite{WILLIAMS2015,KLUS2016b,KLUS2018b}:
\begin{Proposition}
\label{prop:ergodicity}
If \eqref{eq:strong_law_sampling} holds, then almost surely
\begin{align}
\label{eq:ergodicity_Cpsi}
C(\psi) &= \lim_{m \rightarrow \infty} \widehat{C}(\psi),  & C(\phi) &= \lim_{m \rightarrow \infty} \widehat{C}(\phi), & A(\psi, \phi) &= \lim_{m \rightarrow \infty} \widehat{A}(\psi, \phi).
\end{align}
\end{Proposition}

\begin{Remark}~
\begin{itemize}
 \item[(i)] A standard way to generate the samples $(x_k, y_k)$ is to draw $x_k$ i.i.d.\ from $\rho_0$, and $y_k$ is then obtained by integrating the dynamics $\mathcal{X}_t$ over time $\tau$, starting from $x_k$.
 \item[(ii)] By ergodic theory \cite{Chacon1962}, the conclusions of Proposition \ref{prop:ergodicity} also hold if $\mathcal{X}_t$ is stationary with unique invariant density $\mu = \rho_0 = \rho_1$. In this case, $x_k$ is chosen as the $k$-th step of any discretized trajectory of $\mathcal{X}_t$, and $y_k$ is obtained as the $k$-th step of the same trajectory shifted by $\tau$. 
\end{itemize}
\end{Remark}

By means of the empirical estimates \eqref{eq:empirical_estimates}, we can obtain data-based projections $\widehat{\mathcal{K}}_\tau(\mathbb{V}, \mathbb{W})$, $\widehat{\mathcal{T}}_\tau(\mathbb{V}, \mathbb{W})$ and $\widehat{\mathcal{F}}_\tau(\mathbb{V}, \mathbb{W})$ of the evolution operators. Their matrix representations with respect to $\psi, \,\phi$ are the same as \eqref{eq:Galerkin_Projection_Mat}, only using empirical estimates, i.e.,
\begin{equation} \label{eq:Galerkin_Projection_Mat_Emp}
\begin{split}
\widehat{K}_\tau(\psi, \phi) &= \big(\widehat{C}(\psi)\big)^{-1} \widehat{A}(\psi, \phi), \\
\widehat{T}_\tau(\psi, \phi) &= \big(\widehat{C}(\phi)\big)^{-1} \widehat{A}(\psi, \phi)^\top, \\ 
\widehat{F}_\tau(\psi, \phi) &= \widehat{K}_\tau(\psi, \phi) \widehat{T}_\tau(\psi, \phi).
\end{split}
\end{equation}
Just as we did for the analytical Galerkin projections in Section \ref{subsec:galerkin_proj}, we can use empirically orthonormal bases $\widehat{\eta}, \, \widehat{\zeta}$ (i.e., $\widehat{C}(\widehat{\eta}) = \mathrm{Id}$ and $\widehat{C}(\widehat{\zeta}) = \mathrm{Id}$) to simplify these matrix approximations, obtaining the empirical Koopman matrix as $\widehat{K}_\tau(\widehat{\eta}, \widehat{\zeta}) = \widehat{A}(\widehat{\eta}, \widehat{\zeta})$. Spectral decompositions of the empirical Gramians can be used to find orthonormal bases just as described above, and truncations of the spectral decompositions lead to appropriate reduced matrices $\widehat{M}_{\tau, r}$. The data-driven matrices just introduced form the basis for a number of well-known numerical methods to analyze evolution operators. In particular, \emph{extended dynamic mode decomposition (EDMD)} \cite{WILLIAMS2015,KLUS2016b} and the {variational approach to conformational dynamics (VAC)} \cite{NOE2013,NUESKE2014} apply to the stationary case, while \emph{canonical correlation analysis} \cite{HARDOON2004,KHM19} has been formulated for the non-stationary setting. See also \cite{KLUS2018b} for an overview of the nomenclature.

\subsubsection{The AMUSE Algorithm}
\label{subsec:amuse}
In the data-driven setting, calculation of the empirical Gramians can be entirely avoided. Let rank-$r$ singular value decompositions (SVDs) of the transformed data matrices be given by
\begin{align*}
\Psi(X) &= \widehat{U}_{X, r} \widehat{\Sigma}_{X, r} \widehat{V}_{X, r}^\top + \widehat{E}_{X, r}, & \Phi(Y) &= \widehat{U}_{Y, r} \widehat{\Sigma}_{Y, r} \widehat{V}_{Y, r}^\top + \widehat{E}_{Y, r},
\end{align*}
where $\widehat{E}_{X, r}, \, \widehat{E}_{Y,r}$ are the errors resulting from truncation of the SVD at rank $r \leq n$. Because of~\eqref{eq:empirical_estimates}, the basis sets $\widehat{\eta}^\top_r = \sqrt{m} \psi^\top (\widehat{U}_{X, r} \widehat{\Sigma}_{X, r}^{-1})$, $\widehat{\zeta}^\top_r = \sqrt{m} \phi^\top (\widehat{U}_{Y, r} \widehat{\Sigma}_{Y, r}^{-1})$ are empirically orthonormal. Moreover, their empirical Koopman matrix can be obtained directly from the above SVDs, by observing that
\begin{align*}
\widehat{M}_{\tau, r} = \widehat{A}(\widehat{\eta}_r, \widehat{\zeta}_r) = (\sqrt{m} \widehat{\Sigma}_{X, r}^{-1} \widehat{U}_{X, r}^\top) \frac{1}{m} \Psi(X) \Phi(Y)^\top \sqrt{m} (\widehat{U}_{Y, r} \widehat{\Sigma}_{Y, r}^{-1}) = \widehat{V}_{X, r}^\top \widehat{V}_{Y, r}.
\end{align*}
Quantities of interest, such as singular triplets, can be calculated directly from $\widehat{M}_{\tau, r}$. The method is summarized in Algorithm \ref{alg:AMUSE_FB} \cite{TONG1990,KLUS2018b,WU2020}. It should be kept in mind, though, that $\widehat{M}_{\tau, r}$ only serves as empirical approximation of the Koopman operator on the subspaces spanned by $\widehat{\eta}_r, \widehat{\zeta}_r$. Further below, these spaces will be analyzed in more detail.
\begin{algorithm}
  \caption{AMUSE}
  \label{alg:AMUSE_FB}
  \setlength{\tabcolsep}{.5ex}
  \begin{tabular}{ll}
    \textbf{Input:} & transformed data matrices $\Psi(X)$ and $\Phi(Y)$\\
    \textbf{Output:} & reduced Koopman matrix $\widehat{M}_{\tau, r} = \widehat{K}_\tau(\widehat{\eta}_r, \widehat{\zeta}_r)$, \\&  approximate singular values $\widehat{\sigma}_i$ and singular vectors  $\widehat{\xi}_i, \, \widehat{\chi}_i$ of $\widehat{\mathcal{K}}_\tau(\mathbb{V}, \mathbb{W})$.
  \end{tabular}
  \hrule\vspace{0.2cm}
  \begin{algorithmic}[1]
    \State Compute reduced SVDs of $\Psi(X)$ and $\Phi(Y)$, i.e., $\Psi(X) \approx \widehat{U}_{X, r}  \widehat{\Sigma}_{X, r}  \widehat{V}_{X, r}^\top$ and $\Phi(Y) \approx \widehat{U}_{Y, r} \widehat{\Sigma}_{Y, r}  \widehat{V}_{Y, r}^\top.$
    \State Compute $\widehat{M}_{\tau, r} = \widehat{V}_{X, r}^\top \widehat{V}_{Y, r}$.\label{algline:AMUSE_FB_M}
    \State Compute singular values $\widehat{\sigma}_i$ of $\widehat{M}_{\tau, r}$, and left and right singular vectors $\widehat{v}_i, \,\widehat{w}_i$.
    \State Express singular vectors with respect to original bases: $\widehat{\xi}_i =  \widehat{U}_{X, r} \,\widehat{\Sigma}_{X, r}^{-1} \, \widehat{v}_i$, $\widehat{\chi}_i = \widehat{U}_{Y, r} \, \widehat{\Sigma}_{Y, r}^{-1} \, \widehat{w}_i$.\label{algline:AMUSE_FB_xi}
  \end{algorithmic}
\end{algorithm}

If the subspaces $\mathbb{V}, \, \mathbb{W}$ are identical, it is often desirable to use a single reduced basis in Algorithm \ref{alg:AMUSE_FB}. The resulting modification using just an SVD of $\Psi(X)$ is shown in Algorithm \ref{alg:AMUSE_SINGLE}. If the process $\mathcal{X}_t$ is stationary ($\rho_0 = \rho_1 = \mu$), then the reduced subspaces based on SVDs of $\Psi(X)$ and $\Psi(Y)$ converge to the same limit with infinite data, hence Algorithms \ref{alg:AMUSE_FB} and \ref{alg:AMUSE_SINGLE} are asymptotically equivalent. For finite data however, their outputs will generally be different. We also note that the computation of the reduced matrix in Algorithm \ref{alg:AMUSE_SINGLE} does not break down to a single matrix product, which will be important when comparing the tensor-based versions of both algorithms in Section~\ref{sec:tensor_based_EDMD}.

\begin{algorithm}
  \caption{Single Basis AMUSE }
  \label{alg:AMUSE_SINGLE}
  \setlength{\tabcolsep}{.5ex}
  \begin{tabular}{ll}
    \textbf{Input:} & transformed data matrices $\Psi(X)$ and $\Psi(Y)$\\
    \textbf{Output:} & reduced Koopman matrix $\widehat{M}_{\tau, r} = \widehat{K}_\tau(\widehat{\eta}_r, \widehat{\eta}_r)$,
    \\&  approximate singular triplets $(\widehat{\sigma}_i, \widehat{\xi}_i, \, \widehat{\chi}_i)$ or eigenpairs $(\widehat{\lambda}_i, \widehat{\xi_i})$ of $\widehat{\mathcal{K}}_\tau(\mathbb{V}, \mathbb{V})$.
  \end{tabular}
  \hrule\vspace{0.2cm}
  \begin{algorithmic}[1]
    \State Compute reduced SVD of $\Psi(X)$, i.e., $\Psi(X) \approx \widehat{U}_{X, r} \widehat{\Sigma}_{X, r} \widehat{V}_{X, r}^\top$.
    \State Compute $\widehat{M}_{\tau, r} = \widehat{V}_{X, r}^\top \Psi(Y)^\top \widehat{U}_{X, r} \widehat{\Sigma}_{X, r}^{-1}$.\label{algline:AMUSE_SINGLE_M}
    \State Compute singular triplets $(\widehat{\sigma}_i, \widehat{v}_i, \,\widehat{w}_i)$ or eigenpairs $(\widehat{\lambda}_i, \widehat{w}_i)$ of $\widehat{M}_{\tau, r}$.
    \State Express singular vectors or eigenvectors w.r.t.\ original basis as in Algorithm \ref{alg:AMUSE_FB}.
  \end{algorithmic}
\end{algorithm}

\subsection{Low-rank Tensor Representations}
\label{subsec:TT}

Tensors are multi-dimensional arrays $ \T \in \R^{N}$, where $N = n_1 \times \dots \times  n_p$. Here, $p$ is called the order of a tensor, while the dimensions of the elementary vector spaces (the so-called modes) are $n_k$, $k = 1, \dots , p$. Tensor entries are sometimes represented by multi-indices $\mathbf{i} = (i_1,\dots,i_p) $ with $ i_k \in \{1, \dots, n_k\} $, i.e., $\T_\mathbf{i} = \T_{i_1, \dots, i_p}$. The single-index representation of the multi-index $\mathbf{i} $ is denoted by $\overline{\mathbf{i}} \in \{1, \dots, \prod_{k=1}^p n_k\}$. Conversely, the multi-index representation of the single-index $i \in \mathbb{N}$ is represented by $\underline{i}$. The tensor product is denoted by $\otimes$. We use bold capital letters ($\T, \mathbf{U}$, etc.) to denote tensors, capital letters ($U, V$, etc.) for matrices, and vectors are represented by lower case letters ($x, x_k$, etc.).  For $1\leq k \leq p-1$, the \emph{mode-$k$ unfolding} of the tensor $\T$ is the matrix 
\begin{equation*}
  \mat{\mathbf{T}}{n_1, \dots, n_k}{n_{k+1}, \dots, n_p} \in \R^{(n_1 \cdot \ldots \cdot n_k) \times (n_{k+1} \cdot \ldots \cdot n_p)}.
\end{equation*}
In short, we use the notation
\begin{equation*}
  \mathbf{T}\vert_k =  \mat{\mathbf{T}}{n_1, \dots, n_k}{n_{k+1}, \dots, n_p}
\end{equation*}
if the modes of $\mathbf{T}$ are clear.

\subsubsection{Tensor Train Format}
\label{subsec:tensor_train_format}
We start by introducing the tensor train (TT) format, where a high-dimensional tensor is represented by the contraction of multiple low-dimensional tensors \cite{OSELEDETS2009,OSELEDETS2011}.

\begin{Definition}
  A tensor $\mathbf{T} \in \R^N$ is said to be in the \emph{TT format} if
  \begin{equation*}
    \mathbf{T} = \sum_{l_0=1}^{r_0} \cdots  \sum_{l_p=1}^{r_p} \bigotimes_{k=1}^{p} \mathbf{T}^{(k)}_{l_{k-1},:,l_k} = \sum_{l_0=1}^{r_0} \cdots  \sum_{l_p=1}^{r_p}  \mathbf{T}^{(1)}_{l_0,:,l_1} \otimes \dots \otimes  \mathbf{T}^{(p)}_{l_{p-1},:,l_p}.
  \end{equation*}
  The tensors $\mathbf{T}^{(k)} \in \R^{r_{k-1} \times n_k \times r_k}$ of order 3 are called \emph{TT cores} and the numbers $r_k$ are called \emph{TT ranks}. It holds that $r_0 = r_p =1$ and $r_k \geq 1$ for $k=1, \dots, p-1$.
\end{Definition}

The TT ranks $r_0, \dots, r_p$ have a strong influence on the capability of representing a given tensor as a tensor train and determine the storage consumption of a tensor in the TT format. Figure~\ref{fig: tensor train} shows the graphical representation of a tensor train, which is also called Penrose notation, see~\cite{PENROSE1971}. 

\begin{figure}[htbp]
  \centering
  \begin{tikzpicture}
    \draw[black] (0,0) -- node [label={[shift={(0,-0.15)}]$r_1$}] {} ++ (1,0) ;
    \draw[black] (1,0) -- node [label={[shift={(0,-0.15)}]$r_2$}] {} ++ (1,0) ;
    \draw[black] (2,0) -- node [label={[shift={(0,-0.15)}]$r_3$}] {} ++ (1,0) ;
    \draw[black, dotted] (3,0) -- ++ (1,0) ;
    \draw[black] (4,0) -- node [label={[shift={(0,-0.15)}]$r_{p-2}$}] {} ++ (1,0) ;
    \draw[black] (5,0) -- node [label={[shift={(0,-0.15)}]$r_{p-1}$}] {} ++ (1,0) ;
    \draw[black] (0,0) -- node [label={[shift={(0,-1)}]$n_1$}] {} ++ (0,-0.7) ;
    \draw[black] (1,0) -- node [label={[shift={(0,-1)}]$n_2$}] {} ++ (0,-0.7) ;
    \draw[black] (2,0) -- node [label={[shift={(0,-1)}]$n_3$}] {} ++ (0,-0.7) ;
    \draw[black] (5,0) -- node [label={[shift={(0,-1)}]$n_{p-1}$}] {} ++ (0,-0.7) ;
    \draw[black] (6,0) -- node [label={[shift={(0,-1)}]$n_p$}] {} ++ (0,-0.7) ;
    \node[draw,shape=circle,fill=Gray, scale=0.65] at (0,0){};
    \node[draw,shape=circle,fill=Gray, scale=0.65] at (1,0){};
    \node[draw,shape=circle,fill=Gray, scale=0.65] at (2,0){};
    \node[draw,shape=circle,fill=Gray, scale=0.65] at (5,0){};
    \node[draw,shape=circle,fill=Gray, scale=0.65] at (6,0){};
  \end{tikzpicture}
  \caption{Graphical representation of a tensor train: A core is depicted by a circle with different arms indicating the modes of the tensor and
  the rank indices. The first and the last TT core are regarded as matrices due to the
  fact that $r_0 = r_p = 1$.}
  \label{fig: tensor train}
\end{figure}

We also represent TT cores as two-dimensional arrays containing vectors as elements. For a given tensor train $\mathbf{T} \in \R^N$ with cores $\mathbf{T}^{(k)} \in \R^{r_{k-1} \times n_k \times r_k}$, a single core is written as
\begin{equation*}
  \left\llbracket \mathbf{T}^{(k)} \right\rrbracket = 
  \left\llbracket
  \begin{matrix}
    & \mathbf{T}^{(k)}_{1,:,1} & \cdots & \mathbf{T}^{(k)}_{1,:,r_k} & \\
    & & & & \\
    & \vdots & \ddots & \vdots & \\
    & & & & \\
    & \mathbf{T}^{(k)}_{r_{k-1},:,1} & \cdots & \mathbf{T}^{(k)}_{r_{k-1},:,r_k} &
  \end{matrix}\right\rrbracket.
  \label{eq: core notation - single core}
\end{equation*}
We then use the notation $\mathbf{T} = \left\llbracket \mathbf{T}^{(1)}\right\rrbracket \otimes \dots \otimes \left\llbracket \mathbf{T}^{(p)}\right\rrbracket$ for representing tensor trains $\mathbf{T}$, cf.~\cite{GELSS2016, GELSS2019, GELSS2017}. This notation can be regarded as a generalization of the standard matrix multiplication. The difference is that we here compute the tensor products of the corresponding elements -- which are vectors instead of scalar values -- and then sum over the columns and rows, respectively. A core of a tensor train is \emph{left-orthonormal} if
\begin{equation*}
  \left( \mat{\T^{(k)}}{r_{k-1}, n_k}{r_k} \right)^\top \cdot \left(\mat{\T^{(k)}}{r_{k-1}, n_k}{r_k}\right) = \mathrm{Id} \in \R^{r_k \times r_k}.
\end{equation*}

\subsubsection{Basis Decompositions}
\label{subsec:basis decomp}

Within the context of data-driven approximation of evolution operators, see Section~\ref{subsec:data_driven_approx}, tensors arise if trial spaces are chosen as tensor products of elementary function spaces. We consider a data matrix $X \in \R^{d \times m}$ originating from a stochastic process $\mathcal{X}_t$, and a set of basis functions $\psi_1, \dots, \psi_p$ with $\psi_k \colon \R^d \rightarrow \R^{n_k}$ where $n_k \in \mathbb{N}$ for $k = 1, \dots, p$. Let $\mathbb{V}^k = \mathrm{span}\{\psi_{k,1}, \dots, \psi_{k,n_k}\}$ denote the $n_k$-dimensional subspaces spanned by the elementary basis functions $\psi_{k,1}, \dots, \psi_{k,n_k} \in L^2_{\rho_0}$. We consider the Galerkin projection \eqref{eq:Galerkin_Projection_Op} on the tensor product $\mathbb{V} := \mathbb{V}^1\otimes \dots \otimes \mathbb{V}^p\subset L^2_{\rho_0}$, which is a subspace of dimension at most $n_1 \cdot \ldots \cdot n_p$. Equivalently, choosing elementary basis functions in $L^2_{\rho_1}$ for a data matrix $Y$ yields a tensor space $\mathbb{W} \subset L^2_{\rho_1}$.

The tensor-based counterparts of the transformed data matrices given in \eqref{eq: transformed data matrices} are denoted by $\PSI(X)$ and $\PHI(Y)$, respectively. These \emph{transformed data tensors} can then be used to obtain empirical estimates of the Galerkin tensors for $\mathbb{V}, \mathbb{W}$:
\begin{align*}
\widehat{\mathbf{C}}(\PSI) &= \frac{1}{m} \PSI(X) \PSI(X)^\top, & \widehat{\mathbf{C}}(\PHI) &= \frac{1}{m} \PHI(Y) \PHI(Y)^\top, & \widehat{\mathbf{A}}(\PSI, \PHI) &= \frac{1}{m} \PSI(X) \PHI(Y)^\top,
\end{align*}
where the transposes of $\PSI(X)$ and $\PHI(Y)$ result from index permutations such that \linebreak ${\PSI(X)^\top, \PHI(Y)^\top \in \R^{m \times n_1 \times \dots \times n_p}}$. The multiplication above then denotes the contraction of the last mode of the first tensor with the first mode of the second tensor. As we have shown in \cite{GELSS2019}, the tensor train format can be used to represent transformed data tensors. In what follows, we will focus on the construction of $\PSI(X)$, the case for $\PHI(Y)$ is analogous. We start by considering rank-one tensors of the form
\begin{equation}
  \label{eq: basis decomposition}
  \PSI(x) = \psi_1(x) \otimes \dots \otimes \psi_p(x) = \begin{bmatrix} \psi_{1,1} (x) \\ \vdots \\ \psi_{1, n_1}(x) \end{bmatrix} \otimes \dots \otimes \begin{bmatrix} \psi_{p,1} (x) \\ \vdots \\ \psi_{p,n_p}(x) \end{bmatrix} \in \R^{n_1 \times n_2 \times \dots \times n_p}.
\end{equation}
Note that the so-called \emph{coordinate-major} and \emph{function-major} basis decompositions, introduced in~\cite{GELSS2019}, are special cases of the more general decomposition given in \eqref{eq: basis decomposition}. The transformed data tensor $\PSI(X) \in \R^{n_1 \times \dots \times n_p \times m}$ is then given by adding the rank-one decompositions \eqref{eq: basis decomposition} for all vectors $x_1, \dots, x_m$ and taking the tensor product with an additional unit vector. The result is the following TT decomposition:
\begin{equation}\label{eq: transformed data tensor}
  \begin{split}
    \PSI(X)  &= \sum_{k=1}^m \PSI(x_k) \otimes e_k \\
    &= \sum_{k=1}^m \psi_1(x_k) \otimes \dots \otimes \psi_p(x_k) \otimes e_k \\
    &= \left \llbracket \begin{matrix}
    \psi_1 (x_1) & \cdots & \psi_1 (x_m) 
    \end{matrix} \right\rrbracket \otimes
    \left\llbracket\begin{matrix}
    \psi_2 (x_1) & & 0 \\
    & \ddots & \\
    0 & & \psi_2 (x_m)
    \end{matrix} \right\rrbracket \otimes \cdots \\
    & \qquad \cdots \otimes
    \left \llbracket \begin{matrix}
    \psi_p (x_1) & & 0 \\
    & \ddots & \\
    0 & & \psi_p (x_m) 
    \end{matrix} \right\rrbracket \otimes
    \left \llbracket \begin{matrix}
    e_1 \\
    \vdots \\
    e_m
    \end{matrix} \right\rrbracket \\
     &=: \left \llbracket \PSI ^{(1)} (X)\right \rrbracket \otimes \left \llbracket \PSI ^{(2)} (X)\right \rrbracket \otimes \dots \otimes \left \llbracket \PSI    ^{(p)} (X)\right \rrbracket \otimes \left \llbracket \PSI ^{(p+1)} (X)\right \rrbracket,
  \end{split}
\end{equation}
where $e_k$, $k = 1, \dots , m$, denote the unit vectors of the standard basis in the $m$-dimensional Euclidean space. The matrix-based counterpart of $\PSI(X)$, see \eqref{eq: transformed data matrices}, would be given by the mode-$p$ unfolding
\begin{equation}
\label{eq: mode_p_unfolding_data_tensor}
  \Psi(X) = \PSI(X)\vert_p = \mat{\PSI(X)}{n_1, \dots, n_p}{m},
\end{equation}
that is, modes $n_1 , \dots, n_p$ represent row indices of the unfolding, and mode $m$ is the column index.\footnote{The relation~\eqref{eq: mode_p_unfolding_data_tensor} would still be satisfied if the position of the unit vectors in~\eqref{eq: transformed data tensor} is changed. Moreover, contracting the TT core of unit vectors with any other core results in the so-called \emph{block TT format}~\cite{DOLGOV2014}. For the application of HOCUR and HOSVD, however, we stick to the decomposition scheme given in~\eqref{eq: transformed data tensor}.} 

\subsubsection{Global SVD}
\label{subsec:global_svd}
In order to apply the Algorithms \ref{alg:AMUSE_FB} and \ref{alg:AMUSE_SINGLE} to tensor product bases, SVDs of the mode-$p$ unfoldings of $\PSI(X)$ and $\PHI(Y)$ are required. A multi-linear analogue of the standard SVD was presented for the tensor train format in \cite{KLUS2018}. For a given tensor train $\PSI(X)$, the method provides an orthonormal, tensor-train structured segment $\widehat{\mathbf{U}}_{X,r} \in \R^{n_1 \times \dots \times n_p \times r}$, a diagonal coupling matrix $\widehat{\Sigma}_{X,r} \in \R^{r \times r}$, and an orthonormal matrix $\widehat{V}_{X,r} \in \R^{m \times r}$ representing the last core. The dimension $r$ is the TT rank between the two last cores. Just like a standard SVD of a matrix, $\widehat{\mathbf{U}}_{X,r}, \widehat{\Sigma}_{X,r}$, and $\widehat{V}_{X,r}$ then satisfy the following properties:
\begin{itemize}
  \item[(i)] $\PSI(X) = \widehat{\mathbf{U}}_{X,r} \widehat{\Sigma}_{X,r}  \widehat{V}_{X,r}^\top$,
  \item[(ii)] $(\widehat{\mathbf{U}}_{X,r}\vert_p)^\top  \widehat{\mathbf{U}}_{X,r} \vert_p = \widehat{V}_{X,r}^\top  \widehat{V}_{X,r} = \mathrm{Id} \in \R^{r \times r}$,
  \item[(iii)] $\widehat{\Sigma}_{X,r} \in \R^{r \times r}$ is a diagonal matrix.
\end{itemize}

\begin{algorithm}
  \caption{Global SVD}
  \label{alg:SVD}
  \setlength{\tabcolsep}{.5ex}
  \begin{tabular}{ll}
    \textbf{Input:} & transformed data tensor $\PSI(X) \in \R^{n_1 \times \dots \times n_{p} \times m}$ in TT format\\
    \textbf{Output:} & global SVD of $\PSI(X)$ in the form of $\widehat{\mathbf{U}}_{X,r} \widehat{\Sigma}_{X,r}  \widehat{V}_{X,r}^\top$
  \end{tabular}
  \hrule\vspace{0.2cm}
  \begin{algorithmic}[1]
    \For {$k=1, \dots, p-1$}
    \State Compute (truncated) SVD of $\mat{\PSI^{(k)}(X)}{2}{~}$, i.e., $\mat{\PSI^{(k)}(X)}{2}{~} = U \Sigma V^\top + E$, $\Sigma \in \R^{r_k \times r_k}$. \label{algline:svd_step}
    \State Set $\PSI^{(k)}(X) $ to a reshaped version of $U$.
    \State Set $\PSI^{(k+1)}(X)$ to a reshaped version of $ \Sigma V^\top  \mat{\PSI^{(k+1)}(X)}{1}{~}$.\label{algline:update_svd}
    \EndFor
    \State Compute (truncated) SVD of $\mat{\PSI^{(p)}(X)}{2}{~}$, i.e., $\mat{\PSI^{(p)}(X)}{2}{~} = U \Sigma V^\top + E$, $\Sigma \in \R^{r_p \times r_p}$.
    \State Set $\PSI^{(p)}(X)$ to a reshaped version of $U$.
    \State Set $\PSI^{(p+1)}(X) $ to $ V^\top$.
    \State Define $r = r_p$, $\widehat{\mathbf{U}}_{X,r} = \left\llbracket \PSI^{(1)} (X) \right\rrbracket \otimes \ldots \otimes \left\llbracket \PSI^{(p)} (X) \right\rrbracket$, $\widehat{\Sigma}_{X,r}=\Sigma$, and $\widehat{V}_{X,r} = V$.
  \end{algorithmic}
\end{algorithm}
The method described in Algorithm~\ref{alg:SVD} proceeds as follows: 
Similarly to the TT-SVD algorithm proposed in~\cite{OSELEDETS2011}, we left-orthonormalize the TT cores $\PSI^{(1)}(X), \dots, \PSI^{(p-1)}(X)$ using (truncated) SVDs.
Then, we decompose the TT core $\PSI^{(p)}(X)$, but this time retain the diagonal matrix containing the singular values and only shift the right-orthonormal matrix to the last core. 
This provides the components of the global SVD as shown in Figure~\ref{fig: SVD}.
Note that we do not require any right-orthonormalization of the last core as stated in~\cite{KLUS2018} since $\PSI^{(p+1)}(X)$ is simply a reshaped identity matrix.

\begin{figure}[htbp]
  \centering
    \begin{tikzpicture}
      \draw[black] (0,0) -- ++ (1,0) ;
      \draw[black] (1,0) -- ++ (1,0) ;
      \draw[black] (2,0) -- ++ (1,0) ;
      \draw[black, dotted] (3,0) -- ++ (1,0) ;
      \draw[black] (4,0) -- ++ (1,0) ;
      \draw[black] (5,0) -- node [label={[shift={(0,-0.15)}]$r$}] {} ++ (1,0) ;
      \draw[black] (6,0) -- node [label={[shift={(0,-0.15)}]$r$}] {} ++ (1,0) ;
      \draw[black] (0,0) -- node [label={[shift={(0,-1)}]$n_1$}] {} ++ (0,-0.7) ;
      \draw[black] (1,0) -- node [label={[shift={(0,-1)}]$n_2$}] {} ++ (0,-0.7) ;
      \draw[black] (2,0) -- node [label={[shift={(0,-1)}]$n_3$}] {} ++ (0,-0.7) ;
      \draw[black] (5,0) -- node [label={[shift={(0,-1)}]$n_p$}] {} ++ (0,-0.7) ;
      \draw[black] (7,0) -- node [label={[shift={(0,0.2)}]$m$}] {} ++ (0,0.7) ;
      \node[draw,shape=semicircle,rotate=135,fill=white, anchor=south,inner sep=2pt, outer sep=0pt, scale=0.75] at (0,0){}; 
      \node[draw,shape=semicircle,rotate=315,fill=Blue, anchor=south,inner sep=2pt, outer sep=0pt, scale=0.75] at (0,0){};
      \node[draw,shape=semicircle,rotate=135,fill=white, anchor=south,inner sep=2pt, outer sep=0pt, scale=0.75] at (1,0){}; 
      \node[draw,shape=semicircle,rotate=315,fill=Blue, anchor=south,inner sep=2pt, outer sep=0pt, scale=0.75] at (1,0){};
      \node[draw,shape=semicircle,rotate=135,fill=white, anchor=south,inner sep=2pt, outer sep=0pt, scale=0.75] at (2,0){}; 
      \node[draw,shape=semicircle,rotate=315,fill=Blue, anchor=south,inner sep=2pt, outer sep=0pt, scale=0.75] at (2,0){};
      \node[draw,shape=semicircle,rotate=135,fill=white, anchor=south,inner sep=2pt, outer sep=0pt, scale=0.75] at (5,0){}; 
      \node[draw,shape=semicircle,rotate=315,fill=Blue, anchor=south,inner sep=2pt, outer sep=0pt, scale=0.75] at (5,0){};
      \node[draw,shape=circle,fill=Orange, scale=0.65] at (6,0){};
      \node[draw,shape=semicircle,rotate=135,fill=Green, anchor=south,inner sep=2pt, outer sep=0pt, scale=0.75] at (7,0){}; 
      \node[draw,shape=semicircle,rotate=315,fill=white, anchor=south,inner sep=2pt, outer sep=0pt, scale=0.75] at (7,0){};
      \node[anchor=north] at (2.5,-1.2) {$\underbrace{\hspace*{5.8cm}}_{\widehat{\mathbf{U}}_{X,r}^{\vphantom{T}}}$};
      \node[anchor=north] at (6,-1.2) {$\underbrace{\hspace*{0.8cm}}_{\widehat{\Sigma}_{X,r}^{\vphantom{T}}}$};
      \node[anchor=north] at (7,-1.2) {$\underbrace{\hspace*{0.8cm}}_{\widehat{V}_{X,r}^\top}$};
    \end{tikzpicture}
  \caption{Global SVD in TT format: Global SVD of a transformed data tensor $\PSI(X)$ with components $\widehat{\mathbf{U}}_{X,r}$ (half-filled circles in blue), $\widehat{\Sigma}_{X,r}$ (orange circle), and $\widehat{V}_{X,r}$ (half-filled circle in green). The first $p$ modes depict the row indices while the last mode depicts the column index. We transpose the last TT core in order to account for the different row and column dimensions of the transformed data matrices, cf.~Section~\ref{subsec:basis decomp}.} 
  \label{fig: SVD}
\end{figure}

\section{AMUSE on Tensors}
\label{sec:tensor_based_EDMD}

In this section, we combine the ideas of the AMUSE algorithm developed in Section \ref{subsec:amuse} with the global SVD for tensor-structured bases. Recall that the idea of AMUSE was to use truncated SVDs of the data matrices in order to determine empirically orthonormal bases of reduced subspaces of the trial spaces $\mathbb{V}, \mathbb{W}$. We have seen that a matrix representation of the empirical Koopman operator on these spaces can then be found without computing or inverting any of the Gramian matrices. If $\mathbb{V}$ and $\mathbb{W}$ are tensor product spaces, this procedure rapidly becomes infeasible, as we would have to calculate SVDs of the mode-$p$ unfoldings of $\PSI(X),\, \PHI(Y)$, which grow exponentially in size. However, we have seen in Section \ref{subsec:global_svd} that an approximation to these SVDs can be obtained by the multi-linear global SVD algorithm.

In the following section, we derive a corresponding multi-linear AMUSE algorithm, which we will call AMUSEt (\emph{AMUSE on tensors}). This method only requires operations on individual cores of the TT representation \eqref{eq: transformed data tensor} and the contraction of a tensor network. The analysis of this method will then be presented in Section \ref{sec:analysis_amuset}.

\subsection{The AMUSEt Algorithm}
\label{subsec:amuset}

Suppose we have given TT representations of the transformed data tensors $\PSI(X)$ and $\PHI(Y)$. In order to construct the reduced matrix $\widehat{M}_{\tau, r}$ in Algorithm~\ref{alg:AMUSE_FB}, we first compute global SVDs of $\PSI(X)$ and $\PHI(Y)$, i.e.,
\begin{equation*}
 \PSI(X) = \widehat{\mathbf{U}}_{X,r} \widehat{\Sigma}_{X,r} \widehat{V}_{X,r}^\top \qquad \text{and} \qquad \PHI(Y) = \widehat{\mathbf{U}}_{Y,r} \widehat{\Sigma}_{Y,r} \widehat{V}_{Y,r}^\top,
\end{equation*}
by applying Algorithm~\ref{alg:SVD}. Analogously to the matrix case, the reduced matrix can simply be written as $\widehat{M}_{\tau, r} = \widehat{V}_{X,r}^\top \widehat{V}_{Y,r}$, see line~\ref{algline:AMUSE_FB_M} of Algorithm~\ref{alg:AMUSE_FB}. In what follows, we therefore focus on single basis AMUSEt. The construction of $\widehat{M}_{\tau, r}$ as in line~\ref{algline:AMUSE_SINGLE_M} of Algorithm~\ref{alg:AMUSE_SINGLE} leads to the tensor network shown in Figure~\ref{fig: AMUSE}. Note that here only the contractions of the TT cores of $\widehat{\mathbf{U}}_{X,r}$ and of $\widehat{\Sigma}_{X,r}$ with its inverse cancel out.

\begin{figure}[htbp]
    \centering
    \begin{tikzpicture}

        \def\x{0.45}
        \node[] at (\x,2) {$\Big\{$};
        \node[] at (\x,1) {$\Big\{$};
        \node[] at (\x,0) {$\Big\{$};
        \node[] at (\x,-1) {$\Big\{$};
        \def\x{0.1}
        \node[anchor=east] at (\x,2) {$\widehat{\mathbf{U}}_{X,r} \widehat{\Sigma}_{X,r}^{-1}$};
        \node[anchor=east] at (\x,1) {$\PSI(Y)^\top$};
        \node[anchor=east] at (\x,0) {$\PSI(X)$};
        \node[anchor=east] at (\x,-1) {$\widehat{\Sigma}_{X,r}^{-1} \widehat{\mathbf{U}}_{X,r}^\top$};
    
        \def\y{0}
        \draw[black] (1,\y) -- ++ (1.66,0) ;
        \draw[black, dotted] (2.66,\y) -- ++ (0.66,0) ;
        \draw[black] (3.33,\y) -- ++ (2.66,0) ;
        \draw[black] (1,\y) -- ++ (0,-1) ;
        \draw[black] (2,\y) -- ++ (0,-1) ;
        \draw[black] (4,\y) -- ++ (0,-1) ;
        \draw[black] (6,\y) -- ++ (0,1) ;
        \node[draw,shape=semicircle,rotate=135,fill=white, anchor=south,inner sep=2pt, outer sep=0pt, scale=0.75] at (1,\y){}; 
        \node[draw,shape=semicircle,rotate=315,fill=Blue, anchor=south,inner sep=2pt, outer sep=0pt, scale=0.75] at (1,\y){};
        \node[draw,shape=semicircle,rotate=135,fill=white, anchor=south,inner sep=2pt, outer sep=0pt, scale=0.75] at (2,\y){}; 
        \node[draw,shape=semicircle,rotate=315,fill=Blue, anchor=south,inner sep=2pt, outer sep=0pt, scale=0.75] at (2,\y){};
        \node[draw,shape=semicircle,rotate=135,fill=white, anchor=south,inner sep=2pt, outer sep=0pt, scale=0.75] at (4,\y){}; 
        \node[draw,shape=semicircle,rotate=315,fill=Blue, anchor=south,inner sep=2pt, outer sep=0pt, scale=0.75] at (4,\y){};
        \node[draw,shape=circle,fill=Orange, scale=0.65] at (5,\y){};
        \node[draw,shape=semicircle,rotate=135,fill=Green, anchor=south,inner sep=2pt, outer sep=0pt, scale=0.75] at (6,\y){}; 
        \node[draw,shape=semicircle,rotate=315,fill=white, anchor=south,inner sep=2pt, outer sep=0pt, scale=0.75] at (6,\y){};

        \def\y{-1}
        \draw[black] (1,\y) -- ++ (1.66,0) ;
        \draw[black, dotted] (2.66,\y) -- ++ (0.66,0) ;
        \draw[black] (3.33,\y) -- ++ (2.66,0) ;
        \node[draw,shape=semicircle,rotate=45,fill=white, anchor=south,inner sep=2pt, outer sep=0pt, scale=0.75] at (1,\y){}; 
        \node[draw,shape=semicircle,rotate=225,fill=Blue, anchor=south,inner sep=2pt, outer sep=0pt, scale=0.75] at (1,\y){};
        \node[draw,shape=semicircle,rotate=45,fill=white, anchor=south,inner sep=2pt, outer sep=0pt, scale=0.75] at (2,\y){}; 
        \node[draw,shape=semicircle,rotate=225,fill=Blue, anchor=south,inner sep=2pt, outer sep=0pt, scale=0.75] at (2,\y){};
        \node[draw,shape=semicircle,rotate=45,fill=white, anchor=south,inner sep=2pt, outer sep=0pt, scale=0.75] at (4,\y){}; 
        \node[draw,shape=semicircle,rotate=225,fill=Blue, anchor=south,inner sep=2pt, outer sep=0pt, scale=0.75] at (4,\y){};
        \node[draw,shape=circle,fill=Red, scale=0.65] at (5,\y){};

        \def\y{1}
        \draw[black] (1,\y) -- ++ (1,0);
        \draw[black] (2,\y) -- ++ (0.66,0) ;
        \draw[black, dotted] (2.66,\y) -- ++ (0.66,0) ;
        \draw[black] (3.33,\y) -- ++ (2.66,0) ;
        \draw[black] (1,\y) -- ++ (0,1) ;
        \draw[black] (2,\y) -- ++ (0,1) ;
        \draw[black] (4,\y) -- ++ (0,1) ;
        \node[draw,shape=circle,fill=gray, scale=0.65] at (1,\y){};
        \node[draw,shape=circle,fill=gray, scale=0.65] at (2,\y){};
        \node[draw,shape=circle,fill=gray, scale=0.65] at (4,\y){};
        \node[draw,shape=circle,fill=gray, scale=0.65] at (6,\y){};

        \def\y{2}
        \draw[black] (1,\y) -- ++ (1,0) ;
        \draw[black] (2,\y) -- ++ (0.66,0) ;
        \draw[black, dotted] (2.66,\y) -- ++ (0.66,0) ;
        \draw[black] (3.33,\y) -- ++ (2.66,0) ;
        \node[draw,shape=semicircle,rotate=135,fill=white, anchor=south,inner sep=2pt, outer sep=0pt, scale=0.75] at (1,\y){}; 
        \node[draw,shape=semicircle,rotate=315,fill=Blue, anchor=south,inner sep=2pt, outer sep=0pt, scale=0.75] at (1,\y){};
        \node[draw,shape=semicircle,rotate=135,fill=white, anchor=south,inner sep=2pt, outer sep=0pt, scale=0.75] at (2,\y){}; 
        \node[draw,shape=semicircle,rotate=315,fill=Blue, anchor=south,inner sep=2pt, outer sep=0pt, scale=0.75] at (2,\y){};
        \node[draw,shape=semicircle,rotate=135,fill=white, anchor=south,inner sep=2pt, outer sep=0pt, scale=0.75] at (4,\y){}; 
        \node[draw,shape=semicircle,rotate=315,fill=Blue, anchor=south,inner sep=2pt, outer sep=0pt, scale=0.75] at (4,\y){};
        \node[draw,shape=circle,fill=Red, scale=0.65] at (5,\y){};

        \node[] at (7,0.5) {$=$};

        \def\y{0}
        \draw[black] (8,\y) -- ++ (1,0);
        \draw[black] (9,\y) -- ++ (0.66,0) ;
        \draw[black, dotted] (9.66,\y) -- ++ (0.66,0) ;
        \draw[black] (10.33,\y) -- ++ (2.66,0) ;
        \draw[black] (8,\y) -- ++ (0,1) ;
        \draw[black] (9,\y) -- ++ (0,1) ;
        \draw[black] (11,\y) -- ++ (0,1) ;
        \node[draw,shape=circle,fill=gray, scale=0.65] at (8,\y){};
        \node[draw,shape=circle,fill=gray, scale=0.65] at (9,\y){};
        \node[draw,shape=circle,fill=gray, scale=0.65] at (11,\y){};
        \node[draw,shape=circle,fill=gray, scale=0.65] at (11.66,\y){};
        \node[draw,shape=semicircle,rotate=90,fill=white, anchor=south,inner sep=2pt, outer sep=0pt, scale=0.75] at (12.33,\y){}; 
        \node[draw,shape=semicircle,rotate=270,fill=Green, anchor=south,inner sep=2pt, outer sep=0pt, scale=0.75] at (12.33,\y){};
        
        \def\y{1}
        \draw[black] (8,\y) -- ++ (1,0) ;
        \draw[black] (9,\y) -- ++ (0.66,0) ;
        \draw[black, dotted] (9.66,\y) -- ++ (0.66,0) ;
        \draw[black] (10.33,\y) -- ++ (2.66,0) ;
        \node[draw,shape=semicircle,rotate=135,fill=white, anchor=south,inner sep=2pt, outer sep=0pt, scale=0.75] at (8,\y){}; 
        \node[draw,shape=semicircle,rotate=315,fill=Blue, anchor=south,inner sep=2pt, outer sep=0pt, scale=0.75] at (8,\y){};
        \node[draw,shape=semicircle,rotate=135,fill=white, anchor=south,inner sep=2pt, outer sep=0pt, scale=0.75] at (9,\y){}; 
        \node[draw,shape=semicircle,rotate=315,fill=Blue, anchor=south,inner sep=2pt, outer sep=0pt, scale=0.75] at (9,\y){};
        \node[draw,shape=semicircle,rotate=135,fill=white, anchor=south,inner sep=2pt, outer sep=0pt, scale=0.75] at (11,\y){}; 
        \node[draw,shape=semicircle,rotate=315,fill=Blue, anchor=south,inner sep=2pt, outer sep=0pt, scale=0.75] at (11,\y){};
        \node[draw,shape=circle,fill=Red, scale=0.65] at (12,\y){};


    \end{tikzpicture}
    \caption{Graphical representation of the reduced matrix: Given $\PSI(X)$ (in form of a global SVD) and $\PSI(Y)$, the reduced matrix $\widehat{M}_{\tau, r} \in \R^{r \times r}$ is computed by contracting the above tensor network. Half-filled circles in blue depict the cores of $\widehat{\mathbf{U}}_{X,r}$, $\widehat{\Sigma}_{X,r}$ and $\widehat{\Sigma}_{X,r}^{-1}$ are represented by orange and red circles, respectively, and the half-filled circle in green depicts the matrix $\widehat{V}_{X,r}$ of the global SVD. Since we here do not assume any further properties on $\PSI(Y)$, we simply represent its TT cores by gray circles.}
    \label{fig: AMUSE}
\end{figure}

The complexity of AMUSEt is mainly determined by the computational cost of the global SVD, which can be estimated as $O(p \min\{rn m^2, r^2n^2 m\})$, where $n$ is the maximum mode size. However, the complexity of AMUSEt can be reduced in many cases. Oftentimes, we consider snapshot matrices $X, Y \in \R^{d \times m}$, which are extracted from a trajectory data matrix ${Z\in \R^{d \times \tilde{m}}}$, $\tilde{m} >m $, and share a large number (i.e., close to $m$) of common snapshot vectors. Instead of constructing the transformed data tensors separately, we can construct the TT decomposition of $\PSI(Z)$ and then simply restrict the last TT core to the respective time steps in order to obtain the representations for $\PSI(X)$ and $\PSI(Y)$. That is, given the data matrix $Z$, and index sets $I_X, I_Y$ such that $X = Z_{:,I_X}$ and $Y = Z_{:, I_Y}$, we directly form $\PSI(Z)$ according to \eqref{eq: transformed data tensor}:
\begin{equation*}
  \PSI(Z) = \left \llbracket \PSI ^{(1)} (Z)\right \rrbracket \otimes \left \llbracket \PSI ^{(2)} (Z)\right \rrbracket \otimes \dots \otimes \left \llbracket \PSI ^{(p)} (Z)\right \rrbracket \otimes \left \llbracket \PSI ^{(p+1)} (Z)\right \rrbracket.
\end{equation*}
Then, the tensor trains $\PSI(X)$ and $\PSI(Y)$ are given by 
\begin{align*}
  \PSI(X) &= \left \llbracket \PSI ^{(1)} (Z)\right \rrbracket \otimes \left \llbracket \PSI ^{(2)} (Z)\right \rrbracket \otimes \dots \otimes \left \llbracket \PSI ^{(p)} (Z)\right \rrbracket \otimes \left \llbracket \left( \PSI ^{(p+1)} (Z) \right)_{:, I_X, 1}\right \rrbracket \\
  \PSI(Y) &= \left \llbracket \PSI ^{(1)} (Z)\right \rrbracket \otimes \left \llbracket \PSI ^{(2)} (Z)\right \rrbracket \otimes \dots \otimes \left \llbracket \PSI ^{(p)} (Z)\right \rrbracket \otimes \left \llbracket \left( \PSI ^{(p+1)} (Z) \right)_{:, I_Y, 1}\right \rrbracket.
\end{align*}
In order to arrive at a global SVD of $\PSI(X)$, we first left-orthonormalize the tensor train $\PSI(Z)$, from which $\PSI(X)$ and $\PSI(Y)$ can still be recovered by restricting the last core to the corresponding indices. Additionally, we compute a rank-$r$ SVD of only the last core of $\PSI(X)$, thereby completing the global SVD of $\PSI(X)$ in the form $\widehat{\mathbf{U}}_{X,r}  \widehat{U}_{X,r}  \widehat{\Sigma}_{X,r}  \widehat{V}_{X,r}^\top$. These transformations are visualized in Figure~\ref{fig: TT transformation}.

\begin{figure}[htbp]
    \centering
    \begin{tikzpicture}

        \def\x{0}
        \node[] at (\x,2) {$\PSI(Z)$}; 
        \node[] at (\x,1) {$\PSI(X)$};
        \node[] at (\x,0) {$\PSI(Y)$};
        
        \def\x{1}
        \node[] at (\x,2) {$=$};
        \node[] at (\x,1) {$=$};
        \node[] at (\x,0) {$=$};
        
        \def\x{2.85}
        \node[] at (\x,2) {$\widehat{\mathbf{U}}_{X,r} M_Z$};
        \node[] at (\x,1) {$\widehat{\mathbf{U}}_{X,r} \widehat{U}_{X,r} \widehat{\Sigma}_{X,r} \widehat{V}_{X,r}^\top$};
        \node[] at (\x,0) {$\widehat{\mathbf{U}}_{X,r} \widehat{M}_{Y,r}$};
        
        \def\x{4.7}
        \node[] at (\x,2) {$=$};
        \node[] at (\x,1) {$=$};
        \node[] at (\x,0) {$=$};
    
        \def\x{5.7}
        
        \def\y{2}
        \draw[black] (\x,\y) -- ++ (1.66,0) ;
        \draw[black, dotted] (\x+1.66,\y) -- ++ (0.66,0) ;
        \draw[black] (\x+2.33,\y) -- ++ (1.66,0) ;
        \draw[black] (\x,\y) -- ++ (0,-0.5) ;
        \draw[black] (\x+1,\y) -- ++ (0,-0.5) ;
        \draw[black] (\x+3,\y) -- ++ (0,-0.5) ;
        \draw[black] (\x+4,\y) -- ++ (0,0.5) ;
        \node[draw,shape=semicircle,rotate=135,fill=white, anchor=south,inner sep=2pt, outer sep=0pt, scale=0.75] at (\x,\y){}; 
        \node[draw,shape=semicircle,rotate=315,fill=Blue, anchor=south,inner sep=2pt, outer sep=0pt, scale=0.75] at (\x,\y){};
        \node[draw,shape=semicircle,rotate=135,fill=white, anchor=south,inner sep=2pt, outer sep=0pt, scale=0.75] at (\x+1,\y){}; 
        \node[draw,shape=semicircle,rotate=315,fill=Blue, anchor=south,inner sep=2pt, outer sep=0pt, scale=0.75] at (\x+1,\y){};
        \node[draw,shape=semicircle,rotate=135,fill=white, anchor=south,inner sep=2pt, outer sep=0pt, scale=0.75] at (\x+3,\y){}; 
        \node[draw,shape=semicircle,rotate=315,fill=Blue, anchor=south,inner sep=2pt, outer sep=0pt, scale=0.75] at (\x+3,\y){};
        \node[draw,shape=circle,fill=Gray, scale=0.65] at (\x+4,\y){};
        
        \def\y{1}
        \draw[black] (\x,\y) -- ++ (1.66,0) ;
        \draw[black, dotted] (\x+1.66,\y) -- ++ (0.66,0) ;
        \draw[black] (\x+2.33,\y) -- ++ (3.66,0) ;
        \draw[black] (\x,\y) -- ++ (0,-0.5) ;
        \draw[black] (\x+1,\y) -- ++ (0,-0.5) ;
        \draw[black] (\x+3,\y) -- ++ (0,-0.5) ;
        \draw[black] (\x+6,\y) -- ++ (0,0.5) ;
        \node[draw,shape=semicircle,rotate=135,fill=white, anchor=south,inner sep=2pt, outer sep=0pt, scale=0.75] at (\x,\y){}; 
        \node[draw,shape=semicircle,rotate=315,fill=Blue, anchor=south,inner sep=2pt, outer sep=0pt, scale=0.75] at (\x,\y){};
        \node[draw,shape=semicircle,rotate=135,fill=white, anchor=south,inner sep=2pt, outer sep=0pt, scale=0.75] at (\x+1,\y){}; 
        \node[draw,shape=semicircle,rotate=315,fill=Blue, anchor=south,inner sep=2pt, outer sep=0pt, scale=0.75] at (\x+1,\y){};
        \node[draw,shape=semicircle,rotate=135,fill=white, anchor=south,inner sep=2pt, outer sep=0pt, scale=0.75] at (\x+3,\y){}; 
        \node[draw,shape=semicircle,rotate=315,fill=Blue, anchor=south,inner sep=2pt, outer sep=0pt, scale=0.75] at (\x+3,\y){};
        \node[draw,shape=semicircle,rotate=90,fill=white, anchor=south,inner sep=2pt, outer sep=0pt, scale=0.75] at (\x+4,\y){}; 
        \node[draw,shape=semicircle,rotate=270,fill=Blue, anchor=south,inner sep=2pt, outer sep=0pt, scale=0.75] at (\x+4,\y){};
        \node[draw,shape=circle,fill=Orange, scale=0.65] at (\x+5,\y){};
        \node[draw,shape=semicircle,rotate=135,fill=Green, anchor=south,inner sep=2pt, outer sep=0pt, scale=0.75] at (\x+6,\y){}; 
        \node[draw,shape=semicircle,rotate=315,fill=white, anchor=south,inner sep=2pt, outer sep=0pt, scale=0.75] at (\x+6,\y){};
        
        \def\y{0}
        \draw[black] (\x,\y) -- ++ (1.66,0) ;
        \draw[black, dotted] (\x+1.66,\y) -- ++ (0.66,0) ;
        \draw[black] (\x+2.33,\y) -- ++ (1.66,0) ;
        \draw[black] (\x,\y) -- ++ (0,-0.5) ;
        \draw[black] (\x+1,\y) -- ++ (0,-0.5) ;
        \draw[black] (\x+3,\y) -- ++ (0,-0.5) ;
        \draw[black] (\x+4,\y) -- ++ (0,0.5) ;
        \node[draw,shape=semicircle,rotate=135,fill=white, anchor=south,inner sep=2pt, outer sep=0pt, scale=0.75] at (\x,\y){}; 
        \node[draw,shape=semicircle,rotate=315,fill=Blue, anchor=south,inner sep=2pt, outer sep=0pt, scale=0.75] at (\x,\y){};
        \node[draw,shape=semicircle,rotate=135,fill=white, anchor=south,inner sep=2pt, outer sep=0pt, scale=0.75] at (\x+1,\y){}; 
        \node[draw,shape=semicircle,rotate=315,fill=Blue, anchor=south,inner sep=2pt, outer sep=0pt, scale=0.75] at (\x+1,\y){};
        \node[draw,shape=semicircle,rotate=135,fill=white, anchor=south,inner sep=2pt, outer sep=0pt, scale=0.75] at (\x+3,\y){}; 
        \node[draw,shape=semicircle,rotate=315,fill=Blue, anchor=south,inner sep=2pt, outer sep=0pt, scale=0.75] at (\x+3,\y){};
        \node[draw,shape=circle,fill=white, scale=0.65] at (\x+4,\y){};

    \end{tikzpicture}
    \caption{Construction of transformed data tensors: The tensor trains $\PSI(X)$ and $\PSI(Y)$ are extracted from the left-orthonormalized tensor train $\PSI(Z)$. Additionally, $\PSI(X)$ is represented by its global SVD. Again, the cores of $\widehat{\mathbf{U}}_{X,r}$ (as well as $\widehat{U}_{X,r}$) are represented by half-filled circles in blue, $\widehat{\Sigma}_{X,r}$ by an orange circle, and $\widehat{V}_{X,r}$ by a half-filled green circle. The last cores of $\PSI(Z)$ and $\PSI(Y)$ are depicted by a gray and a white circle, respectively. }
    \label{fig: TT transformation}
\end{figure}
By construction, $\PSI(X)$ and $\PSI(Y)$ share the same segment $\widehat{\mathbf{U}}_{X,r}$. Since the cores of $\widehat{\mathbf{U}}_{X,r}$ (as well as $\widehat{U}_{X,r}$) are left-orthonormal, most of the contractions cancel out, and only four matrices remain, namely $\widehat{V}_{X,r}$, $\widehat{M}_{Y,r}$, $\widehat{U}_{X,r}$, and $\widehat{\Sigma}_{X,r}^{-1}$. Thus, the reduced matrix is then simply given by
\begin{equation}\label{eq:reduced_ev_problem}
    \widehat{M}_{\tau,r} = \widehat{\Sigma}_{X,r}^{-1}  \widehat{U}_{X,r}^\top  \widehat{\mathbf{U}}_{X,r}^\top  \PSI(X)  \PSI(Y)^\top  \widehat{\mathbf{U}}_{X,r}  \widehat{U}_{X,r}   \widehat{\Sigma}_{X,r}^{-1} = \widehat{V}_{X,r}^\top \widehat{M}_{Y,r}^\top  \widehat{U}_{X,r} \widehat{\Sigma}_{X,r}^{-1}.
\end{equation}
Finally, let us discuss the expression of singular vectors or eigenvectors computed from $\widehat{M}_{\tau, r}$ with respect to the full tensor product basis. For instance, given the $q$ leading eigenvectors of the reduced matrix $\widehat{M}_{\tau,r}$ in form of a matrix $\widehat{W} = [\widehat{w}_1, \dots, \widehat{w}_q]$, the approximate eigentensors of $\widehat{\mathcal{K}}_\tau(\mathbb{V}, \mathbb{V})$ can be expressed as a tensor train $\mathbf{\Xi}$ with $\mathbf{\Xi} = \widehat{\mathbf{U}}_{X,r} \widehat{U}_{X,r} \widehat{\Sigma}_{X,r}^{-1} \widehat{W}$, see line \ref{algline:AMUSE_SINGLE_M} of Algorithm~\ref{alg:AMUSE_SINGLE}. The evaluations of the associated eigenfunctions at all snapshots are then given by the matrix $\mathbf{\Xi}^\top \PSI(X)$, see~\cite{KLUS2018b}. The corresponding tensor network also breaks down to a simple matrix product, as is shown in Figure~\ref{fig: AMUSE - eigenfunctions}.

\begin{figure}[htbp]
    \centering
    \begin{tikzpicture}

        \def\x{2}
    
        \def\y{2}
        \node[anchor=east] at (\x-1.6,\y) {$\mathbf{\Xi}$}; 
        \node[] at (\x-1,\y) {$=$};
        \draw[black] (\x,\y) -- ++ (1.66,0) ;
        \draw[black, dotted] (\x+1.66,\y) -- ++ (0.66,0) ;
        \draw[black] (\x+2.33,\y) -- ++ (3.66,0) ;
        \draw[black] (\x,\y) -- ++ (0,-0.5) ;
        \draw[black] (\x+1,\y) -- ++ (0,-0.5) ;
        \draw[black] (\x+3,\y) -- ++ (0,-0.5) ;
        \draw[black] (\x+6,\y) -- ++ (0,0.5) ;
        \node[draw,shape=semicircle,rotate=135,fill=white, anchor=south,inner sep=2pt, outer sep=0pt, scale=0.75] at (\x,\y){}; 
        \node[draw,shape=semicircle,rotate=315,fill=Blue, anchor=south,inner sep=2pt, outer sep=0pt, scale=0.75] at (\x,\y){};
        \node[draw,shape=semicircle,rotate=135,fill=white, anchor=south,inner sep=2pt, outer sep=0pt, scale=0.75] at (\x+1,\y){}; 
        \node[draw,shape=semicircle,rotate=315,fill=Blue, anchor=south,inner sep=2pt, outer sep=0pt, scale=0.75] at (\x+1,\y){};
        \node[draw,shape=semicircle,rotate=135,fill=white, anchor=south,inner sep=2pt, outer sep=0pt, scale=0.75] at (\x+3,\y){}; 
        \node[draw,shape=semicircle,rotate=315,fill=Blue, anchor=south,inner sep=2pt, outer sep=0pt, scale=0.75] at (\x+3,\y){};
        \node[draw,shape=semicircle,rotate=90,fill=white, anchor=south,inner sep=2pt, outer sep=0pt, scale=0.75] at (\x+4,\y){}; 
        \node[draw,shape=semicircle,rotate=270,fill=Blue, anchor=south,inner sep=2pt, outer sep=0pt, scale=0.75] at (\x+4,\y){};
        \node[draw,shape=circle,fill=Red, scale=0.65] at (\x+5,\y){};
        \node[draw,shape=rectangle,fill=white, scale=0.9] at (\x+6,\y){}; 
    
        \def\y{0}
        \node[anchor=east] at (\x-1.6,\y) {$\mathbf{\Xi}^\top \PSI(X)$}; 
        \node[] at (\x-1,\y) {$=$};
        
        \draw[black] (\x,\y+0.5) -- ++ (1.66,0) ;
        \draw[black, dotted] (\x+1.66,\y+0.5) -- ++ (0.66,0) ;
        \draw[black] (\x+2.33,\y+0.5) -- ++ (3.66,0) ;
        \draw[black] (\x,\y+0.5) -- ++ (0,-0.5) ;
        \draw[black] (\x+1,\y+0.5) -- ++ (0,-0.5) ;
        \draw[black] (\x+3,\y+0.5) -- ++ (0,-0.5) ;
        \draw[black] (\x+6,\y+0.5) -- ++ (0,+0.5) ;
        \node[draw,shape=semicircle,rotate=135,fill=white, anchor=south,inner sep=2pt, outer sep=0pt, scale=0.75] at (\x,\y+0.5){}; 
        \node[draw,shape=semicircle,rotate=315,fill=Blue, anchor=south,inner sep=2pt, outer sep=0pt, scale=0.75] at (\x,\y+0.5){};
        \node[draw,shape=semicircle,rotate=135,fill=white, anchor=south,inner sep=2pt, outer sep=0pt, scale=0.75] at (\x+1,\y+0.5){}; 
        \node[draw,shape=semicircle,rotate=315,fill=Blue, anchor=south,inner sep=2pt, outer sep=0pt, scale=0.75] at (\x+1,\y+0.5){};
        \node[draw,shape=semicircle,rotate=135,fill=white, anchor=south,inner sep=2pt, outer sep=0pt, scale=0.75] at (\x+3,\y+0.5){}; 
        \node[draw,shape=semicircle,rotate=315,fill=Blue, anchor=south,inner sep=2pt, outer sep=0pt, scale=0.75] at (\x+3,\y+0.5){};
        \node[draw,shape=semicircle,rotate=90,fill=white, anchor=south,inner sep=2pt, outer sep=0pt, scale=0.75] at (\x+4,\y+0.5){}; 
        \node[draw,shape=semicircle,rotate=270,fill=Blue, anchor=south,inner sep=2pt, outer sep=0pt, scale=0.75] at (\x+4,\y+0.5){};
        \node[draw,shape=circle,fill=Orange, scale=0.65] at (\x+5,\y+0.5){};
        \node[draw,shape=semicircle,rotate=135,fill=Green, anchor=south,inner sep=2pt, outer sep=0pt, scale=0.75] at (\x+6,\y+0.5){}; 
        \node[draw,shape=semicircle,rotate=315,fill=white, anchor=south,inner sep=2pt, outer sep=0pt, scale=0.75] at (\x+6,\y+0.5){};
        
        \draw[black] (\x,\y-0.5) -- ++ (1.66,0) ;
        \draw[black, dotted] (\x+1.66,\y-0.5) -- ++ (0.66,0) ;
        \draw[black] (\x+2.33,\y-0.5) -- ++ (3.66,0) ;
        \draw[black] (\x,\y-0.5) -- ++ (0,+0.5) ;
        \draw[black] (\x+1,\y-0.5) -- ++ (0,+0.5) ;
        \draw[black] (\x+3,\y-0.5) -- ++ (0,+0.5) ;
        \draw[black] (\x+6,\y-0.5) -- ++ (0,-0.5) ;
        \node[draw,shape=semicircle,rotate=45,fill=white, anchor=south,inner sep=2pt, outer sep=0pt, scale=0.75] at (\x,\y-0.5){}; 
        \node[draw,shape=semicircle,rotate=225,fill=Blue, anchor=south,inner sep=2pt, outer sep=0pt, scale=0.75] at (\x,\y-0.5){};
        \node[draw,shape=semicircle,rotate=45,fill=white, anchor=south,inner sep=2pt, outer sep=0pt, scale=0.75] at (\x+1,\y-0.5){}; 
        \node[draw,shape=semicircle,rotate=225,fill=Blue, anchor=south,inner sep=2pt, outer sep=0pt, scale=0.75] at (\x+1,\y-0.5){};
        \node[draw,shape=semicircle,rotate=45,fill=white, anchor=south,inner sep=2pt, outer sep=0pt, scale=0.75] at (\x+3,\y-0.5){}; 
        \node[draw,shape=semicircle,rotate=225,fill=Blue, anchor=south,inner sep=2pt, outer sep=0pt, scale=0.75] at (\x+3,\y-0.5){};
        \node[draw,shape=semicircle,rotate=90,fill=white, anchor=south,inner sep=2pt, outer sep=0pt, scale=0.75] at (\x+4,\y-0.5){}; 
        \node[draw,shape=semicircle,rotate=270,fill=Blue, anchor=south,inner sep=2pt, outer sep=0pt, scale=0.75] at (\x+4,\y-0.5){};
        \node[draw,shape=circle,fill=Red, scale=0.65] at (\x+5,\y-0.5){};
        \node[draw,shape=rectangle,fill=white, scale=0.9] at (\x+6,\y-0.5){}; 
        
        \node[] at (\x+7,\y) {$=$};
        
        \draw[black] (\x+8,\y-1) -- ++ (0,2) ;
        \node[draw,shape=semicircle,rotate=135,fill=Green, anchor=south,inner sep=2pt, outer sep=0pt, scale=0.75] at (\x+8,\y+0.5){}; 
        \node[draw,shape=semicircle,rotate=315,fill=white, anchor=south,inner sep=2pt, outer sep=0pt, scale=0.75] at (\x+8,\y+0.5){};
        \node[draw,shape=rectangle,fill=white, scale=0.9] at (\x+8,\y-0.5){};

    \end{tikzpicture}
    \caption{Graphical representation of the eigentensors and eigenfunctions: The tensor train $\mathbf{\Xi}$ is built by the contraction of $\widehat{\mathbf{U}}_{X,r}$, $\widehat{U}_{X,r}$, $\widehat{\Sigma}_{X,r}^{-1}$, and $\widehat{W}$ (depicted by the square). The matrix $\mathbf{\Xi}^\top \PSI(X)$ comprising the evaluations of the eigenfunctions at the given snapshots is constructed by multiplying the tensors $\PSI(X)$ and $\mathbf{\Xi}$. Similar to the construction of the reduced matrix only a two cores remain since the orthonormal cores cancel out and $\widehat{\Sigma}_{X,r}$ is multiplied by its inverse.}
    \label{fig: AMUSE - eigenfunctions}
\end{figure}

\begin{Remark}
  Note that we do not need the cores of $\widehat{\mathbf{U}}_{X,r}$ after the orthonormalization procedure if we are only interested in the approximated eigenfunctions. For both the construction of the reduced matrix $\widehat{M}_{\tau,r}$ as well as the eigenfunction evaluations $\mathbf{\Xi}^\top \PSI(X)$, see Figure~\ref{fig: AMUSE - eigenfunctions}, the TT segment $\widehat{\mathbf{U}}_{X,r}$ is not required due to its orthonormality. This significantly reduces the storage consumption since we are able to construct the left-orthonormalized version of $\PSI(X)$ step by step, i.e., we only need to store two TT cores in memory at the same time.
\end{Remark}

\section{HOCUR-based Approach}
\label{sec: HOCUR}

In practice, the direct construction of transformed data tensors as described in Section~\ref{subsec:basis decomp} may be infeasible due to a large number of basis functions or snapshots. Thus, an alternative isolation technique for the modes of $\PSI(X)$ and $\PHI(Y)$ has to be used in order to apply AMUSEt. Our idea to circumvent this problem is a combination of different techniques from \cite{OSELEDETS2010, GOREINOV2010, GOREINOV2001}, specifically adapted to transformed data tensors as described in Section \ref{subsec:basis decomp}. Based on so-called \emph{CUR decompositions}, i.e., representing a matrix in terms of appropriate row and column subsets, we propose an iterative technique in Section~\ref{sec: HOCUR for TDT}. The aim is to construct a low-rank TT decomposition of the transformed data tensor without storing the complete representation of $\PSI(X)$ as given in~\eqref{eq: transformed data tensor}.

\subsection{Higher-order CUR Decomposition}
\label{subsec:HOCUR_general}

For a matrix $M \in \R^{m \times n}$, a CUR decomposition consists of index sets $I, J$, as well as submatrices $C = M_{:,J}$, $U = M_{I,J}$, and $R=M_{I,:}$, such that $M \approx C \cdot U^{-1} \cdot R$, see Figure \ref{fig: CUR}.

\begin{figure}[htbp]
  \centering
  \begin{tikzpicture}
    \def\h{0.4}
    \def\cross#1#2{
      \draw[color=Red, line width=0.04cm] (#1-0.25*\h,#2+0.0225+0.25*\h) -- ++ (0.5*\h,-0.5*\h);
      \draw[color=Red, line width=0.04cm] (#1+0.25*\h,#2+0.0225+0.25*\h) -- ++ (-0.5*\h,-0.5*\h);
    }
    \draw[draw=Blue, line width=0.04cm] (\h,0.5*\h) -- ++ (0,-5*\h);
    \draw[draw=Blue, line width=0.04cm] (5*\h,0.5*\h) -- ++ (0,-5*\h);
    \draw[draw=Blue, line width=0.04cm] (7*\h,0.5*\h) -- ++ (0,-5*\h);
    \draw[draw=Green, line width=0.04cm] (-0.5*\h,-\h+0.0225) -- ++ (9*\h,0);
    \draw[draw=Green, line width=0.04cm] (-0.5*\h,-2*\h+0.0225) -- ++ (9*\h,0);
    \draw[draw=Green, line width=0.04cm] (-0.5*\h,-4*\h+0.0225) -- ++ (9*\h,0);
    \cross{\h}{-\h}
    \cross{\h}{-2*\h}
    \cross{\h}{-4*\h}
    \cross{5*\h}{-\h}
    \cross{5*\h}{-2*\h}
    \cross{5*\h}{-4*\h}
    \cross{7*\h}{-\h}
    \cross{7*\h}{-2*\h}
    \cross{7*\h}{-4*\h}
    \node[] at (-\h,-2*\h+0.033) {$\left[ \vphantom{\scalebox{4}{A}} \right.$};
    \node[] at (9*\h,-2*\h+0.033) {$\left. \vphantom{\scalebox{4}{A}} \right]$};
    \foreach \i in {0,...,8}
      \foreach \j in {0,...,4}
  \node[] at (\i*\h,-\j*\h) {\scalebox{2}{$\cdot$}};   
    \node[] at (10*\h,-2*\h) {$\approx$};
    \node[] at (11*\h,-2*\h+0.033) {$\left[ \vphantom{\scalebox{4}{A}} \right.$};
    \node[] at (15*\h,-2*\h+0.033) {$\left. \vphantom{\scalebox{4}{A}} \right]$};
    \draw[draw=Blue, line width=0.04cm] (12*\h,0) -- ++ (0,-4*\h);
    \draw[draw=Blue, line width=0.04cm] (13*\h,0) -- ++ (0,-4*\h);
    \draw[draw=Blue, line width=0.04cm] (14*\h,0) -- ++ (0,-4*\h);
    \node[] at (16*\h,-2*\h) {$\cdot$};
    \node[] at (17*\h,-2*\h) {$\left[ \vphantom{\scalebox{3}{A}} \right.$};
    \node[] at (21*\h+0.2,-2*\h+0.0375) {$\left. \vphantom{\scalebox{3}{A}} \right]^{-1}$};
    \foreach \i in {0,...,2}
      \foreach \j in {0,...,2}
  \cross{18*\h+\i*\h}{-\h-\j*\h}
    \node[] at (22*\h+0.2,-2*\h) {$\cdot$};
    \node[] at (23*\h+0.2,-2*\h) {$\left[ \vphantom{\scalebox{3}{A}} \right.$};
    \node[] at (33*\h+0.2,-2*\h) {$\left. \vphantom{\scalebox{3}{A}} \right]$};
    \draw[draw=Green, line width=0.04cm] (24*\h+0.2,-\h) -- ++ (8*\h,0);
    \draw[draw=Green, line width=0.04cm] (24*\h+0.2,-2*\h) -- ++ (8*\h,0);
    \draw[draw=Green, line width=0.04cm] (24*\h+0.2,-3*\h) -- ++ (8*\h,0);
  \end{tikzpicture}
  \caption{CUR decomposition: The matrix on the left-hand side is approximated by the matrix product $C \cdot U^{-1} \cdot R$, where $C$ (blue lines) is a column subset, $R$ (green lines) is a row subset, and $U$ (red crosses) is the intersection matrix.}
  \label{fig: CUR}
\end{figure}

There are different methods to find optimal sets of rows and columns, cf.~\cite{BOUTSIDIS2014}. An important subproblem is the following: given a set of column indices $J=\{j_1, \dots, j_r\}$ with $r \leq \min(m,n)$ (and $M_{:,J}$ having full column rank), find an optimal subset of row indices $I = \{i_1, \dots, i_r\}$. This NP-hard problem can be approximately solved by applying the \emph{maximum-volume principle} to $M_{:,J}$, so that the infinity norm of $M-M_{:,J} \cdot M_{I,J}^{-1} \cdot M_{I,:}$ is minimized over $I$, see~\cite{GOREINOV2010, GOREINOV2001}. We refer to this algorithm as \textrm{Maxvol} from now on.

As described in \cite{OSELEDETS2010}, the CUR decomposition can be generalized to a tensor if its mode-$k$ unfoldings are successively decomposed using CURs. This method presents an alternative to the decomposition \eqref{eq: transformed data tensor} if applied to the transformed data tensor $\PSI(X)$. However, the procedure requires pre-defined row and column subsets for each unfolding. The authors of \cite{OSELEDETS2010} also suggested an iterative algorithm to circumvent this problem: after initializing row and column subsets in some way, the method alternates between updating the column subsets while all row sets are fixed, and vice versa. This algorithm is the basis of the method we will present in Section \ref{sec: HOCUR for TDT}.

\subsection{HOCUR for Transformed Data Tensors}\label{sec: HOCUR for TDT}

As suggested in \cite{OSELEDETS2010}, Algorithm~\ref{alg: HOCUR} successively updates the row sets of the unfolded residual tensors during a forward loop, while all column sets are fixed. Then, column sets are updated during a backward loop, with all row sets fixed, and the entire procedure is repeated until convergence. The key insight, used in lines \ref{alg: HOCUR Line submatrix 1}--\ref{alg: HOCUR Line MAXVOL1} and \ref{alg: HOCUR Line submatrix 2}--\ref{alg: HOCUR Line MAXVOL2}, is that each update only operates on a small subtensor which is easily evaluated. Assume we are given a row set $\mathbf{I} = \mathbf{I}_q = \{\mathbf{i}_1 , \dots, \mathbf{i}_{r_q}\}$ of multi-indices comprising modes $n_1,\ldots,n_q$, and a column set $\mathbf{J} = \mathbf{J}_{q+2} = \{\mathbf{j}_1 , \dots, \mathbf{j}_s\}$ of multi-indices comprising modes $n_{q+2},\ldots, n_{p}, m$:
\begin{equation*}
  \begin{gathered}
    \mathbf{i}_1 , \dots, \mathbf{i}_{r_q} \in \{1, \dots, n_1\} \times \dots \times \{1, \dots, n_{q}\}, \\
    \mathbf{j}_1 , \dots, \mathbf{j}_s \in \{1, \dots, n_{q+2}\} \times \dots \times \{1, \dots, n_{p}\} \times \{1, \dots, m\}.
  \end{gathered}
\end{equation*}
Then, a new extended row set comprising the first $q+1$ modes
\begin{equation*}
    \mathbf{i}_1 , \dots, \mathbf{i}_{r_{q+1}} \in \{1, \dots, n_1\} \times \dots \times \{1, \dots, n_{q+1}\},
\end{equation*}
can be obtained by applying Algorithm \textrm{Maxvol} to the submatrix $\Psi(X)_{|\I,\J} \in \R^{r_q \cdot n_{q + 1 } \times s}$, given by
\begin{equation}
  \label{eq: submatrix}
  \Psi(X)_{| \I,\J} = \begin{bmatrix}
    \PSI(X)_{\mathbf{i}_1, :, \mathbf{j}_1} & \cdots & \PSI(X)_{\mathbf{i}_1, :, \mathbf{j}_s} \\
    \vdots & \ddots & \vdots \\
    \PSI(X)_{\mathbf{i}_{r_q}, :, \mathbf{j}_1} & \cdots & \PSI(X)_{\mathbf{i}_{r_q}, :, \mathbf{j}_s}
  \end{bmatrix}.
\end{equation}
This matrix is easily set up using basis function evaluations. More precisely, given multi-indices $\mathbf{i} = (i_1, \dots, i_q) \in \mathbf{I}$ and $\mathbf{j} = (i_{q+2}, \dots, i_{p}, k) \in \mathbf{J}$, entries $\PSI(X)_{\mathbf{i},:,\mathbf{j}}$ of $\Psi(X)_{| \I, \J}$ are given by
\begin{equation*}
  \begin{split}
    \PSI(X)_{\mathbf{i}, :, \mathbf{j}} &= \PSI(X)_{i_1, \dots, i_q, :, i_{q+2}, \dots i_{p}, k}\\
    &= \underbrace{\psi_{1, i_1}(x_k) \cdot \ldots \cdot \psi_{q, i_q}(x_k)}_{\in \R} \cdot \underbrace{\psi_{q+1 \vphantom{i_{q+1}}} (x_k)}_{\in \R^{n_{q+1}}} \cdot \underbrace{\psi_{q+2, i_{q+2}}(x_k) \cdot \ldots \cdot \psi_{p, i_p}(x_k)}_{\in \R}.
  \end{split} 
\end{equation*}
Note that the last entry of the column index $\mathbf{j}$ determines the snapshot $x_k$ where the product is evaluated.

\begin{algorithm}[htb]
  \caption{Higher-order CUR decomposition.}
  \label{alg: HOCUR}
  \setlength{\tabcolsep}{.5ex}
  \begin{tabular}{ll}
    \textbf{Input:} & data matrix $X = [x_1, \dots, x_m] \in \R^{d \times m}$, basis functions $\psi_{i, j_i}$, $i = 1, \dots , p$, \\
    & $j_i = 1, \dots , n_i$, maximum ranks $r_1, \dots, r_p$ with $r_q \leq n_{q+1}\cdot r_{q+1}$,\\
    & number of iterations $N$, multiplier $\alpha > 1$\\
    \textbf{Output:} & TT approximation of the transformed data tensor $\PSI(X)$
  \end{tabular}
  \hrule\vspace{0.2cm}
  \begin{algorithmic}[1]
    \State Set $n_{p+1} = m$, $r_0 = r_{p+1} = 1$, and $\I_0 = \{ \varnothing \}$.
    \State Define initial multi-index column sets $\J_2 , \dots, \J_{p+2}$.\label{alg: HOCUR Line multiplier}
    \For {$k=1, \dots, N$}
    \For {$l = 1, \dots, p$} \hfill \textit{(First half sweep)}
    \State Extract submatrix $M = \Psi(X)_{|\I_{l-1}, \J_{l+1}}$, see \eqref{eq: submatrix}.\label{alg: HOCUR Line submatrix 1}
    \If {$k=1$} 
    \State Find set of linearly independent columns $J$ of $M$ with $\left|J \right| \leq r_{l}$.\label{alg: HOCUR Line licols}
    \State Set $M$ to $M_{:, J}$ and $r_{l}$ to $\left| J \right|$.
    \EndIf
    \State Apply Algorithm \textrm{Maxvol} to $M$ to extract row set $I$.\label{alg: HOCUR Line MAXVOL1}
    \State Compute multi-index row set $\mathbf{I}_{l}$ from $\I_{l-1}$ and $I$, see \eqref{eq: updated row set}.\label{alg: HOCUR_I}
    \State Define core $\PSI(X)^{(l)}$ as $M \cdot M_{I,:}^{-1}$ reshaped as $\R^{r_{l-1} \times n_{l} \times r_{l}}$.\label{alg: HOCUR Line core 1}
    \EndFor
    \For {$l = p+1, \dots, 2$} \hfill \textit{(Second half sweep)}
    \State Extract submatrix $M = \Psi(X)_{|\I_{l-1}, \J_{l+1}}$ and reshape as $\R^{r_{l-1} \times n_l \cdot r_l}$.\label{alg: HOCUR Line submatrix 2}
    \State Apply Algorithm \textrm{Maxvol} to $M^\top$ to extract column set $J$ and set $r_{l-1} = \left| J \right|$.\label{alg: HOCUR Line MAXVOL2}
    \State Compute multi-index column set $\J_l$ from $J$ and $\J_{l+1}$.\label{alg: HOCUR_J}
    \State Define core $\PSI(X)^{(l)}$ as $M_{:,J}^{-1} \cdot M$ reshaped as $\R^{r_{l-1} \times n_l \times r_l}$.\label{alg: HOCUR Line core 2}
    \EndFor
    \EndFor
    \State Define first core $\PSI(X)^{(1)}$ as $\Psi(X)_{|\I_0, \J_2}$ reshaped as $\R^{1 \times n_1 \times r_1}$.
  \end{algorithmic}
\end{algorithm}

Let us elaborate on a few more details of Algorithm~\ref{alg: HOCUR}. First, note that multi-index sets for the construction of the submatrices $\Psi(X)_{| \I,\J}$ are nested sets by construction. In line \ref{alg: HOCUR_I}, after each application of Algorithm \textrm{Maxvol}, the resulting single-index set $I = \{i_1, \dots, i_{r_{q+1}}\}$ needs to be converted into a multi-index row set for modes $n_1,\ldots,n_{q+1}$. Given a multi-index row set $\mathbf{I}_{q} = \{\mathbf{i}_{1}, \dots , \mathbf{i}_{r_q}\}$ as above, each row of $\Psi(X)_{| \I_q, \J_{q+2}}$ can naturally be associated with a multi-index $(k_1, k_2 ) \in \{1, \dots, r_q\} \times \{1, \dots, n_{q + 1}\}$. Hence, we map each single-index $i_k\in I$ to a multi-index $\underline{i}_k = (\underline{i}_{k,1}, \underline{i}_{k,2})$ in $\{1, \dots, r_{q}\} \times \{1, \dots, n_{q+1}\}$, and then define the extended multi-index row set
\begin{equation}
  \label{eq: updated row set}
  \mathbf{I}_{q+1} = \{ (\mathbf{i}_{\underline{i}_{1,1}}, \underline{i}_{1,2}), \dots, (\mathbf{i}_{\underline{i}_{r_{q+1},1}}, \underline{i}_{r_{q+1},2}) \}.
\end{equation}
Column sets are updated analogously in line~\ref{alg: HOCUR_J} of Algorithm~\ref{alg: HOCUR}.

Second, the algorithm requires initial column sets which are generated in line \ref{alg: HOCUR Line multiplier}. While it was suggested in \cite{OSELEDETS2010} to pick these columns at random, we build them up recursively to ensure the column sets are also nested. Starting from $\mathbf{J}_{p+2} = \{\varnothing\}$, column set $\mathbf{J}_{q}$ is obtained by simply selecting the first $\min(\alpha \cdot r_{q}, n_{q+1} \cdot r_{q+1})$ indices out of the index set $\{ 1, \dots, n_{q+1} \cdot r_{q+1}\}$, and then joining them with multi-index column set $\mathbf{J}_{q+1}$ as described above in \eqref{eq: updated row set}. In practice, we found it helpful to select a rather large number of columns at this point, as these initial columns would often be highly redundant. The parameter $\alpha$ can be tuned to ensure enough columns are selected during the initialization stage.

Third, we need to find index sets of linearly independent columns of the matrices $\Psi(X)_{| \I,\J}$ during the first iteration, see line~\ref{alg: HOCUR Line licols}. This can again be done by applying QR decompositions with column pivoting. And finally, the cores of the TT approximation of $\PSI(X)$ are updated in lines~\ref{alg: HOCUR Line core 1} and \ref{alg: HOCUR Line core 2}, by multiplication of parts of the determined CUR decomposition. In the notation used in Section~\ref{sec: HOCUR}, the updated cores are given by tensor foldings of $C \cdot U^{-1}$ and $U^{-1} \cdot R$, respectively.

\section{Analysis of AMUSE and AMUSEt}
\label{sec:analysis_amuset}

The goal of this section is to show that AMUSEt as introduced in Section~\ref{subsec:amuset} produces an empirical matrix representation of the Koopman operator on data-dependent subspaces of the tensor spaces $\mathbb{V}, \mathbb{W}$, and to establish the convergence of this representation in the limit of infinite data. We first investigate the standard case in Section~\ref{subsec:spectral_subspace}. We define the spectral subspaces and their empirical counterparts. If carried out at fixed prescribed SVD rank $r$, the standard AMUSE algorithm amounts to approximating the Koopman operator on empirical spectral subspaces. Convergence of these spaces and the corresponding Koopman operator representation is then established. This analysis can be carried over to the tensor case. To this end, we first provide a multi-linear analogue of (empirical) spectral subspaces in Section~\ref{subsec:mult_spectral_convergence}, and also establish convergence in the limit of infinite data. In Section~\ref{subsec:mult_spectral_computation}, we then show that, analogous to the standard AMUSE algorithm, AMUSEt is indeed the algorithmic framework to compute projections of evolution operators on  multi-linear spectral subspaces.

In what follows, the space of bounded linear operators on a Hilbert space $\mathbb{H}$, equipped with the standard operator norm, is denoted by $L(\mathbb{H})$. The orthogonal projector onto a finite-dimensional subspace $\mathbb{V}$ is labeled $\mathcal{P}_\mathbb{V}$. The \emph{distance between subspaces} $\mathbb{V}$ and $\mathbb{W}$ of the same dimension $n$ is given by
\begin{equation*}
d(\mathbb{V}, \mathbb{W}) = \sup_{v \in \mathbb{V}, \|v\| = 1} \inf_{w \in \mathbb{W}} \|w - v \| = \|(\mathrm{Id} - \mathcal{P}_\mathbb{W}) \mathcal{P}_\mathbb{V} \|_{L(\mathbb{H})} = \|\mathcal{P}_\mathbb{V} - \mathcal{P}_\mathbb{W} ||_{L(\mathbb{H})}.
\end{equation*}
Finally, consider the situation that $\mathbb{V}$ is a finite-dimensional Hilbert space of functions on $\mathbb{R}^d$, and $x_k$, $ k\in \mathbb{N}$, is a sequence of $\mathbb{R}^d$-valued random variables such that the first equation in \eqref{eq:ergodicity_Cpsi} holds true almost surely. Then for $m$ large enough, the bilinear form
\begin{align*}
\innerprod{\psi}{\widetilde{\psi}}_{\mathbb{V}}^{\wedge} &:= \frac{1}{m}\sum_{k=1}^m \psi(x_k) \widetilde{\psi}(x_k), \quad \psi, \, \widetilde{\psi} \in \mathbb{V},
\end{align*}
is an inner product on $\mathbb{V}$. Orthogonal projections with respect to this empirical inner product will be labeled $\widehat{\mathcal{P}}_{\cdot}$, accordingly.

\subsection{Spectral Subspaces and AMUSE}
\label{subsec:spectral_subspace}

We begin this section by explicitly stating the definition of the empirical counterpart of the projected Koopman operator $\mathcal{K}_\tau(\mathbb{V}, \mathbb{W})$, independently of the basis sets used for $\mathbb{V}, \mathbb{W}$, see \cite{Korda2018}. Analogous results can be obtained for the Perron--Frobenius and forward-backward operators:

\begin{Proposition}
\label{prop:empirical_koopman}
Let \eqref{eq:strong_law_sampling} hold, and let $\mathbb{V}, \, \mathbb{W}$ be finite-dimensional subspaces of $L^2_{\rho_0}$ and $L^2_{\rho_1}$. For almost all sequences $(x_k, y_k)$, and $m$ large enough, there is a linear operator $\widehat{\mathcal{K}}_\tau(\mathbb{V}, \mathbb{W}) \colon \mathbb{W} \to \mathbb{V}$, satisfying
\begin{align}
\label{eq:def_empirical_Koopman}
\innerprod{\psi}{\widehat{\mathcal{K}}_\tau(\mathbb{V}, \mathbb{W})\phi}_{L^2_{\rho_0}}^{\wedge} &= \frac{1}{m}\sum_{k=1}^m \psi(x_k) \phi(y_k) \quad \forall \psi \in \mathbb{V}.
\end{align}
Its matrix representation with respect to $\psi, \phi$ is $\widehat{K}_\tau(\psi, \phi) = (\widehat{C}(\psi))^{-1}\widehat{A}(\psi, \phi)$. Also, 
\[\|\widehat{K}_\tau(\mathbb{V}, \mathbb{W}) - \mathcal{K}_\tau(\mathbb{V}, \mathbb{W}) \|_{L(\mathbb{W}, \mathbb{V})} \rightarrow 0\]
for $m \rightarrow \infty$, where the topology on $L(\mathbb{W}, \mathbb{V})$ is induced by the standard inner products on $L^2_{\rho_0}, \, L^2_{\rho_1}$.
\end{Proposition}

\begin{proof}
The right-hand side of \eqref{eq:def_empirical_Koopman} is a linear functional on the finite-dimensional space $\mathbb{V}$ with empirical inner product, which ensures existence of the operator $\widehat{K}_\tau(\mathbb{V}, \mathbb{W})$. The matrix representation can be directly verified, and by \eqref{eq:strong_law_sampling}, this representation converges to that of $\mathcal{K}_\tau(\mathbb{V}, \mathbb{W})$ in any matrix norm, which proves the last statement.
\end{proof}

We also observe that for subspaces $\mathbb{F}_1 \subset \mathbb{V}$ and $\mathbb{F}_2 \subset \mathbb{W}$, we have $\widehat{K}_\tau(\mathbb{F}_1, \mathbb{F}_2) = \widehat{\mathcal{P}}_{\mathbb{F}_1} \widehat{K}_\tau(\mathbb{V}, \mathbb{W}) \widehat{\mathcal{P}}_{\mathbb{F}_2}$. Next, we introduce the family of subspaces which serve as reduced trial spaces for the Koopman operator if the standard AMUSE algorithm is employed. We have briefly encountered these spaces before in Section~\ref{subsec:galerkin_proj}.

\begin{Definition}[Spectral Subspace]
\label{def:spectral_subspaces}
Let $\mathbb{V}$ be a finite-dimensional Hilbert space of functions with basis $\psi = \{\psi_j \}_{j=1}^n$. Denote the spectral decomposition of the Gramian matrix by $C(\psi) = U_\psi \Sigma^2_\psi U_\psi^\top$, with eigenvalues arranged in decreasing order. For $r \leq n$ such that $\sigma_r > \sigma_{r+1}$, denote the first $r$ columns of $U_\psi$ by $U_{\psi, r}$, and the upper $r\times r$-block of the the diagonal matrix by $\Sigma_{\psi, r}$. Then the space $\mathbb{G}_r \subset \mathbb{V}$ spanned by orthonormal functions $\eta^\top_r = \psi^\top (U_{\psi, r} \Sigma_{\psi, r}^{-1})$ is called \emph{spectral subspace (of order $r$)} of $\mathbb{V}$. The coefficient vector space associated with $\mathbb{G}_r$ is denoted by $\mathbb{B}_r \subset \mathbb{R}^n$.
\end{Definition}

\begin{Remark} ~
\label{rem:properties_spectral_subspace}
\begin{itemize}
\item[(i)] The condition $\sigma_r > \sigma_{r+1}$ is necessary for $\mathbb{G}_r$ to be well-defined, otherwise it would be unclear how to break up the eigenspace corresponding to $\sigma_r = \sigma_{r+1}$.
\item[(ii)] The spectral subspace depends on the basis chosen for the Hilbert space $\mathbb{V}$. We will emphasize this dependence by writing $\mathbb{G}_r(\psi), \, \mathbb{B}_r(\psi)$ whenever necessary.
\item[(iii)] If an orthogonal change of basis is used, that is, the basis $\psi$ changes to $\widetilde{\psi}^\top = \psi^\top Q$, with $Q^\top Q = \mathrm{Id}_n$, then the Gramian matrix changes to $C(\widetilde{\psi} ) = Q^\top C(\psi) Q$. Hence, the spectral subspaces $\mathbb{G}_r$ remain the same, the associated coefficient vector spaces are $\mathbb{B}_r = \mathrm{span}(Q^\top U_{\psi, r})$.
\item[(iv)] Clearly, empirical spectral subspaces $\widehat{\mathbb{G}}_r(\psi), \, \widehat{\mathbb{B}}_r(\psi)$ can be defined in the same way using the empirical Gramian matrix.
\end{itemize}
\end{Remark}

Empirical spectral subspaces are stable in the limit of infinite data:

\begin{Lemma}
\label{lem:consistency_amuse_p1}
Let $\mathbb{V}$ be a finite-dimensional Hilbert space of functions on $\mathbb{R}^d$ with basis $\psi = \{\psi_j \}_{j=1}^n$. For fixed $r \leq n$, let the spectral subspace $\mathbb{G}_r$ be well-defined. Let $x_k$, $k\in \mathbb{N}$, be a sequence of $\mathbb{R}^d$-valued random variables such that the first equation in \eqref{eq:ergodicity_Cpsi} holds true almost surely. Then we also have with probability one:
\begin{align*}
\lim_{m\rightarrow \infty} d(\mathbb{B}_r, \widehat{\mathbb{B}}_r) &= 0, & \lim_{m\rightarrow \infty} d(\mathbb{G}_r, \widehat{\mathbb{G}}_r) &= 0, & \lim_{m\rightarrow \infty} \|\widehat{\mathcal{P}}_{\widehat{\mathbb{G}}_r(\psi)} - \mathcal{P}_{\mathbb{G}_r(\psi)} \|_{L(\mathbb{V})} &= 0,
\end{align*}
where the topology on $L(\mathbb{V})$ is induced by the standard inner product on $\mathbb{V}$.
\end{Lemma}

\begin{proof}
By assumption, we have $\lim_{m\rightarrow \infty} \| C(\psi) - \widehat{C}(\psi) \|_2 = 0$ almost surely, hence $\widehat{\mathbb{B}}_r,\, \widehat{\mathbb{G}}_r$ can be defined with probability one if $m$ is large enough. As the gap between $\sigma_r$ and $\sigma_{r+1}$ is positive, perturbation theory for symmetric matrices \cite{STEWART1990}[Chapter 5, Thm 3.6] ensures that the distance in $\mathbb{R}^n$ between $\mathbb{B}_r$ and $\widehat{\mathbb{B}}_r$ converges to zero. As these are the coefficient vector spaces corresponding to $\mathbb{G}_r, \, \widehat{\mathbb{G}}_r$, Lemma \ref{lem:equivalences_subspaces} (iii) in Appendix \ref{app:aux_results} implies the second claim. Finally, we choose an orthonormal basis (ONB) $U \in \mathbb{R}^{n \times r}$ of $\mathbb{B}_r(\psi)$, and invoke Lemma \ref{lem:convergence_procrustes} to choose ONBs $\widehat{U}$ of $\widehat{\mathbb{B}}_r(\psi)$ such that $\| \widehat{U} - U \|_2 \rightarrow 0$ holds almost surely. Using the representations in Lemma \ref{lem:equivalences_subspaces}(i) and (ii), this implies
\begin{align*}
& \quad \|\widehat{\mathcal{P}}_{\widehat{\mathbb{G}}_r(\psi)} - \mathcal{P}_{\mathbb{G}_r(\psi)} \|_{L(\mathbb{V})} \\
&= \|C(\psi)^{1/2}\left[\widehat{U}(\widehat{U}^\top \widehat{C}(\psi) \widehat{U})^{-1} \widehat{U}^\top \widehat{C}(\psi) - U(U^\top C(\psi) U)^{-1} U^\top C(\psi)\right] C(\psi)^{-1/2} \|_2 \rightarrow 0,
\end{align*}
using that $\widehat{C}(\psi) \rightarrow C(\psi)$ and $\widehat{U} \rightarrow U$.
\end{proof}

We can now use the previous result to complete our analysis of AMUSE. We have already seen in Section~\ref{subsec:amuse} that, if the SVD rank $r$ is fixed, the basis sets $\widehat{\eta}_r, \, \widehat{\zeta}_r$ appearing in the standard AMUSE algorithm are spanning the empirical spectral subspaces $\widehat{\mathbb{G}}_r(\psi), \, \widehat{\mathbb{G}}_r(\phi)$ in $\mathbb{V}$ and $\mathbb{W}$. We show that the resulting empirical projection of the Koopman operator converges in concert with these spaces.

\begin{Proposition}
\label{prop:stability_rt}
Let $\mathcal{X}_t$ be a dynamical system such that \eqref{eq:strong_law_sampling} holds true. Let $\mathbb{V} \subset L^2_{\rho_0}$ and $\mathbb{W} \subset L^2_{\rho_1}$ be $n$-dimensional subspaces with bases $\psi, \, \phi$. For $r \leq n$ such that $\mathbb{G}_r(\psi)$ and $\mathbb{G}_r(\phi)$ are both well-defined, we conclude that almost surely:
\[ \|\widehat{\mathcal{K}}_\tau(\widehat{\mathbb{G}}_r(\psi), \widehat{\mathbb{G}}_r(\phi)) - \mathcal{K}_\tau(\mathbb{G}_r(\psi), \mathbb{G}_r(\phi)) \|_{L(\mathbb{W}, \mathbb{V})} \rightarrow 0. \]
\end{Proposition}
\begin{proof}
We have already seen that 
\[\widehat{\mathcal{K}}_\tau(\widehat{\mathbb{G}}_r(\psi), \widehat{\mathbb{G}}_r(\phi)) = \widehat{\mathcal{P}}_{\widehat{\mathbb{G}}_r(\psi)} \widehat{\mathcal{K}}_\tau(\mathbb{V}, \mathbb{W}) \widehat{\mathcal{P}}_{\widehat{\mathbb{G}}_r(\phi)}, \quad \mathcal{K}_\tau(\mathbb{G}_r(\psi), \mathbb{G}_r(\phi)) = \mathcal{P}_{\mathbb{G}_r(\psi)} \mathcal{K}_\tau(\mathbb{V}, \mathbb{W}) \mathcal{P}_{\mathbb{G}_r(\phi)}. \]
Since $\widehat{\mathcal{P}}_{\widehat{\mathbb{G}}_r(\psi)} \rightarrow \mathcal{P}_{\mathbb{G}_r(\psi)}$ and $\widehat{\mathcal{P}}_{\widehat{\mathbb{G}}_r(\phi)} \rightarrow \mathcal{P}_{\mathbb{G}_r(\phi)}$ by Lemma \ref{lem:consistency_amuse_p1}, and $\widehat{\mathcal{K}}_\tau(\mathbb{V}, \mathbb{W}) \rightarrow \mathcal{K}_\tau(\mathbb{V}, \mathbb{W})$ by Proposition \ref{prop:empirical_koopman}, the claim follows.
\end{proof}

\subsection{Multi-linear Spectral Subspaces}
\label{subsec:mult_spectral_convergence}
For the rest of this section, we assume that a vector of fixed ranks $\mathbf{r} = [r_1, \ldots, r_p]$ is given. Also, $\mathbb{V} = \bigotimes_{k=1}^p \mathbb{V}^k$ is a tensor product space of functions on $\mathbb{R}^d$, where $\mathbb{V}^k = \mathrm{span}\{\psi_{k, i_k}\}_{i_k = 1}^{n_k}$. The full tensor product basis is denoted by $\PSI$. We also introduce the symbols $\mathbb{V}^{:k} = \bigotimes_{l=1}^k \mathbb{V}^l$ for the partial tensor product up to mode $k$, the corresponding basis is denoted $\PSI_{:k}$. Similarly, if $\mathbf{T}$ is a tensor train of order $p$, the partial tensor train up to mode $k \leq p$ is denoted $\mathbf{T}^{:k} = \llbracket \mathbf{T}^{(1)} \rrbracket \otimes \ldots \otimes \llbracket \mathbf{T}^{(k)} \rrbracket$.

The results of the previous section can now be generalized to tensor-structured subspaces. The construction of multi-linear spectral subspaces, to be described below, is inspired by the global SVD. Recall from Section~\ref{subsec:global_svd} that the first step of Algorithm~\ref{alg:SVD} applied to $\PSI(X)$ and $\PHI(Y)$ is the same as if we were applying standard AMUSE just to $\mathbb{V}^1$ and $\mathbb{W}^1$. For $k = 1$, and focusing just on $\PSI(X)$ for the sake of illustration, the matrix of singular vectors $U$ in line \ref{algline:svd_step} of Algorithm~\ref{alg:SVD} encodes a basis of the $r_1$-dimensional spectral subspace $\widehat{\mathbb{G}}_{r_1}(\psi_1)$. As we show in Section~\ref{subsec:mult_spectral_computation}, the next step amounts to computing a basis of the $r_2$-dimensional spectral subspace for a specific basis of the tensor product $\widehat{\mathbb{G}}_{r_1}(\psi_1) \otimes \mathbb{V}^2$. This procedure is then repeated through all steps of the method. Some care, however, needs to be taken when choosing the basis set for each of these spectral subspaces.

\begin{Definition}[Multi-linear Spectral Subspaces]
\label{def:rank_r_subspaces}
Define $\mathbb{G}^{:0} = \mathrm{span}\{1\}$ and $\theta_0 = 1$. Recursively for $k = 1, \ldots, p$, consider the Gramian matrix $C(\mathbf{vec}(\theta_{k-1} \otimes \psi_k))$ of $\mathbb{G}^{:k-1} \otimes \mathbb{V}^k$, and denote its eigenvalues by $\sigma^2_{k, l},\, 1 \leq l \leq r_{k-1} n_k$. If $\sigma_{k, r_k} > \sigma_{k, r_k + 1}$, denote its $r_k$-dimensional spectral subspace by $\mathbb{G}^{:k}$. Choose any orthonormal basis $U_k \in \mathbb{R}^{r_{k-1} n_k \times r_k}$ of the associated coefficient vector space in $\mathbb{R}^{r_{k-1}n_k}$, and define a basis set for $\mathbb{G}^{:k}$ as $(\theta_k)^\top = (\mathbf{vec}(\theta_{k-1} \otimes \psi_k))^\top U_k$.

The subspaces $\mathbb{G}^{:k}$ are then called \emph{multi-linear spectral subspaces}. Their associated coefficient vector spaces in $\mathbb{R}^{\prod_{l=1}^k n_l}$, with respect to the basis $\PSI_{:k}$, are denoted $\mathbb{B}^{:k}$.
\end{Definition}

\begin{Remark}~
\label{rem:empirical_gramian_subspace}
\begin{itemize}
 \item[(i)] Definition \ref{def:rank_r_subspaces} can be repeated almost verbatim to define empirical counterparts of all quantities introduced above. As before, we will use hats $\,\widehat{\cdot}\,$ to denote these quantities.
 \item[(ii)] Note that the ONBs $U_k$ used in the construction are coefficient vectors with respect to the previous level of recursion, that is, to a basis of $\mathbb{G}^{:k-1} \otimes \mathbb{V}^k$. The coefficient vector space $\mathbb{B}^{:k}$ refers to the full $k$-fold tensor product basis $\PSI_{:k}$.
\end{itemize}
\end{Remark}

The spaces $\mathbb{G}^{:k}$ are independent of the specific orthonormal bases used in each step of the construction. In fact, the following results can be proven:

\begin{Lemma}
\label{lem:recursive_relations_gramian_subspaces}
Compile the orthonormal bases $U_k$ in the construction of Definition~\ref{def:rank_r_subspaces} into a tensor train $\mathbf{U}$, that is:
\begin{align*}
\mathbf{U} &= \llbracket \mathbf{U}^{(1)} \rrbracket \otimes \ldots \otimes \llbracket \mathbf{U}^{(p)} \rrbracket, & \mathbf{U}^{(k)}\vert_2 &= U_k, 1\leq k \leq p.
\end{align*}
For $1 \leq k \leq p$, the partial tensor train $\mathbf{U}^{:k} \in \mathbb{R}^{n_1 \times \ldots \times n_k \times r_k}$ is the coefficient tensor of the basis $\theta_k$ with respect to $\PSI_{:k}$, that is $(\theta_k)^\top = (\PSI_{:k}\vert_k)^\top \mathbf{U}^{:k}\vert_k$.
Moreover, let $\widetilde{U}_k$, $1 \leq k \leq p$  be a different sequence of ONBs used in the construction in Definition~\ref{def:rank_r_subspaces}, with resulting tensor train $\widetilde{\mathbf{U}}$ and basis sets $\widetilde{\theta}_k$. Then there exist orthonormal matrices $Q_k \in \mathbb{R}^{r_k \times r_k}, \, 0\leq k \leq p$, such that
\begin{equation}
\label{eq:rel_cores_theta_basis_sets}
\widetilde{\mathbf{U}}^{(k)}\vert_2 = (Q_{k-1} \otimes \mathrm{Id}_{n_k})^\top \mathbf{U}^{(k)}\vert_2 Q_k
\end{equation}
for $1 \leq k \leq p$. In particular, $\theta_k$ and $\widetilde{\theta}_k$ are related by
\begin{align}
\label{eq:rel_theta_basis_sets}
(\theta_k)^\top &= (\PSI_{:k}\vert_k)^\top \mathbf{U}^{:k}\vert_k, &
(\widetilde{\theta}_k)^\top &= (\PSI_{:k}\vert_k)^\top \mathbf{U}^{:k}\vert_k Q_k = (\theta_k)^\top Q_k,
\end{align}
and the spaces $\mathbb{G}^{:k}$ and $\mathbb{B}^{:k}$ are the same in both cases.
\end{Lemma}

\begin{proof}
The fact that each $\mathbf{U}^{:k}$ is the coefficient tensor of the basis $\theta_k$ can be verified directly by an inductive argument. To prove (\ref{eq:rel_cores_theta_basis_sets}--\ref{eq:rel_theta_basis_sets}), we proceed by induction over $k$. The claims are clearly true for $k = 1$ with $Q_0 = 1$. For general $k$, if the claims of the Lemma are true for $k - 1$, then \eqref{eq:rel_theta_basis_sets} for $k-1$, and Remark~\ref{rem:properties_spectral_subspace} (iii), imply that the spectral subspaces of $C(\mathbf{vec}(\theta_{k-1} \otimes \psi_k))$ and $C(\mathbf{vec}(\widetilde{\theta}_{k-1} \otimes \psi_k))$  are the same, and the coefficient vector spaces are related by the transformation $(Q_{k-1} \otimes \mathrm{Id}_{n_k})^\top$.  It follows that there is an orthonormal $Q_k \in \mathbb{R}^{r_k \times r_k}$ such that
\[ \widetilde{U}_k = (Q_{k-1} \otimes \mathrm{Id}_{n_k})^\top U_k Q_k,\]
as claimed. In turn, this implies that
\begin{align*}
\widetilde{\mathbf{U}}^{:k}\vert_k &= \left(\widetilde{\mathbf{U}}^{:k-1}\vert_{k-1} \otimes \mathrm{Id}_{n_k}\right) \widetilde{U}_k \\
&= \left( (\mathbf{U}^{:k-1}\vert_{k-1} Q_{k-1}) \otimes \mathrm{Id}_{n_k} \right) (Q_{k-1} \otimes \mathrm{Id}_{n_k})^\top U_k Q_k \\
&= \left(\mathbf{U}^{:k-1}\vert_{k-1} \otimes \mathrm{Id}_{n_k} \right) U_k Q_k = \mathbf{U}^{:k}\vert_k Q_k.
\end{align*}
This proves \eqref{eq:rel_theta_basis_sets}, and hence $\theta_k$ and $\widetilde{\theta}_k$ span the same space $\mathbb{G}^{:k}$ with coefficient vector space $\mathbb{B}^{:k}$.
\end{proof}

\begin{figure}[htbp]
    \centering
    \begin{tikzpicture}
        
        \def\x{1}
        \def\y{2.8}
        \def\d{0.8}
        \node[anchor=east] at (\x-0.25,\y-0.25) {$(\theta_k(x))^\top$}; 
        \node[] at (\x,\y-0.25) {$=$};
        \def\x{1.65}
        \draw[black] (\x,\y) -- node [label={[shift={(0,-0.15)}]$r_1$}] {} ++ (\d,0) ;
        \draw[black] (\x+\d,\y) -- ++ (0.66*\d,0) ;
        \draw[black, dotted] (\x+1.66*\d,\y) -- ++ (0.66*\d,0) ;
        \draw[black] (\x+2.33*\d,\y) -- ++ (0.66*\d,0) ;

        \draw[black] (\x,\y) -- ++ (0,-0.5) ;
        \draw[black] (\x+\d,\y) -- ++ (0,-0.5) ;
        \draw[black] (\x+3*\d,\y) -- ++ (0,-0.5) ;
        \draw[black] (\x+3*\d,\y) -- node [label={[shift={(0,-0.15)}]$r_k$}] {} ++ (0.75,0) ;
        \node[draw,shape=semicircle,rotate=135,fill=white, anchor=south,inner sep=2pt, outer sep=0pt, scale=0.75] at (\x,\y){}; 
        \node[draw,shape=semicircle,rotate=315,fill=Blue, anchor=south,inner sep=2pt, outer sep=0pt, scale=0.75] at (\x,\y){};
        \node[draw,shape=circle,fill=Green, scale=0.65] at (\x,\y-0.5){};
        \node[draw,shape=semicircle,rotate=135,fill=white, anchor=south,inner sep=2pt, outer sep=0pt, scale=0.75] at (\x+\d,\y){}; 
        \node[draw,shape=semicircle,rotate=315,fill=Blue, anchor=south,inner sep=2pt, outer sep=0pt, scale=0.75] at (\x+\d,\y){};
        \node[draw,shape=circle,fill=Green, scale=0.65] at (\x+\d,\y-0.5){};
        \node[draw,shape=semicircle,rotate=135,fill=white, anchor=south,inner sep=2pt, outer sep=0pt, scale=0.75] at (\x+3*\d,\y){}; 
        \node[draw,shape=semicircle,rotate=315,fill=Blue, anchor=south,inner sep=2pt, outer sep=0pt, scale=0.75] at (\x+3*\d,\y){};
        \node[draw,shape=circle,fill=Green, scale=0.65] at (\x+3*\d,\y-0.5){};
        
        \def\x{5.25}
        \node[] at (\x,\y-0.25) {$=$};
        
        \def\x{6.7}
        \draw[black] (\x-0.8,\y+0.1) -- ++ (0,-1) ;
        \draw[black] (\x-0.3,\y+0.1) -- ++ (0,-1) ;
        \node[] at (\x+0.3,\y-0.3){$\cdots$};
        \draw[black] (\x+0.8,\y+0.1) -- ++ (0,-1) ;
        \draw[black] (\x,\y+0.1) -- ++ (0,0.75) ;
        \node[draw, inner sep=5pt, fill=white, minimum width=2cm] at (\x,\y+0.25) {$\mathbf{U}^{:k} $};
        \node[draw, inner sep=5pt, fill=white, minimum width=2cm] at (\x,\y-0.75) {$\mathbf{\Psi}_{:k}$};
        
        \def\x{8.15}
        \node[] at (\x,\y-0.25) {$=$};
        \def\x{9.4}
        \draw[black] (\x+1,\y) -- ++ (1,0) ;
        \draw[black] (\x+2,\y) -- node [label={[shift={(0,-0.15)}]$r_k$}] {} ++ (0.75,0) ;
        \draw[black] (\x+2,\y) -- ++ (0,-0.5) ;
        \node[draw, inner sep=5pt, fill=white, anchor=east] at (\x+1.25,\y-0.25) {$(\theta_{k-1}(x))^\top$};
        \node[draw,shape=semicircle,rotate=135,fill=white, anchor=south,inner sep=2pt, outer sep=0pt, scale=0.75] at (\x+2,\y){}; 
        \node[draw,shape=semicircle,rotate=315,fill=Blue, anchor=south,inner sep=2pt, outer sep=0pt, scale=0.75] at (\x+2,\y){};
        \node[draw,shape=circle,fill=Green, scale=0.65] at (\x+2,\y-0.5){};
        
        
        \def\x{1}
        \def\y{0.7}
        \def\d{0.625}
        \node[anchor=east] at (\x-0.25,\y-0.25) {$(\widetilde{\theta}_k(x))^\top$}; 
        \node[] at (\x,\y-0.25) {$=$};
        \def\x{1.65}
        \draw[black] (\x,\y) -- ++ (5*\d,0) ;
        \draw[black, dotted] (\x+5*\d,\y) -- ++ (\d,0) ;
        \draw[black] (\x+6*\d,\y) -- ++ (3*\d,0) ;
        \draw[black] (\x+9*\d,\y) -- node [label={[shift={(0,-0.15)}]$r_k$}] {} ++ (0.75,0) ;
        \draw[black] (\x,\y) -- ++ (0,-0.5) ;
        \draw[black] (\x+3*\d,\y) -- ++ (0,-0.5) ;
        \draw[black] (\x+8*\d,\y) -- ++ (0,-0.5) ;
        \node[draw,shape=semicircle,rotate=135,fill=white, anchor=south,inner sep=2pt, outer sep=0pt, scale=0.75] at (\x,\y){}; 
        \node[draw,shape=semicircle,rotate=315,fill=Blue, anchor=south,inner sep=2pt, outer sep=0pt, scale=0.75] at (\x,\y){};
        \node[draw,shape=circle,fill=Green, scale=0.65] at (\x,\y-0.5){};
        \node[draw,shape=semicircle,rotate=90,fill=Gray, anchor=south,inner sep=2pt, outer sep=0pt, scale=0.75] at (\x+\d,\y){}; 
        \node[draw,shape=semicircle,rotate=270,fill=white, anchor=south,inner sep=2pt, outer sep=0pt, scale=0.75] at (\x+\d,\y){};
        \node[draw,shape=semicircle,rotate=90,fill=white, anchor=south,inner sep=2pt, outer sep=0pt, scale=0.75] at (\x+2*\d,\y){}; 
        \node[draw,shape=semicircle,rotate=270,fill=Gray, anchor=south,inner sep=2pt, outer sep=0pt, scale=0.75] at (\x+2*\d,\y){};
        \node[draw,shape=semicircle,rotate=135,fill=white, anchor=south,inner sep=2pt, outer sep=0pt, scale=0.75] at (\x+3*\d,\y){}; 
        \node[draw,shape=semicircle,rotate=315,fill=Blue, anchor=south,inner sep=2pt, outer sep=0pt, scale=0.75] at (\x+3*\d,\y){};
        \node[draw,shape=circle,fill=Green, scale=0.65] at (\x+3*\d,\y-0.5){};
        \node[draw,shape=semicircle,rotate=90,fill=Gray, anchor=south,inner sep=2pt, outer sep=0pt, scale=0.75] at (\x+4*\d,\y){}; 
        \node[draw,shape=semicircle,rotate=270,fill=white, anchor=south,inner sep=2pt, outer sep=0pt, scale=0.75] at (\x+4*\d,\y){};
        \node[draw,shape=semicircle,rotate=135,fill=white, anchor=south,inner sep=2pt, outer sep=0pt, scale=0.75] at (\x+8*\d,\y){}; 
        \node[draw,shape=semicircle,rotate=315,fill=Blue, anchor=south,inner sep=2pt, outer sep=0pt, scale=0.75] at (\x+8*\d,\y){};
        \node[draw,shape=circle,fill=Green, scale=0.65] at (\x+8*\d,\y-0.5){};
        \node[draw,shape=semicircle,rotate=90,fill=white, anchor=south,inner sep=2pt, outer sep=0pt, scale=0.75] at (\x+7*\d,\y){}; 
        \node[draw,shape=semicircle,rotate=270,fill=Gray, anchor=south,inner sep=2pt, outer sep=0pt, scale=0.75] at (\x+7*\d,\y){};
        \node[draw,shape=semicircle,rotate=90,fill=Gray, anchor=south,inner sep=2pt, outer sep=0pt, scale=0.75] at (\x+9*\d,\y){}; 
        \node[draw,shape=semicircle,rotate=270,fill=white, anchor=south,inner sep=2pt, outer sep=0pt, scale=0.75] at (\x+9*\d,\y){};
        \def\x{8.5}
        \node[] at (\x,\y-0.25) {$=$};
        \def\x{9.4}
        \draw[black] (\x+1,\y) -- ++ (1,0) ;
        \draw[black] (\x+2,\y) -- node [label={[shift={(0,-0.15)}]$r_k$}] {} ++ (0.75,0) ;
        \node[draw, inner sep=5pt, fill=white, anchor=east] at (\x+1.25,\y-0.25) {$(\theta_{k}(x))^\top$};
        \node[draw,shape=semicircle,rotate=90,fill=Gray, anchor=south,inner sep=2pt, outer sep=0pt, scale=0.75] at (\x+2,\y){}; 
        \node[draw,shape=semicircle,rotate=270,fill=white, anchor=south,inner sep=2pt, outer sep=0pt, scale=0.75] at (\x+2,\y){};

        
        \def\x{1}
        \def\y{-1.35}
        \node[anchor=east] at (\x-0.25,\y-0.25) {$\theta_{k-1} \otimes \psi_k$}; 
        \node[] at (\x,\y-0.25) {$=$};
        \def\x{1.65}
        \def\d{0.92}
        \draw[black] (\x,\y+0.1) -- ++ (0,-1) ;
        \draw[black] (\x+\d,\y+0.1) -- ++ (0,-1) ;
        \draw[black] (\x+2*\d,\y+0.1) -- ++ (0,-1) ;
        \draw[black] (\x+4*\d,\y+0.1) -- ++ (0,-1) ;
        \draw[black] (\x+6*\d,\y+0.1) -- ++ (0,-1) ;
        \draw[black] (\x+2*\d,\y+0.1) -- ++ (0,0.75) ;
        \draw[black] (\x+6*\d,\y+0.1) -- ++ (0,0.75) ;
        \node[draw, inner sep=5pt, fill=white, minimum width=4cm] at (\x+2*\d,\y+0.25) {$\mathbf{U}^{:k-1}$};
        \node[draw,shape=circle,fill=Green, scale=0.65] at (\x,\y-0.95){};
        \node[draw,shape=circle,fill=Green, scale=0.65] at (\x+\d,\y-0.95){};
        \node[draw,shape=circle,fill=Green, scale=0.65] at (\x+2*\d,\y-0.95){};
        \node[] at (\x+3*\d,\y-0.95){$\cdots$};
        \node[draw,shape=circle,fill=Green, scale=0.65] at (\x+4*\d,\y-0.95){};
        \node[] at (\x+5*\d,\y-0.25) {$\otimes$};
        \node[draw, inner sep=5pt, fill=white] at (\x+6*\d,\y+0.25) {$\mathrm{Id}$};
        \node[draw,shape=circle,fill=Green, scale=0.65] at (\x+6*\d,\y-0.95){};
        \def\x{8}
        \node[] at (\x,\y-0.25) {$=$};
        \def\x{8.95}
        \def\d{0.91}
        \draw[black] (\x,\y+0.1) -- ++ (0,-1) ;
        \draw[black] (\x+\d,\y+0.1) -- ++ (0,-1) ;
        \node[] at (\x+2*\d,\y-0.3){$\cdots$};
        \draw[black] (\x+3*\d,\y+0.1) -- ++ (0,-1) ;
        \draw[black] (\x+1*\d,\y+0.1) -- ++ (0,0.75) ;
        \draw[black] (\x+2*\d,\y+0.1) -- ++ (0,0.75) ;
        \node[draw, inner sep=5pt, fill=white, minimum width=3.6cm] at (\x+1.5*\d,\y+0.25) {$\mathbf{U}^{:k-1} \otimes \mathrm{Id}$};
        \node[draw, inner sep=5pt, fill=white, minimum width=3.6cm] at (\x+1.5*\d,\y-0.75) {$\mathbf{\Psi}_{:k}$};
    
    \end{tikzpicture}
    \caption{Relationship between the basis sets $\theta_k$ and $\widetilde{\theta}_k$ for the construction of multi-linear spectral subspaces. Half-filled blue circles represent the cores of $\mathbf{U}$, green circles the corresponding tensor product bases, and half-filled gray circles the orthonormal matrices $Q_k$. The third row illustrates the recursive expression of the product basis $\theta_{k-1} \otimes \psi_k$ with respect to the coefficient tensor $\mathbf{U}^{:k-1}$ of $\theta_{k-1}$.
    \label{fig:appendix_networks}}
\end{figure}

A pictorial illustration of the relationship between $\theta_k$ and $\widetilde{\theta}_k$ is given in the first two rows of Figure~\ref{fig:appendix_networks}. Now that multi-linear spectral subspaces have been shown to be well-defined, it is time to generalize Lemma~\ref{lem:consistency_amuse_p1} to the multi-linear case:

\begin{Lemma}
\label{lem:conv_spectral_subspace}
Let $\mathbb{V} = \bigotimes_{k=1}^p \mathbb{V}^k$ be a finite-dimensional tensor space of functions on $\mathbb{R}^d$ as described above. For a fixed vector of TT ranks $\mathbf{r}$, let all multi-linear spectral subspaces $\mathbb{G}^{:k}$ be well-defined. Let $x_k, \, k\in \mathbb{N}$ be a sequence of $\mathbb{R}^d$-valued random variables such that the first equation in \eqref{eq:ergodicity_Cpsi} holds true almost surely in $\mathbb{V}$. Then we have for all $1 \leq k \leq p$:
\begin{align*}
d(\mathbb{B}^{:k}, \widehat{\mathbb{B}}^{:k}) &\rightarrow 0, & d(\mathbb{G}^{:k}, \widehat{\mathbb{G}}^{:k}) &\rightarrow 0, & \|\widehat{\mathcal{P}}_{\widehat{\mathbb{G}}^{:k}} -  \mathcal{P}_{\mathbb{G}^{:k}}\|_{L(\mathbb{V})} &\rightarrow 0.
\end{align*}
\end{Lemma}

\begin{proof}
The proof is an inductive application of Lemma \ref{lem:consistency_amuse_p1}. For $k = 1$, the result is the same as that of Lemma \ref{lem:consistency_amuse_p1}, so let us assume the claim stands for $k -1$. We fix the coefficient tensors $\mathbf{U}^{:k}, \, \widehat{\mathbf{U}}^{:k}$. By Lemma \ref{lem:recursive_relations_gramian_subspaces}, the unfoldings $\mathbf{U}^{:k-1}\vert_{k-1}, \, \widehat{\mathbf{U}}^{:k-1}\vert_{k-1}$ provide ONBs of $\mathbb{B}^{:k-1}, \, \widehat{\mathbb{B}}^{:k-1}$. We use the induction hypothesis and Lemma \ref{lem:convergence_procrustes} to determine an orthonormal $\widehat{Q}_{k-1} \in \mathbb{R}^{r_{k-1} \times r_{k-1}}$ such that $\|\mathbf{U}^{:k-1}\vert_{k-1} - \widehat{\mathbf{U}}^{:k-1}\vert_{k-1} \widehat{Q}_{k-1} \|_2 \rightarrow 0$ almost surely. By extension, we also have
\begin{equation}
\label{eq:convergence_basis_hk}
\|\mathbf{U}^{:k-1}\vert_{k-1} \otimes \mathrm{Id}_{n_k} - (\widehat{\mathbf{U}}^{:k-1}\vert_{k-1} \widehat{Q}_{k-1}) \otimes \mathrm{Id}_{n_k} \|_2  \rightarrow 0.
\end{equation}
Denoting the basis encoded by $(\widehat{\mathbf{U}}^{:k-1}\vert_{k-1} \widehat{Q}_{k-1})$ by  $\bar{\theta}_{k-1}$, the matrices in \eqref{eq:convergence_basis_hk} encode the basis sets $\mathbf{vec}(\theta_{k-1} \otimes \psi_k)$ and $\mathbf{vec}(\bar{\theta}_{k-1} \otimes \psi_k)$ with respect to $\PSI_{:k}$, see the third row in Figure~\ref{fig:appendix_networks}. We then find:
\begin{align*}
& \quad \| C(\mathbf{vec}(\theta_{k-1} \otimes \psi_k)) - \widehat{C}(\mathbf{vec}(\bar{\theta}_{k-1} \otimes \psi_k)) \|_2 \\
&= \| \left(\mathbf{U}^{:k-1}\vert_{k-1} \otimes \mathrm{Id}_{n_k}\right)^\top C(\PSI_{:k}) \left(\mathbf{U}^{:k-1}\vert_{k-1} \otimes \mathrm{Id}_{n_k}\right) \\ 
&\quad - \left((\widehat{\mathbf{U}}^{:k-1}\vert_{k-1} \widehat{Q}_{k-1}) \otimes \mathrm{Id}_{n_k} \right)^\top \widehat{C}(\PSI_{:k}) \left((\widehat{\mathbf{U}}^{:k-1}\vert_{k-1} \widehat{Q}_{k-1}) \otimes \mathrm{Id}_{n_k} \right)\|_2  \rightarrow 0 \quad \text{a.s.},
\end{align*}
where convergence follows from \eqref{eq:convergence_basis_hk} and since $\widehat{C}(\PSI_{:k}) \rightarrow C(\PSI_{:k})$. Now, we note that the dominant $r_k$-dimensional eigenspace of $C(\mathbf{vec}(\theta_{k-1} \otimes \psi_k))$ is spanned by $U_k$ (cf. Definition \ref{def:rank_r_subspaces}), while that of $\widehat{C}(\mathbf{vec}(\bar{\theta}_{k-1} \otimes \psi_k))$ is spanned by $(\widehat{Q}_{k-1} \otimes \mathrm{Id}_{n_k})^\top \widehat{U}_k$. By perturbation theory for singular vectors and Lemma~\ref{lem:convergence_procrustes}, we have found that there is an orthonormal $\widehat{Q}_k \in \mathbb{R}^{r_k \times r_k}$ such that $\|U_k -  (\widehat{Q}_{k-1} \otimes \mathrm{Id}_{n_k})^\top \widehat{U}_k \widehat{Q}_k \|_2 \rightarrow 0$. Hence
\begin{align*}
\widehat{\mathbf{U}}^{:k}\vert_k \widehat{Q}_k &= ((\widehat{\mathbf{U}}^{:k-1}\vert_{k-1} \widehat{Q}_{k-1}) \otimes \mathrm{Id}_{n_k})(\widehat{Q}_{k-1} \otimes \mathrm{Id}_{n_k})^\top \widehat{U}_k \widehat{Q}_k \\
&\rightarrow (\mathbf{U}^{:k-1} \otimes \mathrm{Id}_{n_k}) U_k = \mathbf{U}^{:k}\vert_k.
\end{align*}
This proves that $d(\mathbb{B}^{:k}, \widehat{\mathbb{B}}^{:k}) \rightarrow 0$, and the conclusion about $\mathbb{G}^{:k}, \, \widehat{\mathbb{G}}^{:k}$ follows from Lemma~\ref{lem:equivalences_subspaces} (iii). The final statement about projectors can be shown in the same way as in Lemma~\ref{lem:consistency_amuse_p1}.
\end{proof}

With this, we obtain the central result of this section, which is the convergence of Koopman operator representations on multi-linear subspaces:
\begin{Theorem}
\label{thm:consistency_koopman_ml_spectral}
Let $\mathcal{X}_t$ be a dynamical system such that \eqref{eq:strong_law_sampling} holds true. Let $\mathbb{V} = \bigotimes_{k=1}^p \mathbb{V}^k \subset L^2_{\rho_0}$ and $\mathbb{W} = \bigotimes_{k=1}^p \mathbb{W}^k \subset L^2_{\rho_1}$ be tensor product subspaces with bases $\PSI = \bigotimes_{k=1}^p \psi_k$ and $\PHI = \bigotimes_{k=1}^p \phi_k$. For a fixed vector of TT ranks $\mathbf{r}$, denote the multi-linear spectral subspaces of both bases by $\mathbb{G}^{:k}(\PSI), \mathbb{G}^{:k}(\PHI)$, and assume they are all well-defined. Then we have for all $k \leq p$:
\[ \|\widehat{\mathcal{K}}_\tau(\widehat{\mathbb{G}}^{:k}(\PSI), \widehat{\mathbb{G}}^{:k}(\PHI)) - \mathcal{K}_\tau(\mathbb{G}^{:k}(\PSI), \mathbb{G}^{:k}(\PHI)) \|_{L(\mathbb{W}, \mathbb{V})} \rightarrow 0. \]
\end{Theorem}
\begin{proof}
Since $\widehat{\mathcal{P}}_{\widehat{\mathbb{G}}^{:k}(\PSI)} \rightarrow \mathcal{P}_{\mathbb{G}^{:k}(\PSI)}$ in $L(\mathbb{V})$, and $\widehat{\mathcal{P}}_{\widehat{\mathbb{G}}^{:k}(\PHI)} \rightarrow \mathcal{P}_{\mathbb{G}^{:k}(\PHI)}$ in $L(\mathbb{W})$ by Lemma~\ref{lem:conv_spectral_subspace}, the argument is just the same as in the proof of Proposition~\ref{prop:stability_rt}.
\end{proof}

\subsection{Multi-linear Spectral Subspaces and AMUSEt}
\label{subsec:mult_spectral_computation}

We complete this section by showing that the global SVD applied to the transformed data tensor \eqref{eq: transformed data tensor} is the computational tool to compute bases of multi-linear spectral subspaces, analogous to the standard SVD applied to the data matrix. It will follow that AMUSEt provides a representation of the Koopman operator on multi-linear spectral subspaces.

\begin{Proposition}
\label{lem:global_svd_timeseries}
Let $\mathbb{V} = \bigotimes_{k=1}^p \mathbb{V}^k$ be a finite-dimensional tensor space of functions on $\mathbb{R}^d$, with basis $\PSI = \bigotimes_{k=1}^p \psi_k$. Let $x_1, \ldots, x_m \in \mathbb{R}^d$ be data points such that the empirical spectral subspaces $\widehat{\mathbb{G}}^{:k}$ are well-defined. For prescribed ranks $\mathbf{r} = [r_1, \ldots, r_p]$, denote the global SVD of the transformed data tensor \eqref{eq: transformed data tensor} by $\PSI(X) = \widehat{\mathbf{U}}_{X, r} \widehat{\Sigma}_{X, r} V_{X, r}^\top$, with $r = r_p$. Then, the orthonormal part $\widehat{\mathbf{U}}_{X, r}$ provides orthonormal coefficient tensors for each of the spaces $\widehat{\mathbb{B}}^{:k}$, i.e. $(\widehat{\theta}_k)^\top = (\PSI_{:k}\vert_k)^\top \widehat{\mathbf{U}}_{X, r}^{:k}\vert_k$ is a basis for $\widehat{\mathbb{G}}^{:k}$. The remaining parts of the global SVD contain the time series of the basis for the final space $\widehat{\mathbb{G}}^{:p}$, i.e. $\widehat{\Sigma}_{X, r} \widehat{V}_{X, r}^\top = \widehat{\theta}_p(X) = (\widehat{\mathbf{U}}_{X, r}\vert_p)^\top  \PSI(X)\vert_p$.
\end{Proposition}

\begin{proof}
We prove the statement by showing that after $k \leq p - 1$ iterations in Algorithm~\ref{alg:SVD}, the updated core $\PSI^{(k+1)}(X)$ contains a time series of the product basis $\mathbf{vec}(\widehat{\theta}_{k} \, \otimes \, \psi_{k+1})(X)$, where $\widehat{\theta}_k$ is a basis of $\widehat{\mathbb{G}}^{:k}$ encoded by the partial tensor train $\widehat{\mathbf{U}}_{X, r}^{:k}$. This is clearly true for $k = 0$, as $\PSI^{(1)}(X) = \psi_1(X)$ and $\widehat{\theta}_0 = 1$, so we can again resort to an inductive argument. Assume that, after $k - 1$ iterations, $\PSI^{(k)}(X)$ is indeed of the form
\[ \PSI^{(k)}(X)\vert_2 = \mathbf{vec}(\widehat{\theta}_{k-1} \otimes \psi_k)(X), \]
where $(\widehat{\theta}_{k-1})^\top = (\PSI_{:k-1}\vert_{k-1})^\top \, \widehat{\mathbf{U}}_{X, r}^{:k-1}\vert_{k-1}$ is a basis of $\widehat{\mathbb{G}}^{:k-1}$. It follows that the leading $r_k$ left singular vectors $U$ of $\PSI^{(k)}(X)\vert_2$ (see line \ref{algline:svd_step} of Algorithm~\ref{alg:SVD}) are also eigenvectors of the empirical Gramian $\widehat{C}(\mathbf{vec}(\widehat{\theta}_{k-1} \otimes \psi_k))$, as
\[ \frac{1}{m} \PSI^{(k)}(X)\vert_2 (\PSI^{(k)}(X)\vert_2)^\top = \widehat{C}(\mathbf{vec}(\widehat{\theta}_{k-1} \otimes \psi_k)). \]
Hence, the basis set
\[ (\widehat{\theta}_k)^\top = (\mathbf{vec}(\widehat{\theta}_{k-1} \otimes \psi_k))^\top U = (\PSI_{:k}\vert_k)^\top \widehat{\mathbf{U}}^{:k}_{X, r}\vert_k \]
is, by definition, a basis of $\widehat{\mathbb{G}}^{:k}$. It remains to show that the updated core $\PSI^{(k+1)}(X)$ is of the required form. By inspecting the update formula in line \ref{algline:update_svd}, we see that
\begin{align*}
(\Sigma V^\top) \PSI^{(k+1)}(X)\vert_1 &= (U^\top \PSI^{(k)}(X)\vert_2) \PSI^{(k+1)}(X)\vert_1 \\
&= (U^\top (\mathbf{vec}(\widehat{\theta}_{k-1} \otimes \psi_k)(X))) \PSI^{(k+1)}(X)\vert_1 \\
&=\widehat{\theta}_k(X) \PSI^{(k+1)}(X)\vert_1.
\end{align*}
The last expression is a matrix of shape $r_k \times (n_{k+1} m)$. Using the basis decomposition \eqref{eq: transformed data tensor} for $\PSI^{(k+1)}(X)$, we determine its entries as
\begin{align*}
\left(\widehat{\theta}_k(X) \PSI^{(k+1)}(X)\vert_1\right)_{l_k; i_{k+1}, l} &= \sum_{s = 1}^m \widehat{\theta}_{k, l_k}(x_s) \psi_{k+1, i_{k+1}}(x_s) \delta_{s, l} = \widehat{\theta}_{k, l_k}(x_l) \psi_{k+1, i_{k+1}}(x_l),
\end{align*}
which, upon re-shaping, equals $\mathbf{vec}(\widehat{\theta}_k \otimes \psi_{k+1})(X)$, as claimed. Finally, the same arguments also show that $(\widehat{\theta}_p)^\top = (\PSI\vert_p)^\top \widehat{\mathbf{U}}_{X, r}\vert_p$ is a basis of $\widehat{\mathbb{G}}^{:p}$, so $\widehat{\mathbf{U}}_{X, r}\vert_p$ is a basis of $\widehat{\mathbb{B}}^{:p}$, and the time series $\widehat{\theta}_p(X)$ is simply contained in the remaining components $\widehat{\Sigma}_{X, r} \widehat{V}_{X, r}^\top$ of the global SVD.
\end{proof}

We can now put the pieces together. In AMUSEt, with fixed ranks $\mathbf{r}$, we separately apply global SVD to $\PSI(X)$ and $\PHI(Y)$. We observe that, for $r = r_p$,
\begin{align*}
(\widehat{\eta}_{r})^\top &= \sqrt{m} (\PSI\vert_p)^\top \widehat{\mathbf{U}}_{X, r}\vert_p \widehat{\Sigma}_{X, r}^{-1}, & (\widehat{\zeta}_{r})^\top &= \sqrt{m} (\PHI\vert_p)^\top \widehat{\mathbf{U}}_{Y, r}\vert_p \widehat{\Sigma}_{Y, r}^{-1}
\end{align*}
are empirically orthonormal bases of the multi-linear spectral subspaces $\widehat{\mathbb{G}}^{:p}(\PSI)$ and $\widehat{\mathbb{G}}^{:p}(\PHI)$, respectively. Hence, the empirical Koopman operator between these spaces possesses the matrix representation
\begin{align*}
\widehat{K}_\tau(\widehat{\eta}_{r}, \widehat{\zeta}_{r}) &= \widehat{A}(\widehat{\eta}_{r}, \widehat{\zeta}_{r}) = \frac{1}{m}\widehat{\eta}_{r}(X) \widehat{\zeta}_{r}(Y)^\top \\
&= \widehat{\Sigma}_{X, r}^{-1} (\widehat{\mathbf{U}}_{X, r}\vert_p)^\top (\PSI(X)\vert_p) (\PHI(Y)\vert_p)^\top \widehat{\mathbf{U}}_{Y, r}\vert_p \widehat{\Sigma}_{Y, r}^{-1} = \widehat{V}_{X, r}^\top \widehat{V}_{Y, r} = \widehat{M}_{\tau, r},
\end{align*}
which is just the reduced matrix described in Section~\ref{subsec:amuset}. In summary, application of AMUSEt as outlined in Section~\ref{subsec:amuset} at fixed ranks provides a representation of the Koopman operator on the empirical multi-linear spectral subspaces $\widehat{\mathbb{G}}^{:p}(\PSI), \,\widehat{\mathbb{G}}^{:p}(\PHI)$, which consistently approximates the corresponding representation on $\mathbb{G}^{:p}(\PSI), \,\mathbb{G}^{:p}(\PHI)$, by Theorem~\ref{thm:consistency_koopman_ml_spectral}.

\section{Numerical Examples}\label{sec:Numerical Examples}

In this section, we provide numerical illustrations of the algorithmic and theoretical results presented in this study. All of the example systems have been analyzed using standard techniques before, which serve as reference results for our experiments. However, we will show that by means of our algorithms we are able to compute approximations to evolution operators on large trial space which are not amenable to standard treatment. The construction of these trial spaces requires a varying degree of preprocessing between examples. For the molecular examples, we still cannot do without prior knowledge, but the preprocessing pipeline is conceptually much simpler compared to standard methods.   

Our algorithms have been implemented in Python 3.6 and collected in the toolbox Scikit-TT\footnote{\url{https://github.com/PGelss/scikit_tt}}. Furthermore, we used d3s\footnote{\url{https://github.com/sklus/d3s}}, PyEMMA\footnote{\url{http://www.emma-project.org}} \cite{SCHERER2015} as well as scikit-image\footnote{\url{http://www.scikit-image.org/}} for simulating and analyzing the numerical examples.

\subsection{Molecular Dynamics}

We re-analyze two data sets of equilibrium molecular dynamics simulations in explicit water. The first system is the ten residue peptide deca-alanine (see~\cite{NUESKE2016} for the simulation setup). After downsampling, this data set comprises $m = 3\cdot 10^5$ frames at time spacing of $10\,\mathrm{ps}$. As a reference, we built a \emph{Markov state model} (MSM) using 500 discrete states by following a typical protocol from the literature (linear dimension reduction by TICA, followed by k-means clustering in reduced space), see \cite{BOWMAN2014}. The second system is 39 residue protein NTL9. The data were produced by D.~E.~Shaw~Research on the Anton Supercomputer \cite{LINDORFF2011}. The downsampled data set comprises approximately $56,000$ frames at a time spacing of $50 \, \mathrm{ns}$. The MSM analysis presented in \cite{BONINSEGNA2017} serves as reference model.

As both systems are stationary and reversible, we have $\rho_0 = \rho_1 = \mu$, and we can use identical trial spaces $\mathbb{V} = \mathbb{W}$. Since the Koopman operator is self-adjoint in this case, all evolution operators essentially contain the same information, and singular pairs of $\mathcal{K}_\tau$ are in fact eigenpairs. Following standard methodology in the field, we convert eigenvalue estimates $\widehat{\lambda}_i$ into \emph{implied timescales (ITS)} by the formula
\begin{equation}
\label{eq:def_implied_timescale}
  \widehat{t}_i(\tau) = -\frac{\tau}{\log(\widehat{\lambda}_i)}.
\end{equation}
Each ITS bears a unit of time and corresponds to the relaxation timescale of the dynamical process associated to eigenvalue $\lambda_i$. ITS are typically compared across a range of different lag times $\tau$ (implied timescale test), as ITS estimates are known to improve with increasing $\tau$, and observing a plateau indicates convergence of the Koopman model \cite{Sarich2010,Prinz2011,BOWMAN2014}. Consequently, we will also use the implied timescale test to evaluate the performance of tensor-based Koopman models.

\subsubsection{Deca-alanine}
\begin{figure}
  \centering
  \begin{subfigure}{.5\textwidth}
    \centering
    \hspace*{0.35cm}\includegraphics[height=4.75cm]{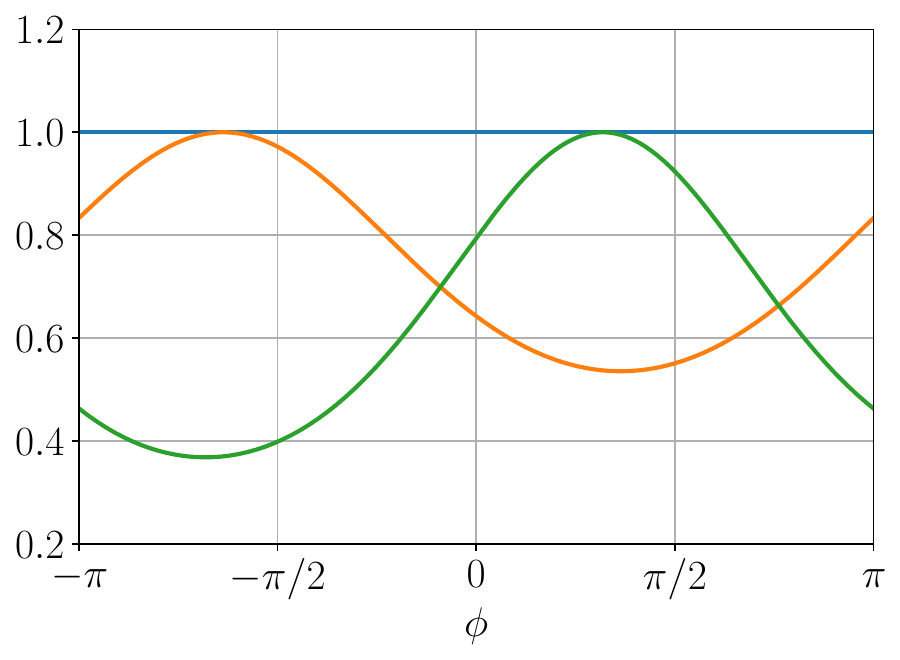}
    \caption{\hspace*{-0.8cm}}
  \end{subfigure}%
  \begin{subfigure}{.5\textwidth}
    \centering
    \hspace*{0.35cm}\includegraphics[height=4.75cm]{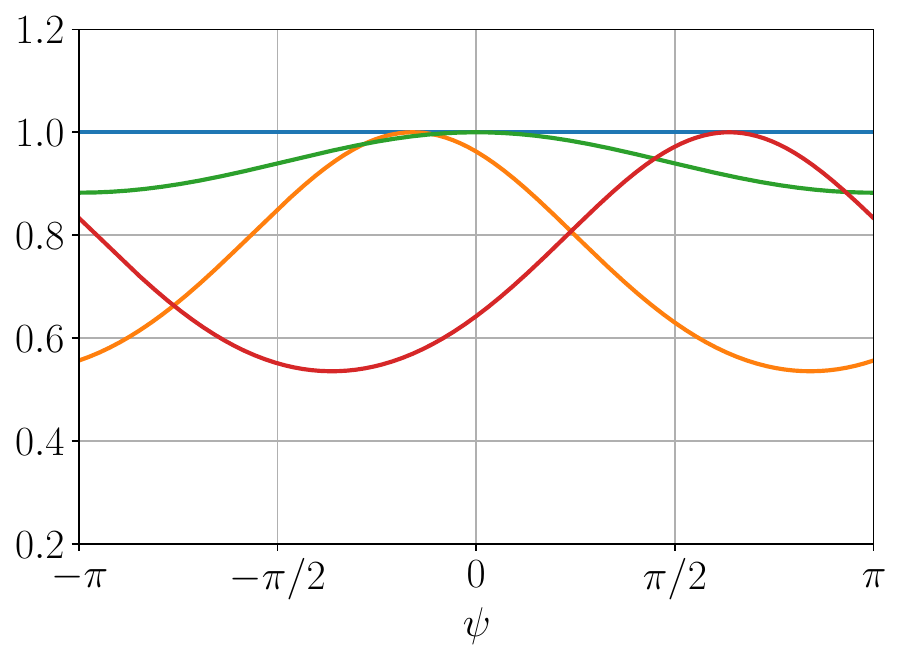}
    \caption{\hspace*{-0.8cm}}
  \end{subfigure}\\[0.3cm]
  \begin{subfigure}{.5\textwidth}
    \centering
    \includegraphics[height=5.50cm]{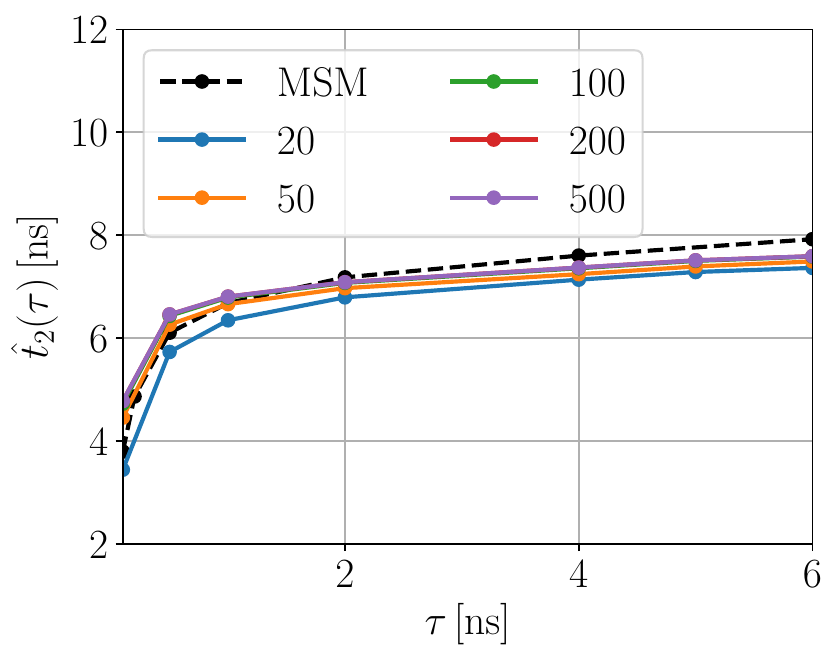}
    \caption{\hspace*{-0.8cm}}
  \end{subfigure}%
  \begin{subfigure}{.5\textwidth}
    \centering
    \includegraphics[height=5.50cm]{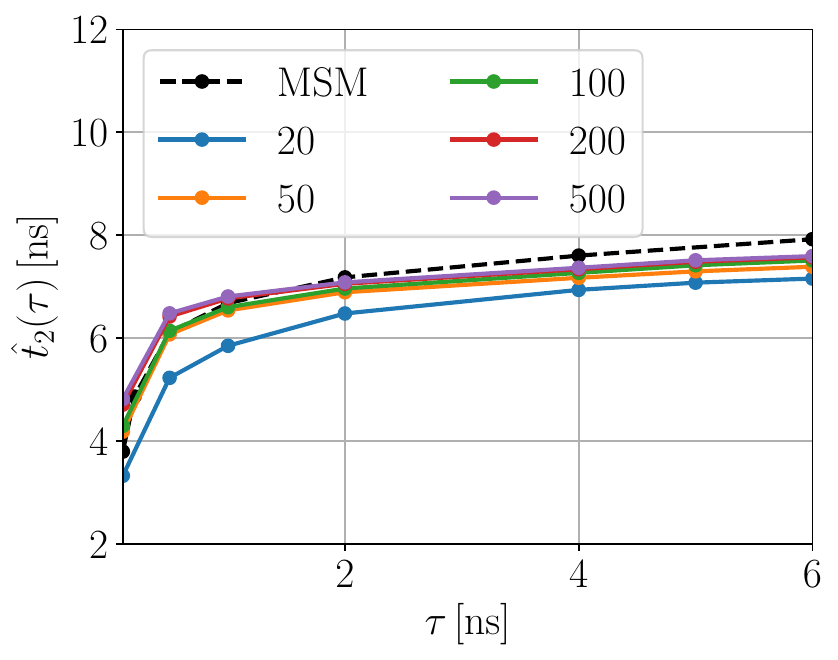}
    \caption{\hspace*{-0.8cm}}
  \end{subfigure}
  \caption{Results for molecular dynamics simulation data of deca-alanine peptide: (a) Univariate basis set used for all $\phi$-dihedral angles, comprised of the constant and periodic Gaussians centered at $c_{i_k} = \{-2, 1\}$, with $s_{i_k} = \{0.8, 0.5\}$. (b) The same for all $\psi$-dihedral angles, where periodic Gaussians are centered at $c_{i_k} = \{-0.5, 0.0, 2.0\}$, with $s_{i_k} = \{0.8, 4.0, 0.8\}$. (c) Slowest timescale $\widehat{t}_2$ obtained from \eqref{eq:reduced_ev_problem} after constructing $\PSI(Z)$ using the TT decomposition \eqref{eq: transformed data tensor}. We show results for different values of the maximal rank allowed during the global SVD, and of the lag time~$\tau$. The reference MSM is represented by the black line. (d) The same if $\PSI(Z)$ is represented by the HOCUR algorithm, for different values of the maximal rank in Algorithm~\ref{alg: HOCUR}.}
  \label{fig: ala10_results}
\end{figure}

For deca-alanine, we construct the trial space $\mathbb{V}$ by first choosing $p = 10$ backbone dihedral angles of the peptide as elementary descriptors, and then defining a subspace $\mathbb{V}^k$ for each of them as the span of either $n_k = 3$ or $n_k = 4$ scalar functions on that dihedral angle. These functions always include the constant and two or three periodic Gaussians of the form
\begin{equation*}
\psi_{k, i_k}(x) = \exp\left[-\frac{1}{2s_{i_k}}\sin^2(0.5(x - c_{i_k}))\right].
\end{equation*}
Their positions and shapes are chosen to align with the typical marginal distribution of protein data along its backbone dihedral angles, see Figures~\ref{fig: ala10_results}~(a) and (b). The full tensor space $\mathbb{V}$ is then of dimension $N = 3^5 \cdot 4^5 \approx 2.5 \cdot 10^5$, which already exceeds what is typically considered feasible as trial space dimension for standard methods. The construction for $\mathbb{V}^k$ used here is clearly informed by physical insight, but only at a fairly basic level.

As the system is in equilibrium, we can generate the data matrices $X,\, Y$ from just a single long trajectory, with $Y$ being obtained by shifting all time steps in $X$ by the lag time $\tau$. We then use the decomposition \eqref{eq: transformed data tensor} and the HOCUR decomposition (Section~\ref{sec: HOCUR for TDT}) to arrive at TT representations of the transformed data tensor $\PSI(Z)$. Subsequently, the procedure outlined in Section~\ref{subsec:amuset} is applied to obtain of the reduced problem \eqref{eq:reduced_ev_problem}. Focusing on just the slowest dynamical process, we monitor the second ITS $\widehat{t}_2(\tau)$ given by \eqref{eq:def_implied_timescale}, as a function of $\tau$.

In line with the theoretical results presented in Section~\ref{sec:analysis_amuset}, we build the direct representation \eqref{eq: transformed data tensor} based on varying amounts of data, ranging between $m = 3000$ and $m = 3\cdot 10^5$ data points. However, we find the resulting timescale estimates to be virtually indistinguishable. The quality of approximation seems to depend more critically on the rank of the TT representation for $\PSI(Z)$. We analyze this dependence by capping the maximal rank $r$ allowed either for the global SVD of \eqref{eq: transformed data tensor}, or for the HOCUR representation of $\PSI(Z)$, at different values. The resulting timescale estimates, as a function of the lag time $\tau$, are shown in Figures~\ref{fig: ala10_results}(c)-(d). In both cases, a maximal rank of $50$ is sufficient to obtain excellent agreement with the Markov model results.

\subsubsection{NTL9}
For NTL9, we also follow an established protocol to arrive at a basic set of descriptors. We consider all closest heavy-atom distances between protein residues, and rank these distances by the fraction of simulation time during which a contact between residues was formed (i.e. their distance is smaller than $0.35 \, \mathrm{nm}$). For each of these distance features, we define a space $\mathbb{V}^k$ as the span of the constant and two Gaussian functions, given by
\begin{equation*}
\psi_{k, i_k}(x) = \exp\left[-\frac{1}{2s_{i_k}}(x - c_{i_k})^2\right],
\end{equation*}
centered at $c_{i_k} \in \{0.285, 0.62\}$, with $s_{i_k} \in \{0.001, 0.01\}$. Again, these parameters were selected to make sure that the marginal distribution of the data along each distance can be reproduced by a linear combination of the Gaussians. Just as in the previous example, the construction of the elementary function spaces relies on some degree of preprocessing, but only at a basic level. Below, we use either the first $p = 10$ or $p = 20$ distance features to construct the full tensor space, which is therefore of dimension $N = 3^{10} \approx 6 \cdot 10^4$ or $N = 3^{20} \approx 3.5 \cdot 10^9$. Both dimensions are beyond what is considered tractable for standard methods.

\begin{figure}
  \centering
  \begin{subfigure}{.5\textwidth}
    \centering
    \includegraphics[height=5.50cm]{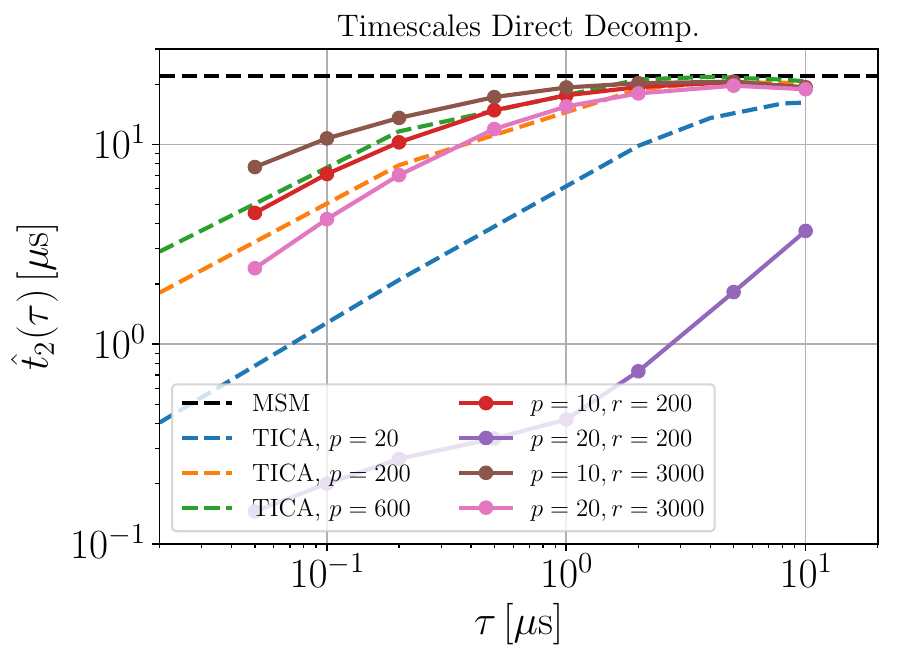}
    \caption{\hspace*{-0.4cm}}
  \end{subfigure}%
  \begin{subfigure}{.5\textwidth}
    \centering
    \vspace*{0.09cm}
    \includegraphics[height=5.50cm]{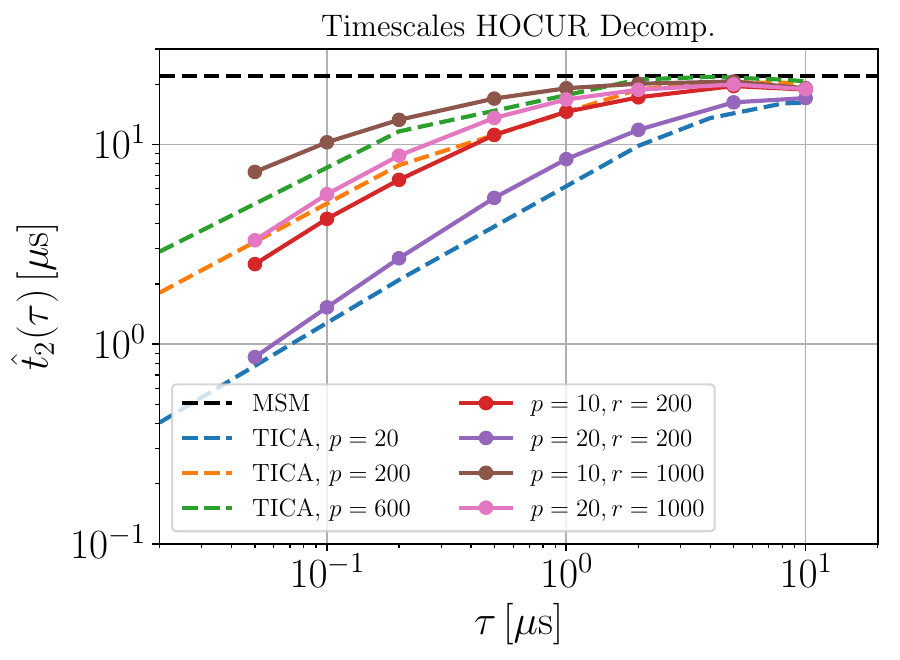}
    \caption{\hspace*{-1cm}}
  \end{subfigure}
  \caption{Results for molecular dynamics simulation data of NTL9 protein: (a) Slowest timescale obtained from \eqref{eq:reduced_ev_problem} after representing the data tensor $\PSI(Z)$ using the exact decomposition \eqref{eq: transformed data tensor}. We show results for different maximal ranks during the global SVD of $\PSI(Z)$, and if either the first $p =10$ or $p = 20$ distance features are used. The MSM-based reference value, as well as the timescales computed by linear TICA are indicated by the dashed lines. (b) The same if $\PSI(Z)$ is represented by the HOCUR decomposition, with different maximal ranks in Algorithm~\ref{alg: HOCUR}.}
  \label{fig: ntl9 1}
\end{figure}

We follow essentially the same protocol as for deca-alanine, using both the direct decomposition~\eqref{eq: transformed data tensor} and the HOCUR iteration to represent the data tensor $\PSI(Z)$. For a series of lag times $\tau$ and various maximal ranks $r$, we apply the procedure outlined in Section~\ref{subsec:amuset} to obtain the reduced matrix $\widehat{M}_{\tau, r}$ in \eqref{eq:reduced_ev_problem}. Estimates for the slowest implied timescale $\widehat{t}_2(\tau)$ are shown in panel (a) of Figure~\ref{fig: ntl9 1} for the decomposition \eqref{eq: transformed data tensor}, and for the HOCUR representation in panel (b) of the same figure.

Two observations stand out: the first is, since our focus is to illustrate the advantage of using a large non-linear model class, a natural comparison for our model is the direct application of linear TICA~\cite{PEREZ2013}, rather than the MSM analysis from \cite{BONINSEGNA2017}, which required a significant level of expertise. Linear TICA is the same as applying the standard AMUSE Algorithm \ref{alg:AMUSE_SINGLE} to the basis set given by all elementary descriptors (i.e. the identity function on each distance feature in our case). Our results show that a tensor-based model on about ten distances, using a moderate TT rank, provides the same performance as linear TICA on several hundreds of distance features. The second observation is that including additional distance features ($p = 20$) does not improve timescale estimates. Rather, larger ranks are required to ensure the same quality of approximation. The HOCUR-based models seem to be less sensitive to this effect for this example. An intuitive explanation would be that larger ranks are required to pass the relevant information down the chain of TT cores. We conclude that efficient and stable ways to compute low-rank representations of the data tensor $\PSI(X)$ clearly remain an important topic for future work.

\begin{figure}
  \centering
  \begin{subfigure}{.5\textwidth}
    \centering
    \includegraphics[height=4.75cm]{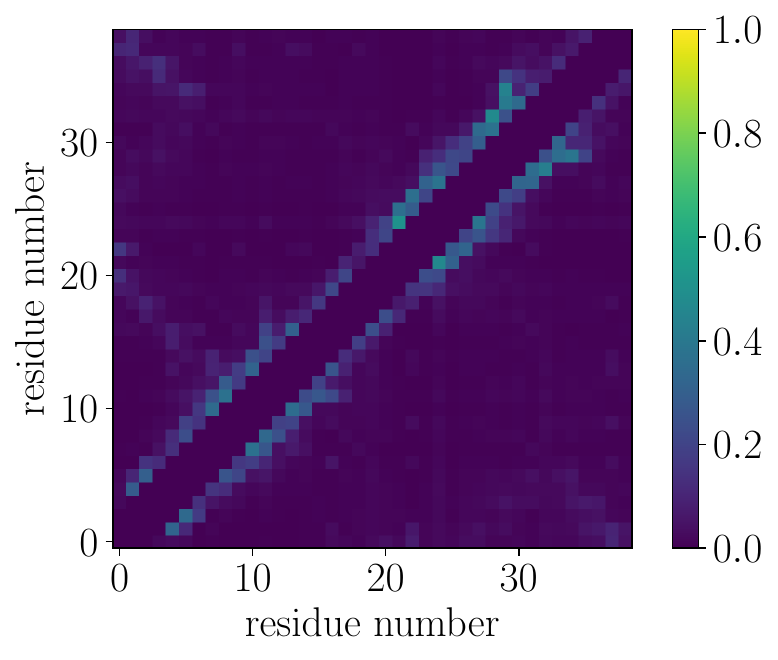}
    \caption{\hspace*{0.1cm}}
  \end{subfigure}%
  \begin{subfigure}{.5\textwidth}
    \centering
    \includegraphics[height=4.75cm]{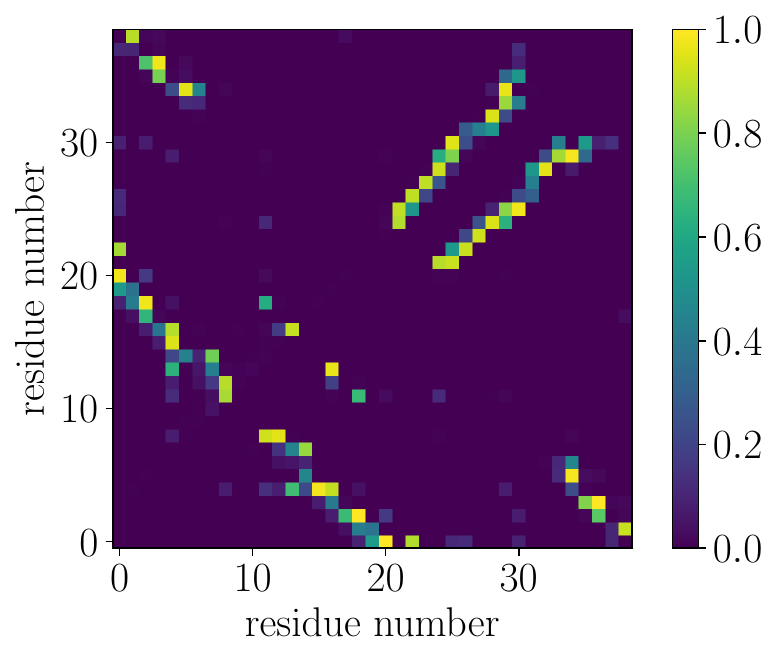}
    \caption{\hspace*{0.1cm}}
  \end{subfigure}\\[0.3cm]
  \begin{subfigure}{.5\textwidth}
    \centering
    \includegraphics[height=4cm]{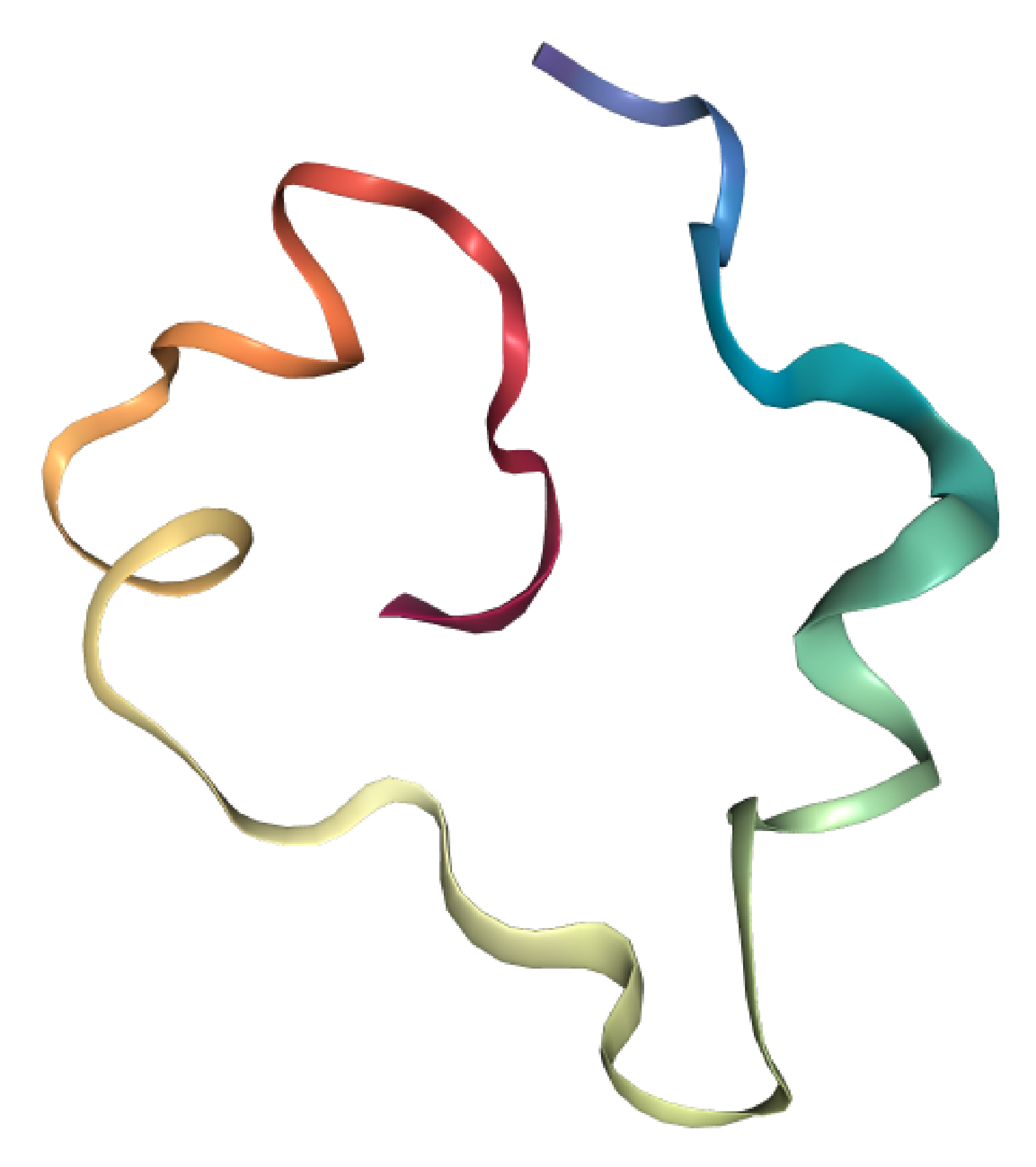}
    \caption{\hspace*{0.1cm}}
  \end{subfigure}%
  \begin{subfigure}{.5\textwidth}
    \centering
    \includegraphics[height=4cm]{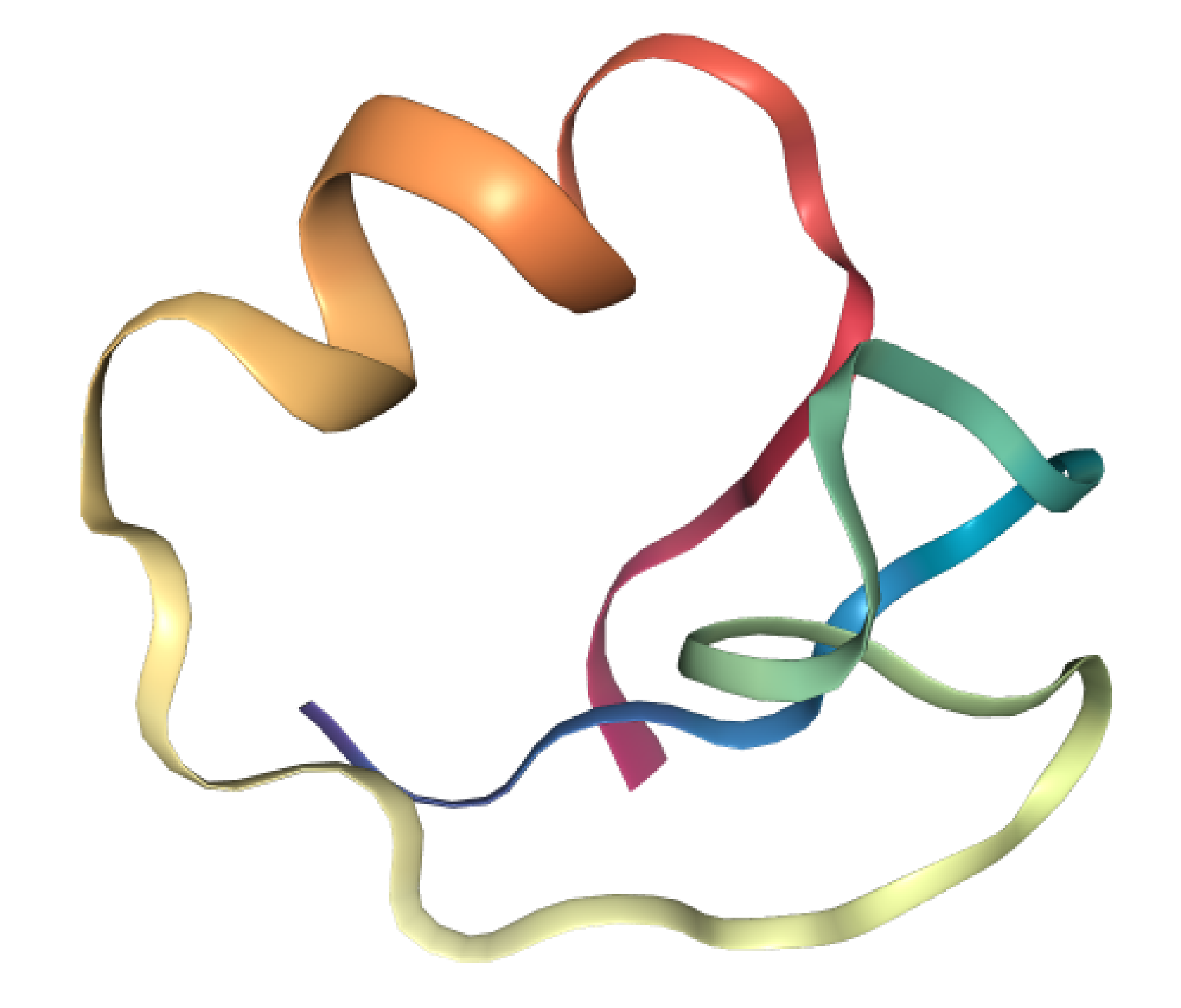}
    \caption{\hspace*{0.1cm}}
  \end{subfigure}
  \caption{States of NTL9: (a) Contact map for the unfolded state computed from the first two eigenfunctions of the HOCUR model corresponding to $p = 10$, $r = 200$, and $\tau = 2\, \mathrm{\mu s}$ (upper left triangle), compared to the corresponding contact map of the reference MSM (lower right triangle). (b) The same for the folded state. (c)/(d) Representative molecular structures for the unfolded and folded state of NTL9.}
  \label{fig: ntl9 2}
\end{figure}

For completeness, we also verify that the eigenfunctions estimated by AMUSEt correctly encode the folding process of NTL9, which is known as the slowest dynamical process for this system. To this end, we apply PCCA~\cite{DEUFLHARD2005} to the time series of the first two eigenfunctions and assign each snapshot to one of two metastable states if the degree of membership exceeds $0.5$. We then calculate the contact frequencies separately for each of the two states. The resulting so-called contact maps are shown in the upper left triangles of Figures~\ref{fig: ntl9 2}~(a) and (b). By comparing to the contact maps provided by the reference Markov model (lower right triangles in Figures~\ref{fig: ntl9 2}~(a) and (b)), we see that there is virtually no difference.

\subsection{ABC Flow}

\begin{figure}[htbp]
    \centering
    \includegraphics[height=7cm]{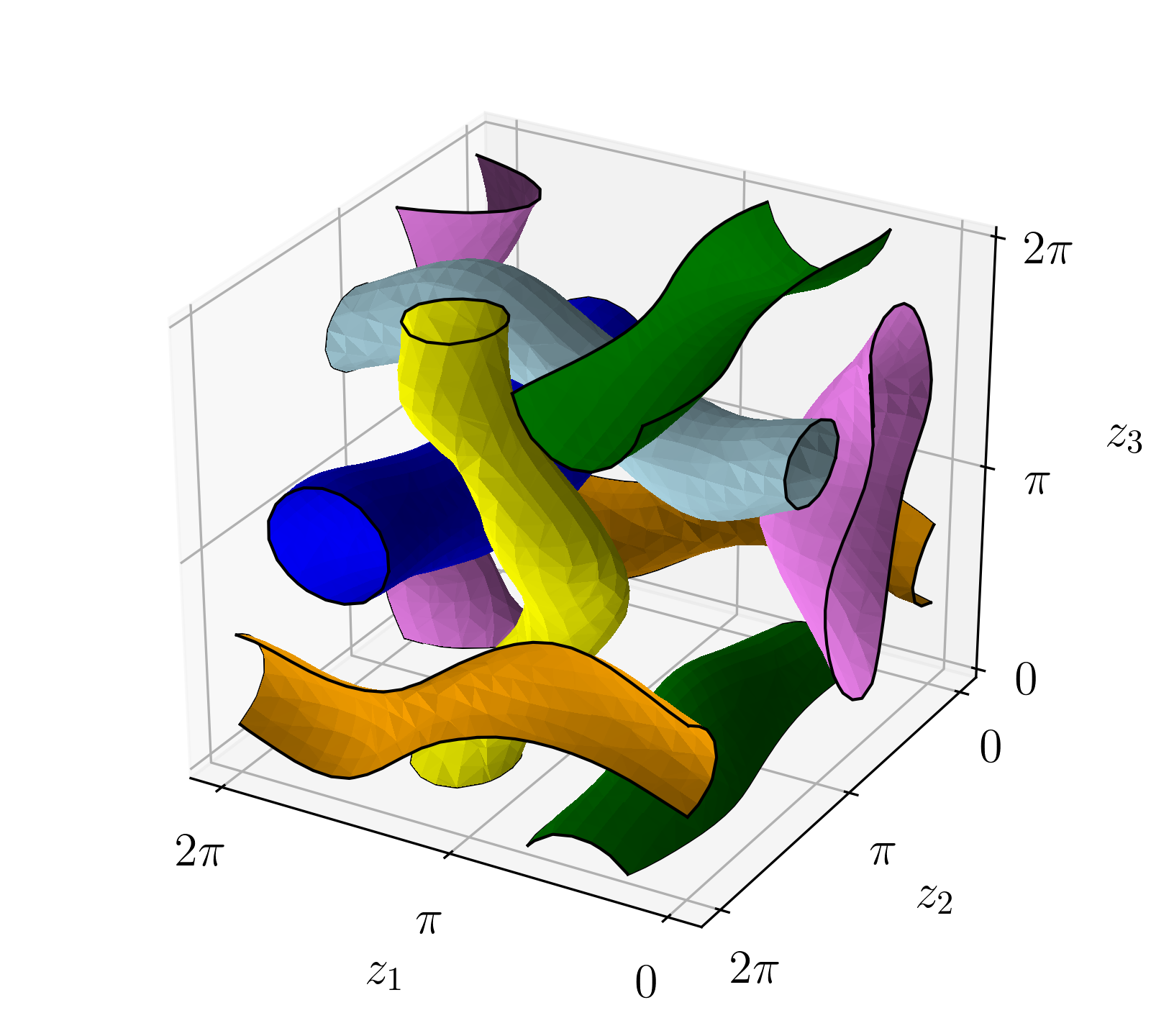}
    \caption{Results for the non-stationary ABC flow: six coherent vortices in the domain $[0, 2\pi]^3$, obtained from approximate eigenfunctions of the forward-backward operator $\widehat{\mathcal{F}}_\tau$, computed by AMUSEt with $10^3$ Gaussian basis functions.}
    \label{fig: ABC}
\end{figure}

Finally, we study a popular toy model for fluid dynamics problems, in order to illustrate the treatment of the non-stationary case and the computation of coherent sets using tensor-based methods. Let us consider the well-known ABC (Arnold--Beltrami--Childress) flow, given by the ordinary differential equation
\begin{align*}
    \dot{z}_1 = A \sin(z_3) + C \cos(z_2), \\
    \dot{z}_2 = B \sin(z_1) + A \cos(z_3), \\
    \dot{z}_3 = C \sin(z_2) + B \cos(z_1),
\end{align*}
with $ A = \sqrt{3} $, $ B = \sqrt{2} $, and $ C = 1 $. The system is defined on the torus, i.e., $ 0 \le z_i \le 2 \pi $ for $ i = 1, 2, 3 $, see~\cite{FP09} for details. We define the lag time to be $ \tau = 5 $ and sample $ 25^3 $ test points $ x_i $ uniformly in $ [0, 2 \pi]^3 $. In order to compute the corresponding points $ y_i $, we use a standard Runge--Kutta integrator with variable step size. Using a coordinate-major decomposition comprising ten Gaussian functions with variance $1$ in each dimension, we arrive at a tensor space $\mathbb{V} = \mathbb{W}$ spanned by $1000$ three-dimensional Gaussians on an equidistant grid. Since $\rho_0 \neq \rho_1$, we are interested in eigenpairs of the projected forward-backward operator $\widehat{\mathcal{F}}_\tau(\mathbb{V, \mathbb{V}})$. We apply AMUSEt as described in Section~\ref{subsec:amuset}, and we use scikit-image in order to extract isosurfaces from the computed eigenfunctions. The coherent sets shown in Figure~\ref{fig: ABC} are consistent with the results presented in \cite{BK17:coherent}.

\section{Summary}
\label{sec:Conclusion and Outlook}

We have presented novel techniques to approximate evolution operators associated with high-dimensional dynamical systems using tensor-structured basis sets. Specifically, we have introduced AMUSEt, a multi-linear version of the AMUSE algorithm, which allows us to derive a reduced matrix representation of evolution operators while operating only on the data tensor. For fixed multi-linear ranks, we have established convergence of AMUSEt in the limit of infinite data. In addition, we have provided a detailed algorithmic description of a novel iterative method to compute a higher-order CUR decomposition of the data tensor, which only requires evaluations of the basis set on the data. We have also presented successful applications to benchmarking data sets of molecular dynamics simulation and fluid dynamics.

\section*{Acknowledgements}

This research has been funded by the National Science Foundation (CHE-1265929, CHE-1738990, CHE-1900374, PHY-1427654) [FN, CC], the Welch Foundation (C-1570) [FN, CC], the Rice University Academy of Fellows [FN], Deutsche Forschungsgemeinschaft (CRC 1114, \emph{``Scaling Cascades in Complex Systems''}) [PG, SK], and the Einstein Foundation Berlin [CC]. The authors are grateful to D. E. Shaw Research for providing the NTL9 simulation data, and to the Paderborn Center for Parallel Computing for computational resources.

{\small{}\bibliographystyle{unsrturl}
\bibliography{references}
}{\small\par}

\appendix

\section{Auxiliary Results}
\label{app:aux_results}

First, we compile some useful results concerning linear operations in finite-dimensional Hilbert spaces:
\begin{Lemma}
\label{lem:equivalences_subspaces}
Let $\mathbb{V}$ be a finite-dimensional Hilbert space with basis $\psi = [\psi_1, \ldots, \psi_n]^\top$ and Gramian matrix $C(\psi)$.
\begin{itemize}
 \item[(i)] Let $\theta_1, \theta_2 \in \mathbb{V}$ with coefficient vectors $a_1, \, a_2$ with respect to $\psi$. Also, let $\mathcal{B}: \mathbb{V} \mapsto \mathbb{V}$ be a linear operator on $\mathbb{V}$ with matrix representation $B$ with respect to $\psi$. Then
\begin{align*}
\innerprod{\theta_1}{\theta_2}_\mathbb{V} &= a_1^\top C(\psi) a_2, & \|\mathcal{B}\|_{L(\mathbb{V})} &= \|C(\psi)^{1/2} B C(\psi)^{-1/2} \|_2.
\end{align*}
 \item[(ii)] Let $\mathbb{F} \subset \mathbb{V}$ be an $r$-dimensional subspace with basis $\theta$. Let the matrix of coefficient vectors associated to $\theta$ be $A \in \mathbb{R}^{n \times r}$. The matrix representation of the projector $\mathcal{P}_\mathbb{F}$ with respect to $\psi$ is
\[P_\mathbb{F} = A (A^\top C(\psi) A)^{-1} A^\top C(\psi). \]
 \item[(iii)] Let $\mathbb{F}_\nu, \, \nu = 1, 2, \ldots$ and $\mathbb{F}$ be subspaces of $\mathbb{V}$ of dimension $r \leq n$, with coefficient vector spaces $\mathbb{A}_\nu, \, \mathbb{A} \subset \mathbb{R}^n$. Then
\begin{align*}
d(\mathbb{A}_\nu, \mathbb{A}) &\rightarrow 0 ~~\Rightarrow~~ d(\mathbb{F}_\nu, \mathbb{F}) \rightarrow 0.
\end{align*}
\end{itemize}
\end{Lemma}
\begin{proof}
(i) and (ii) can be verified directly (see also \cite{KNYAZEV2002}). To prove (iii), we fix an orthonormal basis $U \in \mathbb{R}^{n \times r}$ of $\mathbb{A}$, and then use auxiliary Lemma \ref{lem:convergence_procrustes} below to select ONBs $U_\nu$ of $\mathbb{A}_\nu$ such that $\| U_\nu - U \|_2 \rightarrow 0$ is also satisfied. The orthonormal projectors $\mathcal{P}_{\mathbb{F}_\nu}, \, \mathcal{P}_{\mathbb{F}}$ then have matrix representations with respect to $\psi$ given by
\begin{align*}
P_{\mathbb{F}_\nu} &= U_\nu (U_\nu^\top C(\psi) U_\nu)^{-1} U_\nu^\top C(\psi), & P_{\mathbb{F}} &= U (U^\top C(\psi) U)^{-1} U^\top C(\psi).
\end{align*}
Using part (i), we then find
\begin{align*}
\|\mathcal{P}_{\mathbb{F}_\nu} - \mathcal{P}_{\mathbb{F}} \|_{L(\mathbb{V})} &= \|C(\psi)^{1/2} \left(U_\nu (U_\nu^\top C(\psi) U_\nu)^{-1} U_\nu^\top -  U (U^\top C(\psi) U)^{-1} U^\top \right) C(\psi)^{1/2} \|_2 \\
&\leq \|C(\psi)\|_2 \| U_\nu (U_\nu^\top C(\psi) U_\nu)^{-1} U_\nu^\top -  U (U^\top C(\psi) U)^{-1} U^\top \|_2 \rightarrow 0. \qedhere
\end{align*}
\end{proof}

The following technical result helps us translate convergence of subspaces in Euclidean space into convergence of specific orthonormal bases:
\begin{Lemma}
\label{lem:convergence_procrustes}
For $\nu \in \mathbb{N}$, let $\mathbb{A}_\nu \subset \mathbb{R}^n$ be a subspace of dimension $r \leq n$, each with an orthonormal basis $U_\nu \in \mathbb{R}^{n \times r}$. Let $\mathbb{A} \subset \mathbb{R}^n$ be another subspace of the same dimension, with orthonormal basis $U \in \mathbb{R}^{n \times r}$, and assume $d(\mathbb{A}_\nu, \mathbb{A}) \rightarrow 0$ as $\nu \rightarrow \infty$. For each $\nu$, define $\bar{U}_\nu \in \mathbb{R}^{n \times r}$ as the closest orthonormal matrix to $U$ in the column span of $U_\nu$, i.e., $\bar{U}_\nu$ is obtained by solving the orthogonal Procrustes problem
\begin{align*}
\bar{U}_\nu &= U_\nu R_\nu, & R_\nu = \argmin_{\substack{R\in \mathbb{R}^{r \times r},\\ R^\top\!R = \mathrm{Id}}} \, \|U_\nu R - U \|^2_F.
\end{align*}
Then we have $\|U - \bar{U}_\nu\|_2 \rightarrow 0$.
\end{Lemma}

\begin{proof}
From an SVD of $U_\nu^\top U$, i.e. $U_\nu^\top U  = V_\nu \Sigma_\nu W_\nu^\top$, the optimal orthogonal transformation can be obtained as $R_\nu = V_\nu W_\nu^\top$ \cite{SCHOENEMANN1966}. Moreover, from the above SVD we also obtain a compact singular value decomposition of $U_\nu U_\nu^\top UU^\top$ via
\begin{align*}
U_\nu U_\nu^\top UU^\top = (U_\nu V_\nu) \Sigma_\nu (W_\nu^\top U^\top),
\end{align*}
as the matrices in brackets are orthonormal and the rank of the left-hand side equals $r$. Now, by the convergence of the subspaces $\mathbb{A}_\nu$ towards $\mathbb{A}$, we must have $U_\nu U_\nu^\top UU^\top \rightarrow UU^\top$, and as the latter is an orthogonal projector, we conclude $\Sigma_\nu \rightarrow \mathrm{Id}_{r}$ with $\nu \rightarrow \infty$. Now consider the difference between $U$ and $\bar{U}_\nu$:
\begin{align*}
U - \bar{U}_\nu = U - U_\nu V_\nu W_\nu^\top = U - U_\nu V_\nu \Sigma_\nu^{-1} V_\nu^\top U_\nu^\top U = (\mathrm{Id} - U_\nu V_\nu \Sigma_\nu^{-1} V_\nu^\top U_\nu^\top) U.
\end{align*}
For each $\nu$, the columns of $U_\nu V_\nu \in \mathbb{R}^{n \times r}$ are orthonormal. Let the columns of $Y_\nu \in \mathbb{R}^{n \times n-r}$ form an orthonormal basis of the complement of $\mathrm{span}(U_\nu V_\nu)$. Denote the orthonormal $n \times n$-matrices obtained by padding these two column sets by $S_n = \left[ U_\nu V_\nu \mid Y_\nu\right]$. With this, we can further manipulate the above expression as follows:
\begin{align*}
U - \bar{U}_\nu &= \left[ S_\nu S_\nu^\top -  U_\nu V_\nu \Sigma_\nu^{-1} V_\nu^\top U_\nu^\top\right] U \\
&= \left[ \left[ U_\nu V_\nu \mid Y_\nu\right] \begin{bmatrix}(V_\nu^\top U_\nu^\top) \\ Y_\nu^\top \end{bmatrix} - U_\nu V_\nu \Sigma_\nu^{-1} V_\nu^\top U_\nu^\top\right] U \\
&= (U_\nu V_\nu) (\mathrm{Id} - \Sigma_\nu^{-1} ) (V_\nu^\top U_\nu^\top) U + Y_\nu Y_\nu^\top U.
\end{align*}
As $Y_\nu Y_\nu^\top$ is a matrix representation of the orthogonal projector onto $\mathbb{A}_\nu^\perp$, we must have $Y_\nu Y_\nu^\top U \rightarrow 0$ by the convergence of the subspaces $\mathbb{A}_\nu$. For the first term, we find:
\begin{align*}
\| (U_\nu V_\nu) (\mathrm{Id} - \Sigma_\nu^{-1} ) (V_\nu^\top U_\nu^\top) U \|_2 &\leq \| (U_\nu V_\nu) \|_2 \| (\mathrm{Id} - \Sigma_\nu^{-1} )\|_2 \|(V_\nu^\top U_\nu^\top)\|_2  \|U \|_2 \\
&= \| (\mathrm{Id} - \Sigma_\nu^{-1} )\|_2 \rightarrow 0,
\end{align*}
by our observations above. This completes the proof.
\end{proof}

\end{document}